\documentclass{article}
\usepackage{amssymb, amsmath}
\usepackage[latin9]{inputenc}
\usepackage[francais]{babel}
\usepackage[affil-it]{authblk}
\usepackage{graphicx}
\usepackage{hyperref}
\usepackage{color}
\usepackage[onehalfspacing]{setspace}
\usepackage{ulem}
\makeatletter
\def\imod#1{\allowbreak\mkern10mu({\operator@font mod}\,\,#1)}
\makeatother

\usepackage[babel=true]{csquotes}
\usepackage[T1]{fontenc}
\usepackage{fancyhdr}

\begin{document}

\title{Autour de pratiques algébriques de Poincaré : \\ héritages de la réduction de Jordan}
\author{Frédéric Brechenmacher
\thanks{Electronic address: \texttt{frederic.brechenmacher@euler.univ-artois.fr \\ Ce travail a  bénéficié d'une aide de l'Agence Nationale de la Recherche : projet CaaFÉ (ANR-10-JCJC 0101)}}}
\affil{Université d'Artois \\ Laboratoire de mathématiques de Lens (EA 2462) \\
rue Jean Souvraz S.P. 18, 62307 Lens Cedex France.
 \\ \& \\
École polytechnique \\ Département humanités et sciences sociales \\
91128 Palaiseau Cedex, France. \\}
 \date{}
\maketitle

\begin{abstract}
Cet article propose un regard transversal sur les travaux menés par Poincaré entre 1878 et 1885. 
L'étude est amorcée à partir de la courte note "Sur les nombres complexes" publiée en 1884. Cette note permet de problématiser les cloisonnements disciplinaires rétrospectifs selon lesquels ont souvent été décrits les travaux de  Poincaré. La note a en effet été présentée comme le point d'origine de la reconnaissance du rôle joué par l'algèbre des matrices dans les théories des algèbres associatives et de Lie, ouvrant ainsi la voie aux développements de relations entre ces deux théories et à de premiers résultats sur leurs structures. Nous montrons cependant que cette note n'est pas sous-tendue par une notion ou une théorie mais par une pratique algébrique de classification de groupes linéaires par réduction canonique (dite de Jordan) de leurs substitutions. 

Cette pratique algébrique s'avère jouer un rôle transversal dans les travaux de Poincaré. Nous questionnons les innovations individuelles qu'elle manifeste  au regard de différentes échelles collectives. Nous analysons notamment l'appropriation par Poincaré des travaux de Jordan sur les groupes linéaires au prisme de l'héritage d'Hermite sur la réduction des formes algébriques. Cet héritage mêlé irrigue l'ensemble des travaux de Poincaré entre 1878 et 1885. Il est non seulement au c\oe{}ur du développement de la célèbre théorie des fonctions fuchsiennes mais aussi  de travaux sur la théorie algébrique des formes. Ces derniers n'ont été que très peu commentés jusqu'à présent. Ils ont pourtant joué un rôle clé dans l'attribution par Poincaré d'une place centrale à la notion de groupe. 

En amont des travaux de Poincaré, l'article étudie aussi les circonstances de l'énoncé du théorème de la réduction canonique par Jordan ainsi que la  réception de ce résultat. En aval, nous évoquons le rôle des travaux de Poincaré pour la transmission d'un héritage jordano-hermitien dans les procédés du calcul des Tableaux tels que les développent au début du XX\up{e} siècle des mathématiciens comme Autonne ou Châtelet.
 \bigskip
 
\textbf{Abstract. \\ The legacy of Jordan's canonical form on Poincaré's algebraic practices} \\ 
This paper proposes a transversal overview on Henri Poincaré's early works (1878-1885).  Our investigations start with a case study of a short note published by Poincaré on 1884 : "Sur les nombres complexes". In the perspective of today's mathematical disciplines - especially linear algebra  -, this note seems completely isolated in Poincaré's works. This short paper actually exemplifies that the categories used today for describing some collective organizations of knowledge fail to grasp both  the collective dimensions and individual specificity of Poincarés work. It also highlights the crucial and transversal role played in Poincaré's works by a specific algebraic practice of classification of linear groups by reducing the analytical representation of linear substitution to their Jordan's canonical forms.

We then analyze in details this algebraic practice as well as the roles it plays in Poincaré's works. We first provide a micro-historical analysis of Poincaré's appropriation of Jordan's approach to linear groups through the prism of the legacy of Hermite's works on algebraic forms between 1879 and 1881. This mixed legacy illuminates the interrelations between all the papers published by Poincaré between 1878 and 1885 ;  especially between some researches on algebraic forms and the development of the theory of Fuchsian functions.  Moreover, our investigation  sheds new light on how the notion of group came to play a key role in Poincaré's approach.

The present paper also offers a historical account of the statement by Jordan of his canonical form theorem. Further, we  analyze how Poincaré transformed  this theorem by appealing to Hermite's "calcul des Tableaux" as well as the major role  this reformulation played for the circulation of Jordan's theorem for  the next decades.
\end{abstract}

\tableofcontents

\pagestyle{fancy}
\fancyhead{}
\rhead{ F. BRECHENMACHER}

\newpage

\section*{Introduction}

Henri Poincaré est souvent présenté comme l'un des derniers mathématiciens universels. Ses travaux concernent en effet des domaines variés des sciences mathématiques.\footnote{Pour un regard d'ensemble sur les travaux de Poincaré et leur historiographie, voir \cite{Nabonnand2000}.} Cette variété thématique a souvent été mise au crédit des capacités individuelles exceptionnelles de Poincaré. Elle peut cependant également s'analyser comme témoignant d'organisations mathématiques antérieures aux découpages disciplinaires qui nous sont aujourd'hui familiers. De fait, le savant était contemporain de la montée en puissance de disciplines centrées sur des objets (groupes finis, groupes continus, algèbres associatives etc.) face aux "branches" en lesquelles se divisaient traditionnellement l'organisation des sciences mathématiques en France (Analyse, Géométrie, Mécanique, Physique mathématique etc.). 

Nous proposons dans cet article d'analyser les pratiques algébriques qui se présentent dans la première décennie de travaux de Poincaré. Ces pratiques sont mobilisées au sein de travaux variés, de l'arithmétique des formes quadratiques à la mécanique céleste en passant par les fonctions fuchsiennes. Plus encore, ces pratiques supportent souvent les liens qui sont établis entre ces différents domaines. Elles témoignent ainsi d'une organisation mathématique dans laquelle l'algèbre joue un rôle transversal à différents domaines mais d'une manière différente du rôle qu'ont jouées ultérieurement les structures algébriques abstraites. 

Ces pratiques algébriques ne sont jamais thématisées dans un cadre théorique autonome. S'il arrive fréquemment qu'elles jouent un rôle important dans une publication donnée, elles s'y articulent toujours à d'autres enjeux. La nature algébrique de ces pratiques se manifeste en fait précisément par leur capacité à s'articuler aux différents types de significations portées par des problèmes d'analyse, d'arithmétique, de mécanique etc. 

L'étude de ces pratiques nécessite donc d'être davantage attentif aux relations entre les travaux de Poincaré qu'à des textes isolés. Il nous faudra ainsi analyser globalement le corpus des textes publiés par Poincaré ainsi que leurs liens avec des travaux d'autres auteurs. 

Le terme pratique algébrique, tel que nous l'employons ici, désigne à la fois des procédés opératoires et des aspects culturels propres aux groupes de textes dans lesquels de tels procédés circulent.\footnote{Le sens pris ici par le terme pratique algébrique ne correspond donc pas à aux usages qui peuvent être faits du terme "pratique mathématique" pour distinguer les pratiques de travail de mathématiciens des formes finales données aux résultats mathématiques.} Nous détaillons notamment dans cet article le recours de Poincaré à la réduction canonique des substitutions introduite par Camille Jordan et à l'idée de réduite développée par Charles Hermite en théorie des formes. Nous montrons que l'usage de tels procédés s'inscrit dans certaines cultures algébriques spécifiques. Il s'accompagne ainsi de philosophies de la généralité ainsi que de valeurs épistémiques de simplicité ou d'effectivité. 

L'algèbre figure parmi les principaux points aveugles de l'historiographie de la première décennie de travaux de Poincaré. Elle se prête d'ailleurs mal aux découpages disciplinaires par lesquels ont été organisées les \textit{\oe{}uvres complètes} du savant.\footnote{Ces classements par disciplines, loin d'être neutres, ont été largement basés sur l'analyse faite par Poincaré de ses propres travaux \cite{Poincar1921}, analyse elle même basée sur la forme institutionnelle des candidatures à l'Académie.} En effet, l'algèbre n'avait pas au XIX\up{e} siècle en France un statut de discipline de recherche \cite{Brechenmacher:2010b}. D'un côté, elle renvoyait traditionnellement à une discipline scolaire (l' "algèbre élémentaire") ou intermédiaire (l' "algèbre supérieure") menant vers le point de vue plus élevé de l'Analyse enseignée à l'École polytechnique. D'un autre côté, l'algèbre était souvent considérée par les savants comme un ensemble de procédés opératoires circulant entre les différentes branches des sciences mathématiques. Ceci est illustré par la manière dont Jordan présentait l'algèbre dans l'introduction de son \textit{Cours d'analyse de l'École polytechnique} :
\begin{quote}
Les Mathématiques sont les sciences des quantités. Elles se divisent en plusieurs branches, suivant la nature des grandeurs soumises au calcul. On y distingue principalement l'Arithmétique, la Géométrie, la Mécanique, la Physique mathématique, le Calcul des probabilités.
Ces diverses branches ont pour lien commun l'Algèbre, qu'on pourrait définir [comme] le calcul des opérations.\cite[p.1]{Jordan1882b}.
\end{quote}
Questionner les pratiques algébriques de Poincaré donne donc l'occasion d'étudier la première décennie des travaux du géomètre sans les cloisonner \textit{a priori} dans une discipline particulière. Les procédés opératoires attachés à la \textit{réduction canonique} de Jordan ou au \textit{calcul des Tableaux} traversent notamment \textit{tous} les tomes des \oe{}uvres de Poincaré. Nous montrons que ces procédés donnent un nouvel éclairage sur le lien entre l'émergence de la célèbre théorie des fonctions fuchsiennes et des travaux sur les formes quadratiques et cubiques. Bien qu'il n'ait été que rarement commenté par l'historiographie, ce lien a joué un rôle clé dans l'attribution par Poincaré d'une place centrale à la notion de groupe à partir de 1881-1882. 

Il nous faut également examiner en quoi les pratiques algébriques employées par Poincaré sont spécifiques à ce dernier. Nous contrastons à cette fin les innovations individuelles de Poincaré avec différentes échelles collectives et temporelles. Nous questionnons notamment les rôles joués par les héritages des travaux d'Hermite et de Jordan. Les premiers travaux de Poincaré s'appuient notamment sur les travaux d'Hermite sur la réduction des formes algébriques.\cite[p.391-396]{Goldstein:2007}.  Nous montrons ici que c'est dans cet héritage hermitien que Poincaré s'approprie l'approche de Jordan sur les groupes linéaires. 

La réduction canonique de Jordan joue un rôle clé dans les pratiques algébriques de Poincaré. Il nous faut par conséquent revenir sur les travaux de Jordan des années 1860 ainsi que sur leur réception dans les années 1870-1880. Une étude fine de l'appropriation par Poincaré des travaux de Jordan nous donner aussi l'occasion de préciser l'héritage de ce dernier dans l'évolution de la théorie des groupes et dans l'émergence de la théorie des fonctions fuchsiennes.

Il s'agit donc dans un premier temps d'analyser les dimensions collectives des pratiques algébriques de Poincaré afin d'éclairer dans un second temps leurs spécificités individuelles.

La première partie de cet article prend pour point de départ une courte note publiée par Poincaré en 1884 et intitulée "Sur les nombres complexes". Dans l'objectif d'éviter les écueils des approches rétrospectives, ce texte est pris comme  point d'ancrage à partir duquel nous construisons progressivement une problématique et les différents corpus pertinents pour l'étude de celle-ci.  Nous questionnons notamment le rôle attribué à la note de 1884 par l'historiographie de l'algèbre. Nous montrons que ce rôle fait écho à une lecture rétrospective de ce texte dans le cadre de la théorie des algèbres associatives et de Lie.\footnote{Sur la notion historique de lecture d'un texte, voir \cite{Goldstein:1995}.}  Au sein de l'\oe{}uvre de Poincaré elle même, l'identité de la note n'est en revanche pas portée par une théorie mais par une pratique algébrique permettant la classification de groupes linéaires par réduction canonique de leurs substitutions. 

Nous analysons alors dans une seconde partie l'introduction par Jordan de la réduction canonique des substitutions linéaires dans le cadre de travaux sur les groupes résolubles dans les années 1860. Nous présentons également les usages de la réduction de Jordan dans les années 1870 dans différents cadres théoriques parmi lesquels l'intégration algébrique des équations différentielles linéaires et la théorie algébrique des formes. 

Dans la troisième partie de cet article, nous étudions l'appropriation par Poincaré de l'approche de Jordan. Cette appropriation se mêle notamment au développement de la théorie des fonctions fuchsiennes. Nous montrons que Poincaré change la nature de la réduction de Jordan en intégrant cette dernière dans l'héritage des travaux d'Hermite. Nous analysons alors l'apport individuel de Poincaré quant aux développements ultérieurs de la pratique de réduction des "Tableaux" ainsi qu'au regard d'autres cadres avec lesquels interfère Poincaré comme les matrices de James-Joseph Sylvester. 

  \section{La note "Sur les nombres complexes"}

Notre choix du texte de 1884 comme point d'ancrage tient à ce que ce-dernier témoigne bien d'une pratique algébrique spécifique qui apparaît - et est utilisée - dans des travaux qui ne relèvent pas d'une seule discipline. D'un point de vue actuel, ce texte concerne les théories des algèbres associatives et des algèbres de Lie. Mais nous allons voir que ce point de vue a pour conséquence d'isoler la note de 1884 dans l'\oe{}uvre de son auteur.\footnote{Poincaré ne s'est en effet pas intéressé à la théorie des algèbres associatives et ses principaux travaux sur les groupes continus sont plus tardifs et publiés après la mort de Lie en 1899 (comme le théorème dit de Poincaré-Birkhoff-Witt). Voir \cite{Schmid1982} et \cite{Ton-That1999}.}  Ce texte a d'ailleurs souvent été pris en exemple de la capacité de Poincaré à saisir au vol des éléments essentiels de théories diverses. Nous montrons au contraire que la note s'inscrit parfaitement dans les travaux menés par ce dernier au début des années 1880. Plus encore, la pratique algébrique qui y est mobilisée permet de jeter un nouvel éclairage sur la cohérence des premiers travaux de Poincaré. 

\subsection{Une présentation rapide de la note de Poincaré}
Nous allons tout d'abord présenter rapidement la note d'un point de vue mathématique actuel tout en conservant certains termes du vocabulaire employé par Poincaré.  

La note se propose de "ramener" le problème de la classification des systèmes hypercomplexes, c'est à dire des algèbres associatives unitaires et de dimension finie sur $\mathbb{C}$, à celui de la classification des groupes continus (ou groupes de Lie) :

\begin{quote}
Les remarquables travaux de M. Sylvester sur les matrices ont attiré de nouveau l'attention dans ces derniers temps sur les nombres complexes analogues aux quaternions de Hamilton. Le problème des nombres complexes se ramène facilement au suivant :

\textit{Trouver tous les groupes continus de substitutions linéaires à $n$ variables dont les coefficients sont des fonctions linéaires de $n$ paramètres arbitraires.}
\cite[p.740]{Poincar1884d}
\end{quote}

Considérons un système hypercomplexe de base $e_1$, ..., $e_n$ dont la table de multiplication est définie par les $n^3$ constantes de structures $a_{ij}^k$ telles que
\[
e_ie_j=\sum_{i=1}^{n}a_{ij}^ke_k
\]
À un élément générique $u=\sum_{i=1}^{n}u_ie_i$, on associe la transformation $T_u$ : $x \to y=xu$.\footnote{Il s'agit donc de la représentation régulière d'une algèbre associative.}  Cette transformation est linéaire pour les variables $x_i$ et $y_k$ et pour les paramètres $u_j$ tels que :\footnote{Dans les années 1890, Cartan a désigné sous le nom de "groupes bilinéaires" les monoïdes formés de telles transformations à $n$ variables et $n$ paramètres indépendants dans les équations duquel les paramètres entrent linéairement.} 
\[
y_k=\sum_{i=1}^{n}[\sum_{j=1}^{n}a_{ij}^ku_j]x_i
\]
La loi de composition de ces transformations est associée à la multiplication du système hypercomplexe : $T_uT_v=T_{uv}$. Chaque élément de base $e_i$ de $H$ peut ainsi être représenté par une matrice $E_i$ de $M_n(\mathbb{C})$ et le produit $e_i e_j$ par le produit matriciel $E_iE_j$. Une matrice correspondant à un élément générique $T_u=\sum_{i=1}^{n}u_i E_i$ est non singulière. L'ensemble des matrices non singulières forme un sous-groupe de Lie à $n$ paramètres de $Gl_n(\mathbb{C})$. Les transformations infinitésimales du groupe de Lie $Gl_n(\mathbb{C})$ forment l'algèbre de Lie $M_n(\mathbb{C})$ (le crochet étant le commutateur). 

Ce procédé de représentation régulière permet l'étude de la structure des algèbres associatives par l'étude des structures polynomiales d'invariants comme l'équation caractéristique $det(T_u-sI)=0$.

Jusqu'à la première décennie du XX\up{e} siècle, la plupart des approches sur la classification des algèbres associatives ou des groupes continus s'étaient s'appuyées sur de tels procédés polynomiaux. Ces procédés n'en étaient pas moins variés comme nous le verrons par la suite. La caractérisation de la spécificité des procédés employés par Poincaré est l'un des principaux problèmes que nous posons dans cet article.

La note de 1884 est basée sur l'examen des formes canoniques dites de Jordan auxquelles peuvent être réduites les substitutions linéaires à coefficients complexes. De telles formes dépendent de la multiplicité des racines de l'équation caractéristique. Le texte traite d'abord en quelques lignes du cas des groupes commutatif. Une classification en est donnée en réduisant les substitutions aux types suivants de Tableaux canoniques:\footnote{Ce résultat revient à dire que les matrices commutant avec toutes les matrices commutant avec une matrice donnée $x$ s'expriment comme des polynômes scalaires de $x$. Voir \cite[p.106]{Wedderburn1934} pour une formulation matricielle. }

\begin{quote}
Convenons d'écrire les coefficients d'une substitution quelconque sous la forme d'un Tableau à double entrée. Nous trouverons d'abord que les faisceaux qui donnent naissance à des nombres complexes à multiplication commutative rentrent tous dans des types analogues à ceux qui suivent, pourvu que les variables soient convenablement choisies.

$\begin{vmatrix}
a & 0 & 0 & 0 & 0 \\
0 & b & 0 & 0 & 0\\
0 & 0 & c & 0 & 0\\
0 & 0 & 0 & d & 0 \\
0 & 0 & 0 & 0 & e
\end{vmatrix}$
$\begin{vmatrix}
a & 0 & 0 & 0 & 0 \\
b & a & 0 & 0 & 0\\
c & b & a & 0 & 0\\
d & c & b & a & 0 \\
e & d & c & b & a
\end{vmatrix}$
$\begin{vmatrix}
a & 0 & 0 & 0 & 0 \\
b & a & 0 & 0 & 0\\
0 & 0 & c & 0 & 0\\
0 & 0 & d & c & 0 \\
0 & 0 & 0 & d & c
\end{vmatrix}$
$\begin{vmatrix}
a & 0 & 0 & 0 & 0 \\
b & a & 0 & 0 & 0\\
0 & 0 & c & 0 & 0\\
0 & 0 & d & c & 0 \\
0 & 0 & 0 & 0 & e
\end{vmatrix}$

$a$, $b$, $c$, $d$, $e$ désignant cinq paramètres arbitraires.\cite[p.740]{Poincar1884d}

\end{quote}

Le reste de la note est consacré au cas non commutatif. En s'appuyant à nouveau sur des décompositions de Tableaux à leurs formes canoniques, Poincaré introduit une distinction entre deux types de systèmes de nombres complexes, selon que ceux-ci contiennent ou non les quaternions. En termes actuels, ces deux systèmes prennent respectivement les formes $M_{n_1}(\mathbb{C})$+...+$M_{n_k}(\mathbb{C})$+$N$  et $\mathbb{C}+\mathbb{C}+...+\mathbb{C}+N$ où $N$ est le radical de l'algèbre. Les premiers possèdent une série de composition d'idéaux bilatères et contiennent une sous-algèbre isomorphe à l'algèbre des quaternions.

\begin{quote}
Si l'on considère ensuite un groupe donnant naissance à des nombres complexes à multiplication non commutative, et une substitution quelconque $S$ de ce groupe ; si l'on forme l'équation aux multiplicateurs de cette substitution (équation aux racines latentes des matrices de M. Sylvester), cette équation aura toujours des racines multiples. [...] Supposons maintenant que les variables aient été choisies de telle sorte qu'une substitution $S$ du groupe, non parabolique, soit ramenée à la forme canonique
\[
(x_1, x_2, ..., x_n, \ \lambda_1x_1, \ \lambda_2x_2, ..., \ \lambda_nx_n)
\]
[...] Supposons qu'il y ait $p$  valeurs distinctes de $\lambda$ que nous appellerons $\lambda_1$, $\lambda_2$, ... $\lambda_p$. Nous diviserons les $n$ variables en $p$ systèmes : 
\[
x_{11}, x_{12}, ..., x_{1\alpha}, x_{22}, ..., x_{2\beta}, ..., x_{p1},x_{p-2},..., x_{p\chi},
\]
où $\alpha+ \beta+ ... + \chi=n$ et nous supposerons que la substitution $S$ s'écrive sous la forme
\[
(x_{ik}, \ \lambda_ix_{ik})
\]
le multiplicateur étant ainsi le même pour toutes les variables d'un même système. Cela posé : \\
1° La substitution $(x_{ik}, \ \mu_ix_{ik})$ fera partie du groupe quelles que soient les valeurs des $p$ multiplicateurs $\mu_1$, $\mu_2$, ..., $\mu_p$. \\
2° Écrivons le Tableau à double entrée des coefficients d'une substitution quelconque du groupe, en conservant les mêmes variables dont il vient d'être question. Dans ce Tableau, séparons par des traits verticaux les $\alpha$ premières colonnes, puis les $\beta$ suivantes etc., puis les $\chi$ dernières. Séparons de même par des traits horizontaux les $\alpha$ premières lignes, puis les $\beta$  suivantes, etc., puis les $\chi$  dernières. Nous avons partagé nos coefficients en $p^2$ systèmes. Si l'on choisit  convenablement les $n$ paramètres arbitraires en fonctions desquels tous les coefficients du groupe s'expriment linéairement , un quelconque d'entre eux ne pourra entrer que dans les coefficients d'\textit{un seul} des $p^2$ systèmes.
Il résulte de là : \\
1° Ou bien que les coefficients d'un des $p^2$ systèmes sont tous nuls :  c'est ce qui arrive, par exemple, au groupe à trois variables et trois paramètres
\[ 
\begin{vmatrix}
a & 0 & 0 \\
0 & a & 0 \\
0 & 0 & c 

\end{vmatrix}
\]
2° Ou bien qu'aucune des substitutions du groupe ne peut avoir plus de $\sqrt{n}$ multiplicateurs distincts. C'est ce qui arrive, par exemple, pour les quaternions.\cite[p.741]{Poincar1884d}
\end{quote}

Soulignons dès à présent quelques caractéristiques importantes de la note de 1884 :
\begin{itemize}
\item La référence aux matrices de Sylvester
\item  Le rôle joué par l'équation caractéristique
\item  Le lien entre le problème de la non commutativité et l'occurrence de racines caractéristiques multiples 
\item Les rôles joués par les Tableaux comme moyen d'expression des résultats ainsi que comme procédé opératoire de décomposition des variables et de leurs indices afin de réduire les substitutions à leurs formes canoniques
\item  La tension entre la généralité revendiquée pour des résultats sur $n$ paramètres et leurs présentations pour 3 à 5 paramètres. 
\end{itemize}
Il nous faudra discuter des dimensions collectives de chacun des points ci-dessus afin de dégager la spécificité individuelle de l'approche de Poincaré.

\subsection{Interpréter l'usage des matrices dans la note "Sur les nombres complexes"}
Comme nous l'avons déjà mentionné, une difficulté de notre analyse est de ne pas prendre comme allant de soi les organisations disciplinaires actuelles. Parmi celles-ci, l'algèbre linéaire est particulièrement problématique car, non seulement elle n'existait pas en tant que discipline au XIX\up{e} siècle, mais ses principaux objets n'avaient en outre pas le caractère élémentaire qu'ils revêtent aujourd'hui. Nous montrons dans ce paragraphe qu'interpréter en des termes actuels l'usage fait par Poincaré de la notion de matrice ne permet pas de saisir la place de la note de 1884 dans l'\oe{}uvre de ce dernier. 

\subsubsection{La note de Poincaré comme point d'origine}
 Il nous faut tout d'abord être attentif aux interprétations qui ont été données du texte de 1884 aux XIX\up{e} et XX\up{e} siècles. Comme nous allons le voir, ces interprétations engagent simultanément les mathématiques et l'écriture de leur histoire.

Des travaux consacrés à l'histoire de l'algèbre linéaire ont interprété la note de Poincaré comme une origine des relations mathématiques entre groupes continus et systèmes hypercomplexes. A la suite de Thomas Hawkins [1972], Karen Parshall [1985] a présenté ce texte comme reliant deux traditions distinctes, ancrées en Angleterre et aux États-Unis pour l'une, sur le continent Européen pour l'autre. Pour Parshall, la "tradition anglo-américaine" regroupe des travaux d'auteurs comme Arthur Cayley ou Benjamin Peirce. Elle est caractérisée par la reconnaissance commune d'un acte fondateur dans la découverte des quaternions par William Rowan Hamilton en 1843, ainsi que par une position épistémologique partagée consistant à  considérer les algèbres comme des objets mathématiques satisfaisant certaines propriétés. D'un autre côté, la "tradition de la théorie de Lie" désigne des recherches effectuées sur le continent européen sur les groupes de transformations et ne visant pas à  étudier des entités algébriques pour elles mêmes.

La note du 3 novembre 1884 viendrait ainsi jeter un pont entre deux théories identifiées chacune à une aire géographique et culturelle. Plus encore, la rencontre initiée par Poincaré serait à l'origine d'une unification théorique permettant l'emploi des méthodes de 1898la théorie des groupes pour l'étude de la structure des algèbres associatives :\footnote{Pour un état des lieux des recherches sur les systèmes hypercomplexes au tournant des XIX\up{e} et XX\up{e} siècles, voir \cite{Cartan:1908}. Sur les groupes continus et transformations infinitésimales, voir \cite{Burkhardt:1916}.}
\begin{quote}
By the late 1890s, as a result of the research efforts on both sides of the Atlantic, hypercomplex number systems came to define a distinct area of math investigation [...]. Poincaré recognized that \\
(1) the algebras analogous to the quaternions which Sylvester was studying were algebras of matrices, \\
(2) each element in these algebras defined a linear transformation, and \\
(3) the theory of continuous transformation groups developed by Lie could be applied to these linear transformations.\cite[pp.227 ; 262]{Parshall1985} 
\end{quote}
Ces commentaires situent donc la note dans des évolutions de théories algébriques. C'est dans ce cadre qu'ils se sont attachés à en situer l'origine. Les premières lignes du texte ont été interprétées comme énonçant qu'à un élément $u$ d'une algèbre $A$, il est possible d'associer une transformation, la translation à gauche par $u$, permettant de représenter $A$ comme une sous-algèbre de $M_n(\mathbb{C})$, elle même susceptible d'être munie d'une structure d'algèbre de Lie associée au groupe de Lie $Gl_n(\mathbb{C})$. 

L'origine de la note a alors été attribuée à la reconnaissance par Poincaré du caractère élémentaire de l'algèbre des matrices à la fois en tant qu'algèbre associative et en tant qu'algèbre de Lie. Cette conception des matrices a été qualifiée de "préliminaire indispensable" à la mise en relation de ces deux théories : 
\begin{quote}
It was undoubtly Sylvester's matrix representation of quaternions and nonions that prompted Poincaré's remark. In fact, the examples that Poincaré gave of "bilinear groups" (as we shall call them following E. Cartan) were presented in matrix form [...]. Also of particular importance is the recognition of the special class of hypercomplex systems which we shall term complete matrix algebras. By this we mean systems whose elements can be expressed in the form $\sum a_{ij}e_{ij}$ where the $n^2$ basis elements $e_{ij}$ multiply according to the rule  $e_{ij}e_{kl} = \delta_{jk}e_{il}$. Interest in complete matrix algebras and recognition that the (complex) quaternions are of this type was a necessary preliminary to, and instrumental in, the discovery of the important role they play in the general structure theory of hypercomplex systems.\cite[p.249]{Hawkins1972}
\end{quote}
Il est remarquable que l'attribution à la note de 1884 d'un rôle d'\textit{origine} reflète une architecture mathématique dont nous sommes contemporains, présentant l'algèbre des matrices comme un préliminaire \textit{élémentaire} à des développements plus complexes en algèbre linéaire. 

De fait, l'interprétation \textit{historique} de la note de 1884 comme jetant un pont entre deux théories ou aires géographiques provient en grande partie d'une lecture rétrospective basée sur une organisation \textit{mathématique} attribuant un caractère élémentaire à la notion de matrice. Comme nous allons le voir à présent, cette organisation mathématique s'est mise en place, localement, à partir de 1890 mais n'a fait l'objet d'une culture partagée qu'à partir des années trente du XX\up{e} siècle.\footnote{L'écriture de l'histoire par des textes mathématiques a été étudiée à de nombreuses reprises. Voir, entre autres, \cite{Cifoletti:1995}, \cite{Goldstein:1995}, \cite{Brechenmacher:2006b}, \cite{Brechenmacher:2007a}.}

\subsubsection{Structuration mathématique et structuration historique}

On trouve déjà une interprétation du texte de Poincaré semblable à celle donnée par Hawkins et Parshall dans un mémoire de Georg Scheffers de 1891. Ce travail présente une étude systématique de la classification des systèmes hypercomplexes sur le modèle de celle des groupes continus. Il s'est imposé comme une référence à la fin du XIX\up{e} siècle. Il a notamment contribué à fonder un cadre théorique propre aux études sur les algèbres associatives. Comme souvent lors de la délimitation d'un nouveau domaine des mathématiques, le mémoire de Scheffers n'a pas uniquement proposé une synthèse théorique mais aussi une présentation historique des origines de la théorie des nombres hypercomplexes.

D'un point de vue théorique, Scheffers a pris modèle sur les travaux de Sophus Lie en transposant la distinction entre groupe continu résoluble (intégrable) et non résoluble aux systèmes hypercomplexes. Il a ainsi introduit la distinction entre deux types de tels systèmes selon qu'ils contiennent ou non les quaternions (ou, plus généralement, une algèbre de matrices)\cite[p.293]{Scheffers1891}. Cette distinction a permis de premiers résultats sur la structure des algèbres associatives et a aussi joué un rôle important pour de nombreux travaux ultérieurs comme ceux de Theodor Molien ou d'Élie Cartan.  

C'est aussi à la demande de Lie que la présentation historique de Scheffers a attribué à la note de Poincaré l'origine de la distinction entre les systèmes contenant ou non les matrices. Cette distinction a pourtant été élaborée par Scheffers indépendamment de Poincaré.\footnote{Voir \cite{Hawkins1972}, \cite{Hawkins2000}. Avant Scheffers, Study avait déjà adressé en 1889 un courrier à Poincaré dans lequel il attribuait à la note de 1884 l'origine des relations entre systèmes hypercomplexes et groupes de Lie.\cite{Study1889}} C'est en effet en étudiant de manière empirique les tables de multiplications de systèmes hypercomplexes particuliers que Scheffers a constaté que l'invariant polynomial donné par l'équation minimale du système se décompose en facteurs linéaires lorsque le système ne contient pas de copie des quaternions. 

Examinons plus avant la présentation théorique et historique de Scheffers. L'algèbre des matrices est au c\oe{}ur de l'organisation théorique. Elle donne en effet à la fois l'exemple le plus simple d'algèbre associative et d'algèbre de Lie. Les matrices sont pour cette raison dotées d'un statut élémentaire. 

La présentation historique reflète ce caractère élémentaire en présentant les matrices comme un double point d'origine de la théorie des algèbres. Une première origine du concept de matrice comme nombre hypercomplexe est ainsi attribuée à Hamilton (1852) et Cayley (1858). Des redécouvertes multiples lui succèdent sur le continent, comme chez Edmond Laguerre (1867) ou Georg Frobenius (1878). La note de Poincaré intervient dans le cadre de cette seconde origine en raison de sa reconnaissance des relations établies par les matrices entre algèbres associatives et de Lie. A la suite de la note de Poincaré, les matrices auraient alors permis d'unifier des travaux continentaux de Carl Gauss, Hermann Grassmann, Frobenius ou Lie et des travaux anglo-américains comme ceux de Cayley, Peirce, William Clifford etc.\footnote{Scheffers se situe lui même dans la tradition continentale. Son mémoire fait en effet suite à des travaux de Dedekind, Kronecker et Study sur des systèmes commutatifs qui ont été impulsés par des remarques de Weierstrass en 1884 sur l'arithmétique des entiers de Gauss.}

L'organisation mathématique et historique que l'on trouve localement dans les années 1890 dans des travaux comme ceux de Scheffers a acquis un caractère global avec le développement de la théorie des matrices à une échelle internationale dans les années 1930 \cite{Brechenmacher:2010a}.

Par exemple, le traité intitulé \textit{The theory of matrices} de Cyrus Colton Mac Duffee attribue à Poincaré l'origine du résultat selon lequel toute algèbre associative unitaire simple sur $\mathbb{C}$ peut être représentée comme une sous-algèbre d'une algèbre de matrices. La référence à la note de 1884 intervient dès les premières pages de l'ouvrage de Mac Duffee : le théorème de Poincaré précède la définition même de la notion de matrice à laquelle il vient donner une légitimité théorique. Ce traité contient aussi des notes bibliographiques qui visent à identifier systématiquement les origines des mathématiques exposées. Il a d'ailleurs souvent été cité comme une référence pour l'histoire des matrices avant que ne soient publiés des travaux historiques comme ceux d'Hawkins ou Parshall. Comme chez Scheffers, la référence à la note de Poincaré apparait chez Mac Duffee à la croisée de travaux continentaux sur les groupes de Lie et anglo-américains sur les algèbres associatives. 

Les exemples des synthèses de Scheffers et Mac Duffee illustrent que le miroir, par lequel structurations historique et mathématique de la théorie des algèbres semblent se refléter l'une et l'autre, a précisément été poli par les mathématiciens qui ont conféré à la notion de matrice un caractère élémentaire dans leurs théories. 

\subsubsection{Les problèmes posés par les lectures rétrospectives de la note}
Cette lecture de la note de 1884 pose cependant de nombreux problèmes. Elle rend tout d'abord ce texte très singulier dans l'\oe{}uvre de Poincaré. La note semble en effet l'unique travail consacré par ce dernier aux algèbres associatives. En outre, le caractère fondateur attribué à la note est difficilement conciliable avec le fait que l'approche qui y est mise en \oe{}uvre n'a en réalité pas été reprise par la plupart des mathématiciens qui ont développé la théorie des systèmes hypercomplexes, à commencer par Scheffers. Plus encore, nous avons vu que la notion de matrices est au c\oe{}ur des lectures faites par Scheffers, Mac Duffee, Parshall ou Hawkins. Or, le terme matrice n'est que très rarement employé dans les travaux de Poincaré et n'y est par ailleurs pas mobilisé de manière opératoire. 

Pourtant, et bien que ce dernier distingue nettement les matrices de Sylvester des Tableaux qu'il emploie lui même, ces notions ont par la suite souvent été interprétées toutes deux comme des matrices au sens actuel du terme. Nous allons voir dans le prochain paragraphe que les matrices prennent en réalité chez Poincaré des significations différentes de celles qui leur sont attachées aujourd'hui. 

La question est importante car l'anachronisme terminologique porté par des éléments d'algèbre linéaire comme les matrices n'implique pas seulement des résonnances conceptuelles. Celles-ci se doublent en effet souvent d'un anachronisme sociologique qui fait obstacle à l'analyse des échanges scientifiques ou de la valorisation des résultats.\footnote{Il s'agit là d'un résultat de \cite{Goldstein:2009} qui éclaire la position singulière de Fermat dans une configuration sociale de pratiques et de savoirs dont il domine le fonctionnement et les enjeux propres. Voir aussi \cite[p.189]{Goldstein:1999}} De fait, nous allons montrer que la note de 1884 ne joue à cette date aucun rôle particulier de pont entre grands courants de recherches ou aire culturelles.

Une fois ce constat effectué, il nous faudra alors poser d'une nouvelle manière la question des significations prises par la note de 1884 pour Poincaré lui-même ainsi que celle des cadres collectifs dans lesquels prenait place ce texte au moment de sa parution.

\subsection{La référence de Poincaré aux matrices de M. Sylvester}

\begin{quote}
Si l'on considère [...] un groupe donnant naissance à des nombres complexes à multiplication non commutative, et une substitution quelconque $S$ de ce groupe, si l'on forme l'équation aux multiplicateurs de cette substitution (équation aux racines latentes des matrices de M. Sylvester), cette équation aura toujours des racines multiples.\cite[p.740]{Poincar1884d}
\end{quote}

Ce passage, dans lequel Poincaré explicite la relation entre son propre texte et les travaux de Sylvester, ne vise pas à identifier une notion unificatrice entre théories des groupes continus et des nombres complexes. Le transfert de résultats entre les deux domaines procède au contraire de la reconnaissance d'une identité entre les procédés polynomiaux mobilisés pour ces deux types de problèmes. Tous deux sont en effet abordés par l'intermédiaire d'une équation spécifique dans laquelle nous reconnaissons aujourd'hui l'équation caractéristique d'une matrice. 

Il faut donc insister sur le fait que l'identité algébrique du lien entre les travaux de Poincaré et de Sylvester ne se présente pas sous une forme théorique. Elle ne procède notamment pas du caractère unificateur des notions abstraites dont nous sommes aujourd'hui familier en manipulant espaces vectoriels ou algèbres associatives. C'est une équation spéciale qui permet d'identifier des procédés et des problèmes communs. 

A la fin du XIX\up{e} siècle, cette équation pouvait être désignée de diverses manières : équation en $s$, équation séculaire, équation aux multiplicateurs, équation aux racines latentes etc. La diversité d'appellations s'accompagnait d'une variété de théories mobilisées : formes quadratiques ou bilinéaires, systèmes différentiels, substitutions, matrices etc. Cette variété n'implique cependant pas que ces théories se soient développées en autarcie les unes par rapport aux autres. Elle nécessite au contraire de questionner la manière dont certaines approches mathématiques ont pu circuler à différentes échelles en l'absence d'une organisation disciplinaire globale comme celle que donne aujourd'hui l'algèbre linéaire.

Comme nous allons le voir, la référence de Poincaré à Sylvester donne un exemple typique de point de contact entre des pratiques algébriques différentes. De tels points de contacts éclairent d'un nouveau jour ce que nous percevons rétrospectivement comme des redécouvertes multiples.\footnote{Ainsi, bien que Poincaré ait mis en relation la pratique de réduction canonique et les matrices en 1884, la forme de Jordan n'a été exprimée sous forme matricielle qu'entre 1901 et 1910 et n'a pris l'identité d'un théorème à un niveau international que dans les années 1930.\cite{Brechenmacher:2006a} }

Il faut pour cela envisager une histoire culturelle qui ne se réduise pas à des courants de recherche identifiés à des aires géographiques ou à des institutions. De nombreux travaux ont notamment montré la nécessité d'étudier des espaces de relations intertextuelles - ou réseaux de textes - afin d'analyser des innovations individuelles au regard de dynamiques collectives.\footnote{Voir notamment \cite[p.204-212]{Goldstein:1999},  \cite[p.72-75]{GoldsteinSchappa:2007c}, \cite{Brechenmacher:2007a}.} C'est en effet par la référence à un groupe de textes partagé que Poincaré identifie un noyau commun de nature algébrique entre ses propres procédés et ceux de Sylvester. 

En ce sens, la note de 1884 s'inscrit dans une culture algébrique commune qu'il nous faut caractériser afin de pouvoir, par la suite, questionner la spécificité de l'approche de Poincaré.

\subsubsection{La réapparition de la notion de matrice dans les travaux de Sylvester}

Lorsqu'il évoque les matrices en 1884, Poincaré se réfère implicitement à une série de travaux publiée par Sylvester dans les \textit{Comptes rendus}. 

Publiées en 1882, les deux premières notes de cette série proposent de caractériser les substitutions homographiques $\phi$ de "périodicité donnée",  c'est à dire telles qu'il existe $\mu$ tel que  $\phi^\mu(x)=x$. Le problème est abordé en exprimant les racines $\mu\up{e}$ d'une telle substitution par des fonctions numériques de ses racines caractéristiques. Cette approche était classique depuis le milieu du XIX\up{e} siècle pour le cas de substitutions à deux variables $\phi(x)= \frac{ax+b}{cx+d}$. L'enjeu des travaux de Sylvester est de la généraliser au cas de $n$ variables. Or cet objectif de généralité a immédiatement fait surgir des difficultés liées à l'occurrence de racines caractéristiques multiples dans l'équation $det(\phi-\lambda I)=0$. Ces difficultés ont par la suite été liées par Sylvester à l'étude de la commutativité ou non commutativité de deux matrices. 

Comme nous le détaillons plus loin (paragraphe 1.5), la note de Poincaré généralise elle aussi à $n$ variables des résultats antérieurs sur les fonctions homographiques de 2 variables. On peut donc dors-et-déjà observer trois points communs entre ce texte et les travaux de Sylvester :  
\begin{itemize}
\item une forme de généralisation consistant à passer de 2 à $n$ variables
\item l'utilisation de procédés spécifiques associés à l'équation caractéristique
\item  l'étude de systèmes non commutatifs par la prise en compte des ordres de multiplicités des racines caractéristiques
\end{itemize}

Au-delà de ces points communs, les textes de Sylvester et Poincaré mobilisent des procédés, des notions et des objectifs très différents. Afin de mieux comprendre ce qui réunit et distingue ces travaux, intéressons nous de manière plus approfondie aux travaux de Sylvester.

La première note publiée par ce dernier en 1882 est intitulée "Sur les puissances et les racines de substitutions linéaires".\cite{Sylvester1882a} Le principal résultat de ce texte donne une expression de la puissance d'une substitution $\phi$ par une fonction numérique de ses racines caractéristiques $\lambda$. 

Soit $K$ un coefficient quelconque du déterminant $det(\phi)$ ; soit $K_i$ le coefficient qui occupe la même position dans $det(\phi^i)$ que $K$ dans $det(\phi)$ ; soit $\sigma_i$ les polynômes symétriques élémentaires et $S_i=\sigma_i(\lambda_1, \lambda_2,...,\lambda_n)$. Sylvester énonce en 1882 la formule suivante :
\[
K_i=\sum \frac{K_{n-1}-S_1K_{n-2}+S_2K_{n-3}-... \pm S_{n-1}K_0}{(\lambda_1-\lambda_2)(\lambda_1-\lambda_3)...(\lambda_1-\lambda_n)} \lambda_1^i
\]

Cette formule a cependant été immédiatement contredite par un résultat établi antérieurement par Cayley. Le célèbre mémoire sur la théorie des matrices de ce dernier donne en effet une expression de la racine des substitutions homographiques de trois variables qui est en contradiction avec la formule générale de Sylvester. 

De fait, La formule ci-dessus perd toute validité en cas d'occurrence de racines caractéristiques multiples.\footnote{Voir \cite[p.438]{Cartan:1908} à propos des fonctions analytiques dont les variables sont hypercomplexes (par exemple les exponentielles de matrices). Voir \cite[p.27]{Wedderburn1934}, pour un emploi de la forme de Jordan en cas d'occurrence de valeurs propres multiples.} Sylvester a pour cette raison publié une deuxième note intitulée "Sur les racines des matrices unitaires".\cite{Sylvester1882b} Le problème posé par la multiplicité des racines caractéristiques a ainsi suscité la réapparition de la notion de matrice dans l'\oe{}uvre de Sylvester plus d'une trentaine d'année après son introduction. 

\subsubsection{Les matrices mères des mineurs de 1850 à 1880}

C'était en effet afin de surmonter des difficultés posées par l'occurrence de racines caractéristiques multiples que Sylvester avait introduit la notion de matrice en 1851.\cite{Sylvester1851} Son objectif était alors d'étudier les types d'intersections de coniques. Sylvester avait examiné la décomposition polynomiale du déterminant caractéristique en relation avec la décomposition d'une "matrice" qu'il avait envisagée comme une "forme carrée" afin d'en extraire des "mineurs" comme engendrés d'un parent commun :
\begin{quote}
I have in previous papers defined a "Matrix" as a rectangular array of terms, out of which different systems of determinants may be engendered, as from the womb of a common parent ; these cognate determinants being by no means isolated in their relations to one another, but subject to certain simple laws of mutual dependence and simultaneous deperition. \cite[p.296]{Sylvester1851}.
\end{quote}
Cette définition d'une matrice comme mère des mineurs d'un déterminant a bénéficié d'une circulation large sur le continent européen dès les années 1850. Dans ce contexte, le terme matrice a toujours été associé à des procédés d'extractions de mineurs et aux problèmes spécifiques posés par l'occurrence de racines caractéristiques multiples. Les rares emplois de ce terme par  Poincaré sont  bien représentatifs de cet usage des matrices \cite[p.334]{Poincar1884g}, \cite[p.420]{Poincar1884f}. Jusqu'au début des années 1880, les matrices n'ont au contraire que rarement été envisagées comme des nombres complexes. Poincaré, par exemple, se réfère plutôt à la notion de clef algébrique de Cauchy lorsqu'il utilise un calcul symbolique des transformations sur des systèmes de vecteurs au sens de Grassmann.\cite[p.101]{Poincar1886c}. 

En 1882, le problème de l'expression des fonctions rationnelles d'homographies a cependant amené Sylvester à une nouvelle lecture de la notion de matrice qu'avait développé Cayley en 1858. Rappelons que l'objet  de la théorie des matrices de ce dernier est d'énoncer un "théorème remarquable" supportant une méthode d'expression des substitutions homographiques de périodicité donnée (le théorème dit de Cayley-Hamilton) [Brechenmacher 2006b].\footnote{Ce problème, tout comme l'approche développée par Cayley pour le résoudre, se rattachaient au contexte de l'héritage des opérations symboliques de l'école algébrique de Cambridge\cite{ Durand-Richard:1996}. Des problématiques de calculs de racines de fonctions homographiques ont en particulier été développées par Herschell et Babbage dans le contexte de l'importance donnée à la composition des opérateurs différentiels (puissances, racines, non commutativité) depuis l'introduction en Grande Bretagne du calcul différentiel.} L'expression de la racine carrée d'une matrice par un polynôme de degré inférieur à celui de la matrice elle-même avait amené Cayley à introduire des lois d'opérations sur les matrices. Dans une série de notes publiée de 1882 à 1884, Sylvester a interprété ces lois d'opérations en envisageant désormais les matrices comme des quantités complexes généralisant les quaternions d'Hamilton.

\subsubsection{Racines multiples et commutativité}
L'incorporation de la théorie des quaternions dans celle des matrices nécessite la résolution des équations matricielles du type $mn=-nm$ qui définissent les éléments de base des quaternions (ou de leur généralisation à des systèmes d'ordre 3, les nonions). Dans un premier temps, Sylvester a énoncé que si $mn=nm$ alors, soit l'une des matrices est diagonale, soit il existe une relation fonctionnelle entre $m$ et $n$. Mais Sylvester a par la suite constaté qu'un tel énoncé n'est pas toujours valable en cas d'occurrence de racines multiples. Cette nouvelle difficulté posée par les racines caractéristiques multiples a alors conduit Sylvester à approfondir le problème de la détermination de l'ensemble des matrices commutant à une matrice donnée. Son approche consiste à étudier les relations entre décompositions d'une matrice en mineurs et décompositions polynomiales de l'équation caractéristique.\footnote{C'est dans ce contexte qu'a  notamment émergée la  distinction entre polynômes caractéristique et polynôme minimal d'une matrice.}

Comme nous l'avons évoqué précédemment, cette problématique du lien entre commutativité et racines multiples est précisément celle que la note de 1884 aborde à l'aide de la réduction canonique des Tableaux. 

Poincaré n'est cependant pas le seul à avoir répondu aux problèmes posés par Sylvester dans les \textit{Comptes rendus}. En 1885, Eduard Weyr a notamment débuté ses travaux sur ce que nous désignerions aujourd'hui comme la décomposition d'un espace vectoriel de dimension finie en sous-espaces caractéristiques sous l'action d'un opérateur. Entre 1885 et 1890, Weyr a également appliqué ses travaux à l'étude des systèmes hypercomplexes.\cite{Weyr1890} Son approche a par la suite été reprise dans ce cadre par des auteurs comme Scheffers, Frobenius, Molien ou Kurt Hensel \cite{Brechenmacher:2006a}. 

Nous ne nous attarderons cependant pas plus avant ici sur cet aspect de la réception des travaux de Sylvester. En effet, contrairement à d'autres mathématiciens continentaux, Poincaré n'a pas lui même repris les méthodes développées par Sylvester pour l'étude des quaternions et nonions. De fait, la référence du premier au second ne vise que les trois points que nous avons déjà évoqué plus haut : procédés spécifiques associés à l'équation caractéristique, forme de généralisation consistant à passer de 2 à $n$ variables, étude de systèmes non commutatifs par la prise en compte des ordres de multiplicités des racines caractéristiques.

\subsection{Une culture algébrique commune : l'équation séculaire}
Revenons à présent à la note de 1884 et à la question de la mise en relation des nombres complexes et des groupes continus. Nous avons déjà vu que cette mise en relation ne tient pas à l'identification d'une notion unificatrice entre deux grandes théories. Elle s'appuie en revanche sur l'explicitation d'une identité en partie commune entre deux pratiques distinctes : la pratique des matrices de Sylvester et la pratique de réduction canonique utilisée par Poincaré. 

Nous allons à présent montrer que cette identité commune présente une dimension collective qui dépasse les travaux de Sylvester et Poincaré. Elle relève en fait d'une culture algébrique largement partagée dans la seconde partie du XIX\up{e} siècle et qu'identifie très précisément le titre du mémoire dans lequel Sylvester développe, en 1883, la notion de racine latente d'une matrice : "On the equation to the secular inequalities in the planetary theory".\cite{Sylvester1883} Ce titre reprend d'ailleurs celui d'un article  publié par Sylvester en 1851 lorsque ce dernier a introduit le terme matrice : "Sur l'équation à l'aide de laquelle on détermine les inégalités séculaires des planètes".\cite{Sylvester1852}

De telles références à l'équation séculaire interviennent dans de très nombreux textes du XIX\up{e} siècle. Elles ne manifestent le plus souvent pas de préoccupations en mécanique céleste mais s'avèrent la manière dont était identifiée une pratique algébrique commune entre études des cordes vibrantes et mécanique céleste, systèmes différentiels et géométrie analytique, arithmétique des formes quadratiques et algèbre des invariants. C'est d'ailleurs dans ce cadre que les notions de matrices et de mineurs ont été introduites par Sylvester à la suite de travaux de Cauchy et ont par la suite circulé entre des travaux de Cayley, Sylvester, Hermite, Riemann etc.

On trouve au c\oe{}ur de cette pratique algébrique un procédé consistant à donner une expression polynomiale des solutions $(x_i)$ d'un système linéaire symétrique de $n$ équations par des quotients impliquant, outre le déterminant caractéristique $S$ du système, ses mineurs successifs $P_{1i}$ obtenus par développements par rapport à la 1\up{re} ligne et $i$\up{e} colonne. L'expression polynomiale
\[
(*)\frac{P_{1i}}{\frac{S}{x-s_j}}(x)
\]
donne
\[
x_i^{s_j}=\frac{P_{1i}}{\frac{S}{x-s_j}}(s_j)
\]
où $x_i^{s_j}$ désigne le système de solutions associé à la racine caractéristique $s_j$.\cite{Brechenmacher:2007b}.\footnote{En termes actuels, ces procédés reviennent à donner une expression polynomiale générale des vecteurs propres d'une matrice comme quotients de mineurs extraits du déterminant caractéristique $det(A-\lambda I)$ et d'une factorisation de l'équation caractéristique par un terme linéaire : $(\lambda-\lambda_i)$ (où $\lambda_i$ est une racine caractéristique de $A$). Une telle expression polynomiale s'identifie à un facteur près aux colonnes non nulles de la matrice des cofacteurs de $A-\lambda I$.} 

La pratique algébrique associée à l'équation séculaire ne se réduit cependant pas à l'usage de tels procédés polynomiaux. En effet, l'expression (*) est rarement utilisée en tant que telle mais le plus souvent incorporée au sein de différentes méthodes élaborés dans divers cadres théoriques. Ainsi, dans l'approche de Joseph-Louis Lagrange (1766) sur les systèmes différentiels, cette pratique supporte la  transformation d'un système de $n$ équations différentielles linéaires à coefficients constants à la forme intégrable de $n$ équations indépendantes (un système diagonal en termes actuels). Elle est indissociable d'une représentation mécanique selon laquelle des oscillations d'un fil lesté de $n$ masses peuvent se représenter par la composition des oscillations propres de $n$ fils lestés d'une seule masse. La nature des racines de l'équation caractéristique a ainsi été liée à la condition de stabilité mécanique des oscillations. 

Cette pratique a par la suite sous-tendu une mathématisation des oscillations séculaires des planètes sur leurs orbites. Pour cette raison, l'équation sur laquelle elle est basée (l'équation caractéristique) a été dénommée "l'équation à laquelle on détermine les inégalités séculaires des planètes". 

Dans ce contexte, Pierre-Simon Laplace a cherché à démontrer la stabilité du système solaire en déduisant la nature réelle des racines de l'équation séculaire de la propriété de symétrie des systèmes mécaniques. Plus tard, Cauchy s'est appuyé en 1829 sur l'analogie entre les propriétés de cette équation et celles des équations caractéristiques des coniques et quadriques pour donner aux procédés élaborés par Lagrange une interprétation géométrique. La transformation de systèmes linéaires a alors été envisagée en terme de changements d'axes en géométrie analytique. Plus tard encore, ces mêmes procédés ont été associés aux calculs de résidus de l'analyse complexe, aux suites de Sturm d'une équation algébrique, aux classes d'équivalence des couples de formes quadratiques etc.

Ces divers travaux reconnaissaient cependant une nature algébrique sous-jacente aux différents cadres théoriques mobilisées. Cette identité algébrique tient principalement à trois caractères. 
\begin{itemize}
\item Le premier est la permanence sur plus d'un siècle des procédés opératoires (*). 
\item Le deuxième est l'ambition explicite de généralité qui se manifeste depuis la prise en compte par Lagrange de systèmes linéaires  à $n$ variables. 
\item Le troisième est le caractère spécial de l'équation séculaire dont on a longtemps pensé que les racines devaient être non seulement réelles mais aussi distinctes.
\end{itemize}

A partir des années 1840, les références à l'équation séculaire apparaissent le plus souvent associées aux problèmes posés par l'occurrence de racines multiples dans la dite équation. Dans ce cas, l'expression (*) est en effet susceptible de prendre une valeur $\frac{0}{0}$ . Cette difficulté a suscité des critiques quant à la tendance de l'algèbre à attribuer une généralité excessive aux expressions symboliques. C'est d'ailleurs à partir de ce constat que Cauchy s'est détourné de l'algèbre au profit de l'analyse complexe : le calcul des résidus permet en effet de remplacer le problème de la multiplicité des racines par celui de la détermination des degrés des pôles de l'expression (*) envisagée comme une fonction méromorphe.

À partir des années 1850, plusieurs approches algébriques du problème de la multiplicité des racines ont été développées de manière distinctes : matrices de Sylvester en géométrie analytique, réduites d'Hermite en théorie des formes quadratiques, diviseurs élémentaires de Weierstrass pour les couples de formes quadratiques et bilinéaires, réduction canonique de Jordan pour les systèmes différentiels linéaires à coefficients constants. 

Ces différents développements témoignent de la dislocation progressive de la culture algébrique commune qu'a longtemps porté l'équation séculaire. Jusque dans les années 1880, c'est cependant toujours par l'intermédiaire de cette équation que ces lignes divergentes se sont rencontrées épisodiquement. C'est notamment dans ce cadre que Poincaré a réagi en 1884 aux difficultés posées par les racines multiples dans les travaux de Sylvester, sans pour autant s'intéresser aux problématiques sur les matrices-quantités complexes développées à cette époque par ce dernier. 

\subsection{La place de la note de 1884 au sein de l'\oe{}uvre de Poincaré}

Nous allons à présent chercher à situer la place de la note de 1884 dans l'\oe{}uvre de Poincaré. En complément des références explicites aux matrices de Sylvester, d'autres références sont implicitement portées par le vocabulaire employé dans ce texte. Ainsi, l'étude des groupes à paramètres (référence à Lie) manifeste l'ancrage de la note dans une série de textes publiée par Poincaré depuis 1881 à la suite de travaux de Jordan (termes "substitutions", "faisceaux" et "forme canonique") et de Félix Klein ("substitution parabolique"). Ce vocabulaire révèle des liens vers d'autres travaux de Poincaré contemporains de la note du 3 novembre 1884. Comme nous allons le voir, tous ces textes s'appuient sur des réductions canoniques afin d'aborder des problèmes de classification des sous-groupes finis du groupe linéaire.

\subsubsection{Généralisation de la classification des groupes fuchsiens aux groupes hyperfuchsiens}

Le 11 février 1884, Poincaré a publié aux \textit{Comptes rendus} une note intitulée "Sur les substitutions linéaires". Il y répond à des travaux récents d'Emile Picard\cite{Picard1884} sur les substitutions à trois variables (sur $\mathbb{C}$), qui conservent l'hypersphère $x^2+y^2+z^2=1$, et de la forme :
\[
(x, y, z, \ \frac{ax+by+cz}{a''x+b''y+c''z}, \frac{a'x+b'y+c'z}{a''x+b''y+c''z})
 \]
 Poincaré propose de généraliser à de telles substitutions la classification qu'il a donnée en 1881 aux groupes fuchsiens de substitutions linéaires réelles à deux variables qui n'altèrent pas le cercle fondamental :\footnote{C'est-à-dire, l'étude des sous-groupes discrets de $PSL_2(\mathbb{R})$ et de leur domaine fondamental dans le plan hyperbolique. Le cas des sous-groupes de $PSL_2(\mathbb{C})$ et des polyèdres de l'espace hyperbolique a été dénommé "groupes kleinéens".}
 \[
 \frac{ax+by}{a''x+b''y}
  \]
  Cette problématique de généralisation de la classification des groupes fuchsiens à celles des groupes hyperfuchsiens a été développée dans une nouvelle note le 25 février.\cite{Poincar1884b}
  
La classification procède d'une distinction entre substitutions loxodromiques, hyperboliques, elliptiques et paraboliques à la suite des travaux de Jordan sur les groupes de mouvement\cite{Jordan1868b} puis des classifications données par Klein \cite{Klein1875} et Jordan \cite{Jordan1876a} aux sous-groupes finis de $Gl_2(\mathbb{C})$. Poincaré reprend surtout l'approche par laquelle Jordan a étendu en 1878 cette classification aux sous-groupes finis de $Gl_3(\mathbb{C})$\cite{Jordan1878}. Comme ce dernier, il s'appuie en effet sur une classification des formes canoniques des substitutions en fonction de la multiplicité de leurs racines caractéristiques $\alpha, \beta, \gamma$ :
 
\begin{quote} 
 On peut, par un changement convenable de variables, amener cette substitution à l'une des formes suivantes, que l'on peut appeler formes canoniques : \\
$
(A) \ (x, y, z ; \ \alpha x, \beta y, \gamma z), \\
(B)  \ (x, y, z ; \ \alpha x, \beta y+z, \beta z), \\
(C) \ (x, y, z ;  \ \alpha x, \beta y, \beta z), \\
(D) \ (x, y, z, ; \  \alpha x+y, \alpha y+z, \alpha z), \\
(E) \ (x, y, z ;  \ \alpha x, \alpha y+z, \alpha z),
$
\cite[p.349]{Poincar1884a}
\end{quote}

Cette pratique de classification par réduction canonique est mobilisée dans de nombreux autres travaux de Poincaré. Elle y vise toujours à généraliser des résultats obtenus pour 2 variables à 3 ou $n$ variables. Pour cette raison, la première partie du mémoire "Sur les groupes des équations linéaires", publié en 1884 dans \textit{Acta Mathematica}, est entièrement consacrée à l'exposé de la réduction canonique des substitutions linéaires de $n$ variables.\cite[p.300-313]{Poincar1884e}. Le recours à cette pratique dépasse d'ailleurs le cadre des groupes discontinus. Poincaré l'emploie notamment pour déterminer les formes algébriques homogènes de $n$ variables reproductibles par des groupes continus non commutatifs.\cite{Poincar1883b}
\subsubsection{Equations différentielles linéaires}

Envisagée du point de vue de la pratique de réduction canonique qui y est mise en \oe{}uvre, la note de 1884 sur les nombres complexes est donc loin d'apparaître isolée au sein de l'\oe{}uvre de Poincaré. Elle prend place au contraire dans un corpus de textes dominé par des problématiques de classification des groupes associés aux équations différentielles linéaires.

Rappelons que le sujet des équations différentielles linéaires a fait l'objet de très nombreux travaux dans les années 1870-1880.\footnote{Elles impliquaient notamment des auteurs comme Lazarus Fuchs, Frobenius, Wilhelm Thomé, Lie, Klein, Francesco Brioschi, Jules Tannery, Paul Appel, Georges Halphen, Gaston Floquet, Picard, Édouard Goursat, Felice Casorati etc} Jeremy Gray en a proposé une étude historique qui a bien mis en évidence la diversité des approches des acteurs de l'époque : développement en séries, analyse complexe, invariants et covariants, groupes de substitutions ou de transformations etc.\cite{Gray:2000} Il faut aussi souligner que de mêmes objets mathématiques, comme les groupes de substitutions ou les équations, étaient à cette époque associés à des pratiques diverses. Une variété d'analogies entre l'étude des équations différentielles et des équations algébriques était notamment mobilisée : groupes continus de Lie, notion d'équation irréductible de Frobenius, groupes de monodromie de Jordan, recours par Poincaré à la notion de résolvente de Galois, revendication d'une théorie dite de Galois différentiel par Picard etc.\footnote{Au sujet de la théorie de Galois différentiel à la fin du XIX\up{e} siècle, voir \cite{Archibald:2011}}   

Dans ce contexte, l'usage de la réduction canonique des substitutions linéaires se présente comme l'une des spécificités des travaux de Poincaré. Cette pratique de réduction se distingue notamment des approches par lesquelles la plupart des autres auteurs abordaient les problèmes posés par l'occurrence de racines caractéristiques multiples. Certains, comme Lie, se contentaient d'une approche générique en écartant les cas de multiplicité des racines.\footnote{L'approche générique de Lie a été discutée en détail par Hawkins en comparaison avec les travaux indépendants de Killing qui s'appuyaient quant à eux sur les diviseurs élémentaires.\cite{Hawkins2000} Le remplacement des méthodes de Killing par la réduction canonique de Jordan fait d'ailleurs l'objet d'un texte de Poincaré. \cite[p.216-252]{Poincar1901}} De nombreux autres, comme Frobenius ou Klein, s'appuyaient sur des invariants polynomiaux : les diviseurs élémentaires introduits par Weierstrass en 1868. D'autres encore, comme Brioschi, Fuchs ou Tannery recouraient à des trigonalisations des systèmes linéaires \cite[p.135-138]{Tannery1875}.

\subsubsection{La représentation analytique des substitutions}

Il faut aussi être attentif au fait que les procédés de réduction de Poincaré s'appuient sur une notation spécifique. Celle-ci était désignée au XIX\up{e} siècle  sous le nom de "représentation analytique des substitutions". Cette désignation n'identifie pas seulement une notation mais aussi une approche développée par certains travaux sur les groupes finis. Elle s'accompagne ainsi de problèmes et de procédés algébriques spécifiques. 

Donner une représentation analytique à une substitution $S$ d'un ensemble fini de lettres nécessite tout d'abord d'indexer ces lettres, $a_0, a_1, ..., a_{q-1}$, avec $q=p^n$ où $p$ est un nombre premier. $S $ opère ainsi sur des nombres entiers modulo $p$ (cas $n=1$) ou sur leur généralisation aux "imaginaires de Galois", c'est à dire à des corps finis. Il s'agit ensuite de déterminer une fonction analytique $\phi$ telle  $S(a_i)=a_{\phi(i)}$. L'obtention d'une telle représentation analytique permet alors d'employer des procédés de décomposition polynomiale pour classer les substitutions en différents types.

Bien qu'elle soit passée inaperçue de nombreux travaux historiques, la représentation analytique des substitutions a joué un rôle important au XIX\up{e} siècle. Elle a été mise en avant au début du siècle par Louis Poinsot puis Évariste Galois à la suite des travaux de Gauss sur les équations cyclotomiques.  Dans les travaux de Galois, elle s'est accompagnée de procédures de réductions (ou décomposition) des groupes qui s'articulent, d'une part, aux substitutions homographiques de deux variables des équations modulaires et, d'autre part, aux substitutions linéaires de $n$ variables associées aux groupes primitifs résolubles. Or ces deux types de substitutions ont par la suite sous-tendu deux formes distinctes d'héritages des travaux de Galois : le premier a particulièrement été utilisé par Hermite dans les années 1850,\cite{Goldstein:2011} tandis que le second a été développé par Jordan dans les années 1860.\cite{Brechenmacher:2011}

Nous verrons dans la deuxième partie de cet article que la synthèse donnée par Poincaré aux héritages des travaux d'Hermite et de Jordan constitue un apport crucial du premier tant ces deux héritages relèvent de pratiques distinctes des substitutions linéaires. 

\subsection{Conclusion de la première partie}
Le rôle de point d'origine qui a souvent été conféré à la note "Sur les nombres complexes" est basé sur l'attribution à Poincaré de la mise en relation des théories des systèmes hypercomplexes et des groupes continus par l'intermédiaire de la notion de matrice. En conséquence, cette note a été présentée comme une première rencontre entre deux grandes aires culturelles identifiées de manière géographique. Nous avons insisté sur le fait qu'une telle lecture du texte de Poincaré s'avère elle-même indissociable du statut élémentaire attribué à la notion de matrice au sein d'organisations du savoir mathématique qui se sont mises en place, localement, à partir des années 1890, mais n'ont fait l'objet d'une culture partagée qu'à partir des années 1930 avec la constitution de l'algèbre linéaire comme une discipline. 

Une telle lecture s'est perpétuée jusqu'à nos jours dans des textes mathématiques et historiques. Elle nous est par conséquent contemporaine et peut donner l'illusion d'un caractère universel. De fait, l'universalité prêtée à certaines notions élémentaires d'algèbre linéaire semble avoir conférée à celles-ci un caractère naturel, dénué d'histoire. Pourtant, malgré leurs similitudes visuelles, les procédés de décompositions de formes imagées comme les matrices de Sylvester et les Tableaux de Poincaré renvoient à des pratiques et des dynamiques collectives très différentes. 

La note de 1884 n'a en réalité joué aucun rôle crucial dans la rencontre de grands courants de recherches. D'une part, les matrices-mères des mineurs circulaient déjà sur le continent depuis les années 1850 ; d'autre part l'implantation sur le continent des matrices-quantités complexes n'a pas été amplifiée par la note de Poincaré et a davantage été portée par des réponses d'autres auteurs aux travaux de Sylvester. 

Une attention aux liens entre le texte de Poincaré et les travaux de Sylvester nous a amené à reconnaître que la note de 1884 n'est pas sous-tendue par une notion unificatrice mais par une pratique algébrique de classification des groupes linéaires par réduction canonique de leurs substitutions. En outre, la technicité de tels procédés s'avère indissociable d'aspects culturels propres aux espaces dans lesquels des textes ont circulé et ont été lus à différents niveaux. Or, ces espaces ne s'identifient le plus souvent pas simplement à des aires géographiques ou à des configurations disciplinaires :

Premièrement, à des niveaux interpersonnels et sur des temps courts, Poincaré interagissait avec les travaux de Sylvester sans pour autant s'approprier les pratiques d'opérations sur les matrices qui s'ancraient dans des traditions académiques britanniques et américaines. 

Deuxièmement, Poincaré et Sylvester partageaient une culture algébrique commune qui s'est déployée à une échelle européenne sur le temps long du XIX\up{e} siècle. Les références de Sylvester à l'"équation à l'aide de laquelle on détermine les inégalités séculaires des planètes" manifestent notamment une utilisation tardive d'un mode d'identification de procédés, problèmes et approches spécifiques. Comme nous l'avons vu, les références à cette équation spéciale identifient non seulement des procédés de manipulations des systèmes linéaires basés sur des décompositions polynomiales mais aussi des problèmes de commutativité abordés par l'examen de la multiplicité des racines caractéristiques. C'est par l'intermédiaire de ces procédés et problèmes partagés que Poincaré a réagi en 1884 aux travaux de Sylvester. Le premier s'est en effet appuyé sur l'efficacité de sa pratique de réduction canonique pour aborder les problèmes de commutativité pour lesquels le second avait élaboré sa pratique des matrices. 

L'identification de cette  culture algébrique commune nous permet à présent de questionner la spécificité individuelle des pratiques algébriques de Poincaré. Il va nous falloir dans un premier temps analyser les héritages locaux dont les travaux de Poincaré se nourrissent. Il faut pour cela nous intéresser plus en détail aux travaux de Jordan afin d'en étudier l'appropriation par Poincaré.  

\section{La réduction des substitutions linéaires chez Jordan (1860-1880)}

La réduction canonique des substitutions linéaires sur des entiers mod.$p$ a d'abord été énoncée par Jordan en 1868 dans le cadre d'une étude des sous-groupes résolubles de $Gl_2(F_p)$. Elle a ensuite été présentée sous la forme d'un théorème dans le chapitre du \textit{Traité des substitutions et des équations algébriques} consacré au groupe linéaire $Gl_n(F_p)$.\cite[p.114]{Jordan1870}. Comme nous allons le voir dans cette section, ce théorème s'inscrit plus largement dans une pratique de réduction de la représentation analytique des substitutions qui est sous-jacente à la majorité des travaux menés par Jordan dans les années 1860-1870. 

Dès l'introduction de la thèse qu'il a soutenue en 1860, Jordan a attribué un caractère "essentiel" à une méthode de réduction qu'il a rattaché au cadre de la "théorie de l'ordre" en revendiquant l'héritage de Poinsot.\footnote{Au sujet de la théorie de l'ordre de Poinsot, voir \cite{Boucard:2011}.} Cette théorie n'est pas centrée sur un objet mais vise au contraire l'étude des relations entre des classes d'objets. Elle se présente ainsi comme transversale à la théorie des nombres (cyclotomie, congruences), l'algèbre (équations, substitutions), l'analyse (groupes de monodromie et lacets d'intégration des équations différentielles linéaires), la géométrie/ topologie (cristallographie, symétries des polyèdres et des surfaces - y compris de Riemann) et la mécanique (mouvements des solides). Tous les travaux de Jordan publiés dans les années 1860 se rattachent à ces thèmes.\footnote{Cette articulation de différentes problématiques manifeste une dimension collective à une échelle européenne que nous désignons sous le nom du champ de recherche des solides réguliers afin de mettre en avant un motif commun à un ensemble de textes mêlant arithmétique, algèbre, analyse, géométrie-topologie et cinématique. Les interprétations de Schwarz et Klein sur les polyèdres et les surfaces de Riemann en lien avec les équations différentielles linéaires se rattachent notamment à des travaux de Jordan dans ce champ. L'unité de ce champ se disloque cependant progressivement dans les années 1880-1890 pour disparaître au tournant du siècle avec le développement de la théorie des groupes et de la topologie\cite{Brechenmacher:2011b}.} 

Dans ce cadre, Jordan a développé des procédés de réduction d'un groupe en sous-groupes sur le modèle des procédés d'indexation et de représentation analytique des substitutions. Ces procédés sous-tendent les relations entre domaines que vise la théorie de l'ordre. La réduction d'un groupe a ainsi été envisagée par Jordan comme une sorte de dévissage par analogie avec la décomposition du mouvement hélicoïdal d'un solide en mouvements de rotation et de translation. Cette analogie a notamment été développée dans l'étude de 1868 sur les groupes de mouvement de solides polyédraux pour laquelle Jordan s'est inspiré des travaux de cristallographie de Bravais. 

A partir de la fin des années 1860, Jordan a attribué à la réduction des groupes de substitutions le caractère transversal qu'il avait d'abord associé à la théorie de l'ordre. Il a alors appliqué cette réduction à l'étude des équations différentielles et des formes algébriques.

\subsection{La première thèse de Jordan et l'origine du groupe linéaire}

La première thèse présentée par Jordan à la Faculté des sciences de Paris en 1860 est consacrée au problème du "nombre des valeurs des fonctions". Ce problème est l'une des racines de la théorie des groupes de substitutions. Il a émergé de travaux du XVIII\up{e} siècle sur les équations :  la résolubilité par radicaux d'une équation algébrique de degré $n$ avait été mise en relation avec le nombre de valeurs qu'une fonction résolvante de $n$ variables peut prendre lorsque ses variables sont permutées de toutes les manières possibles.  

Étant donnée une fonction $\phi(x_1, x_2,...,x_n)$ de $n$ lettres, une valeur de $\phi$ est une fonction obtenue par permutation des variables pour toute substitution $\sigma \in Sym(n)$ :
\[
\phi^\sigma(x_1,x_2, ..., x_n)=\phi(x_{1\sigma},x_{2\sigma}, ... ,x_{n\sigma})
\]
En général, une fonction peut ainsi prendre $n!$ valeurs différentes mais il peut arriver que quelques unes de ces valeurs soient identiques. Le problème abordé par la thèse de Jordan consiste à déterminer le nombre de valeurs pouvant être prises pour certaines classes de fonctions. En termes actuels, le problème revient à déterminer tous les ordres possibles des sous-groupes du groupe symétrique $Sym(n)$.

Dans sa thèse, Jordan revendique une approche "générale" du problème par des "réductions successives" en sous-problèmes. La première réduction montre que l'étude des fonctions générales transitives (correspondant à des équations algébriques irréductibles) peut se réduire à l'étude de fonctions sur des "systèmes imprimitifs"  de lettres : l'ensemble de toutes les lettres peut dans ce cas se diviser en blocs\footnote{Jordan désignait en réalité les blocs de lettres par le terme "groupe" et ce que nous appellerions des groupes des substitutions par le terme de système conjugué introduit par Cauchy. Rappelons que chez des auteurs antérieurs comme Galois ou Poinsot, le terme groupe concernait intrinsèquement à la fois des arrangements de lettres en blocs de permutations et les opérations sur de tels arrangements par des substitutions.}   - représentés ci-dessous par une succession de lignes - de manière à ce que les substitutions opèrent soit en permutant entre elles les lettres d'une même ligne soit les lignes les unes avec les autres. 
\[
\begin{matrix}
a_1 & a_2 & ... & a_p\\
b_1& b_2 & ... & b_p\\
c_1 & c_2 & ... & c_p\\
. & . & . & .
\end{matrix}
\]
L'enjeu de cette première réduction est d'introduire une indexation des lettres par des suites d'entiers $(1, 2, ..., m)$ de manière à faire opérer les substitutions sur ces entiers. Cette indexation permet alors de caractériser des groupes de substitutions par les représentations analytiques de ces dernières. 

Après avoir réduit le problème à l'étude des fonctions imprimitives, Jordan montre comment poursuivre la réduction aux fonctions primitives, c'est à dire au cas où il n'est pas possible de subdiviser les lettres en plusieurs lignes $\Gamma_1, \Gamma_2,..., \Gamma_m$ que les substitutions permuteraient en blocs. 

Formulons cette seconde étape en des termes actuels. Soit $G$ un groupe transitif opérant sur un ensemble $V$. Un sous-ensemble $V_1$ de $V$ est appelé un bloc d'imprimitivité si $V_1 \ne \emptyset$ et si pour tout $g \in G$, soit $V_1g=V_1$ soit $V_1g \cap V_i = \emptyset$. Si $V_1$ est un tel bloc et $V_1$, $V_2$, ..., $V_m$ ses orbites distinctes $V_1g$ pour $g \in G$, alors $(V_1, V_2, ..., V_m)$ est une partition de $V$. $G$ est dit primitif s'il n'existe aucun bloc non trivial d'imprimitivité. Soit à présent $G$ un groupe résoluble, transitif et imprimitif. Une partition maximale en systèmes d'imprimitivités $(V_1, V_2, ..., V_m)$ peut être mise en correspondance avec une décomposition de $G$ en sous-groupes $\Gamma, \Gamma', \Gamma'', ..., \Gamma^{(n)}$ tels que :
\begin{itemize}
\item $\Gamma, \Gamma', \Gamma'', ...$ laissent respectivement stable chaque système $V_i$, ces groupes sont tous isomorphes à un groupe $\Gamma$
\item le groupe $\Delta$ de permutation des ensembles $V_i$, est isomorphe au groupe quotient $G/\Gamma^n$.
\end{itemize}
Le choix d'une décomposition maximale garantit la primitivité des groupes $\Delta$ et $\Gamma$ qui permettent de dévisser $G$, c'est-à-dire de l'écrire comme le produit semi-direct de $\Delta$ et $\Gamma^n$. 

Cette formulation actuelle ne permet cependant pas de percevoir le rôle crucial que joue la représentation analytique des substitutions dans la thèse de Jordan. Les substitutions y sont en effet décomposées en deux espèces correspondant à deux formes de représentations analytiques des cycles :
\begin{itemize}
\item  La première espèce permute cycliquement les lettres à l'intérieur de chaque bloc $\Gamma_i$ par des substitutions représentées analytiquement par la forme  $(x \ x+a)$.  Dans le cas plus général de lettres indexées par $n$ indices $a_{x, x', x'',...}$ ces substitutions prennent la forme :
\[
a_{x+ \alpha \ mod.p, \ x'+\alpha' \ mod.p, \ x''+\alpha'' \ mod.p, \ ...}
\]
\item La seconde espèce  permute cycliquement les blocs  $\Gamma_1, \Gamma_2, ..., \Gamma_n $ eux-mêmes par des substitutions du type $(x \ gx)$. Dans le cas de $n$ indices, elles prennent ainsi la forme :

\[
a_{ax+bx'+cx''... mod. p, \ a'x+b'x'+c'x''... mod.p, \ a''x+b''x'+c''x''... mod.p}
\]
\end{itemize}

Le principal théorème de la thèse de Jordan énonce alors :  

\begin{quote}
\textbf{Premier théorème de Jordan}
\begin{itemize}
\item Les systèmes primitifs comportent un nombre de lettres donné par la puissance d'un nombre premier $p^n$.
\item Les substitutions sur ces systèmes ont une forme analytique linéaire : 
\[
(x, x', x'' ... ; \ ax+bx'+cx''+... +d, \ a'x+b'x'+c'x''+...+d', \ a''x+b''x'+c''x''+..+d'', ...)
\]
que l'on peut également dénoter par : 
\[
A=
\begin{vmatrix}
x & ax+bx'+cx''+... \\
x' & a'x+b'x'+c'x''+... \\
x'' & a''x+b''x'+c''x''+... \\
.. & .....................
\end{vmatrix}
\]
\end{itemize}
\end{quote}

Le groupe constitué de telles substitutions a été désigné par Jordan sous le nom de "groupe linéaire". 

En termes actuels,  un groupe primitif résoluble $G$ a pour sous-groupe normal minimal un $p$-groupe abélien élémentaire. Avec les notations de Jordan, un tel groupe est l'ensemble des substitutions de la forme $(x, y, z, ...; x+\alpha, x+\beta, x+\gamma, ...)$. Il forme un produit direct de groupes cycliques et correspond au groupe multiplicatif du corps fini $F_q$. $G$ opère alors sur son groupe normal minimal comme un sous-groupe de $Gl_m(F_q)$.\footnote{Dans les années 1890, ce théorème a été reformulé comme énonçant que le groupe des automorphismes d'un $p$-groupe abélien élémentaire est un groupe linéaire. Voir \cite{Brechenmacher:2011}.}  Le groupe linéaire est ainsi défini comme le normalisateur d'un $p$-groupe abélien élémentaire (lui même un espace vectoriel sur $F_p$). 

La réduction des groupes primitifs à celle des groupes linéaires joue un rôle clé dans la pièce maitresse du \textit{Traité} de 1870, à savoir le Livre IV consacré à la recherche des sous-groupes résolubles maximaux du groupe symétrique. Ce livre met en effet en \oe{}uvre une chaine de réductions de la classe des groupes résolubles de la plus générale à la plus simple : groupes transitifs, primitifs, linéaires, symplectiques etc. 

En termes actuels, cette approche consiste à dévisser un groupe à partir de son socle. Si $F$ est un sous-groupe normal minimal de $G$, le centralisateur $C_G(F)$ est normal (il s'agit du socle de $G$) et $G/ C_G (F)$ est isomorphe à un sous-groupe de $Aut(F)$. Dans le cas étudié par Jordan, $F$ est abélien élémentaire, il est donc son propre centralisateur et s'identifie avec le socle de $G$.\cite{Dieudonne:1962} 

Une fois encore, cette formulation moderne ne permet pas de percevoir le rôle joué par la représentation analytique des substitutions dans l'approche de Jordan. Le groupe linéaire est en effet le premier maillon de la chaîne de réduction du \textit{Traité} dont les substitutions ont une représentation analytique. Or ce sont les procédés de réduction élaborés par Jordan sur ces représentations qui permettent de poursuivre plus avant la chaîne de réduction en étudiant les sous-groupes résolubles du groupe linéaire.

\subsection{Indexation des racines et décomposition de la représentation analytique des substitutions}

\begin{quote}
Ce principe du classement des lettres en divers groupes est le même dont Gauss et Abel ont déjà montré la fécondité dans la théorie des équations : il me semble être dans l'essence même de la question, et sert de fondement à toute mon analyse.\cite[p.5]{Jordan1860}
\end{quote}

Dans sa thèse, Jordan a explicitement présenté son procédé de décomposition de la représentation analytique des substitutions comme prenant modèle sur le procédé d'indexation des racines cyclotomiques à l'aide d'une racine primitive $g$ de l'unité modulo $p$.\cite{Neumann2007} 

L'usage par Gauss de ce procédé a été commenté par Poinsot en 1808 en termes d'"ordre" et de "groupes".\cite{Boucard:2011} Il permet de répartir les racines $x_i$ d'une équation cyclotomique en blocs d'imprimitivité correspondant à une décomposition des substitutions du groupe associé en deux formes représentations des cycles selon que ceux-ci permutent les éléments d'un même bloc, forme $(x \ x+1)$ ou opèrent sur les blocs eux-mêmes, forme $(x \  gx)$.

Comme Jordan l'a constaté en concluant sa thèse, sa méthode de réduction généralise un résultat déjà énoncé et partiellement démontré par Galois. Plus précisément, ce dernier s'est appuyé sur ce qu'il a désigné comme la "méthode de décomposition de M. Gauss" pour réduire le problème de la recherche des groupes résolubles transitifs au cas des groupes primitifs.\footnote{La notion de groupe utilisée par Jordan dans les années 1860 correspond en termes actuels à un groupe fini opérant sur un ensemble. Un groupe transitif $G$ correspond ainsi à une action transitive de $G$ sur un ensemble $X$. Il s'agit donc du cas où cette action possède une et une seule orbite: $X$ n'est pas vide et deux éléments quelconques de $X$ peuvent être envoyés l'un sur l'autre par l'action d'un élément du groupe $G$. La notion de transitivité a émergé de la théorie des équations : le groupe des racines d'une équations irréductible est en effet transitif.}
 
Galois a en effet démontré qu'un groupe primitif permutant un nombre premier $p$ de lettres est résoluble si et seulement si ses substitutions ont une "forme linéaire" $(x \ ax+b)$.\footnote{Galois, comme Jordan, utilise en effet le terme de "groupe linéaire" pour désigner les sous-groupes du groupe affine $AGl_n(F_p)$.}  Plus encore, il a énoncé une généralisation de ce résultat au cas d'un groupe résoluble primitif permutant $q=p^n$ lettres, $n$ quelconque. En 1861, Jordan a tenté de démontrer ce résultat dans un supplément ajouté à sa thèse.\footnote{Voir à ce sujet \cite{Neumann2006}.}  Comme la preuve donnée par Galois à son critère de résolubilité des équations primitives de degré $p$, la démonstration donnée par Jordan pour le cas $p^n$ consiste à engendrer la forme analytique des substitutions linéaires par des compositions de translations $(x \ x+1)$ et dilatations $(x \ gx)$. Plus tard, dans son \textit{Traité} de 1870, Jordan a présenté l'"origine du groupe linéaire" comme résultant du problème de la recherche de la représentation analytique du plus grand groupe ayant pour sous-groupe normal un $p$-groupe abélien élémentaire.\cite[p.91]{Jordan1870}. 

A la suite de sa thèse, Jordan a cherché à caractériser les sous-groupes résolubles de $Gl_n(F_p)$. Le mémoire "Sur la résolution algébrique des équations primitives de degré $p^2$", publié en 1868 dans le \textit{Journal de Liouville}, vise ainsi à donner la classification des sous-groupes résolubles de $Gl_2(F_p)$, une problématique que l'auteur présente comme dépassant les explorations menées par Galois dans son "Fragment d'un second mémoire".\cite[p.113]{Jordan1868a}.

C'est dans ce contexte que Jordan a introduit la réduction canonique des substitutions linéaires pour deux variables. Il s'agissait de représenter des classes de groupes par la forme analytique de leurs substitutions en recherchant un changement d'indices tel que l'action de  $S=(x, y ; \ ax+by, \ a'x+b'y)$ se réduise à multiplier $z=mx+ny$ par un simple facteur constant : $Sz \equiv kz \mod (p)$. Ce problème amène à considérer le système de deux équations linéaires $ma+na'\equiv km$ et $mb+nb'\equiv km$ :
\begin{quote}
[...] qui détermineront le rapport  $\frac{m}{n}$  pourvu que l'on prenne pour $k$ une racine de la congruence 
\[
\begin{vmatrix}
a-k & a' \\
b & b'-k
\end{vmatrix}
\equiv 0  \mod (p)
\]
\end{quote}

Jordan a alors distingué trois formes auxquelles ramener $S$ selon que l'équation caractéristique admette deux racines réelles distinctes $\alpha, \beta$, deux racines imaginaires distinctes $\alpha+\beta i$ ou une racine double :

\[
\begin{vmatrix}
z & \alpha z \\
u & \beta u 
\end{vmatrix}
, \begin{vmatrix}
z & (\alpha + \beta i)z \\
u & (\alpha + \beta i^p)u 
\end{vmatrix}
\begin{vmatrix}
z & \alpha z \\
u & \alpha z+ \gamma u 
\end{vmatrix}
\]
Deux ans plus tard, la réduction canonique a donné corps à l'un des principaux théorèmes du chapitre clé du \textit{Traité} consacré aux substitutions linéaires. Son objet est de ramener une substitution linéaire sur $p^n$ indices "à une forme aussi simple que possible" :
\begin{quote}
\textbf{Théorème de réduction canonique} \\ Cette forme simple
\[
\begin{vmatrix}
y_0, z_0, u_0, ..., y'_0, ... & K_0y_0, K_0(z_0+y_0), ... , K_0y'_0, ... \\
y_1, z_1, u_1, ..., y'_1, ... & K_1y_1, K_1(z_1+y_1), ... , K_1y'_1, ... \\
.... & ... \\
v_0, ... & K'_0v_0, ... \\
... & ... \\
\end{vmatrix}
\]
à laquelle on peut ramener la substitution $A$ par un choix d'indice convenable, sera pour nous sa forme canonique.\cite[p.126]{Jordan1870}
\end{quote}

L'énoncé ci-dessus met en évidence le rôle de modèle joué par la décomposition de la forme linéaire $(x \; ax+b)$ en deux formes de représentations analytiques de cycles $(x \ x+a)$ et $(x \ gx)$. Comme dans le cas de la réduction des groupes imprimitifs opérée par Jordan dans sa thèse, la forme canonique est obtenue en articulant la décomposition polynomiale de la congruence caractéristique à une réduction de la forme analytique de la substitution en fonction de l'action de celle-ci sur des blocs des indices.\footnote{En termes actuels, la démonstration de Jordan procède d'une décomposition d'un espace vectoriel de dimension finie sur un corps fini sous l'action d'un opérateur linéaire. La notion sous-jacente d'espace vectoriel véhicule cependant une interprétation géométrique implicite absente de l'approche de Jordan en 1870. Voir \cite[p.178-192]{Brechenmacher:2006a}}  

Dans son \textit{Traité}, Jordan a appliqué son théorème à la détermination de l'ensemble des substitutions commutant à une substitution donnée et à la caractérisation des sous-groupes commutatifs du groupe linéaire.\cite[p.126-154]{Jordan1870}. Cette application de la réduction canonique conclue la partie de l'ouvrage consacrée aux propriétés générales des groupes linéaires et précède l'étude de groupes linéaires spéciaux (orthogogaux et symplectiques notamment). Plus loin, elle intervient dans la classification des groupes résolubles afin de réduire la recherche des groupes résolubles linéaires à celle des groupes résolubles symplectiques.

\subsection{L'ombre d'un théorème de Weierstrass sur la réception de la réduction de Jordan }
Le théorème énoncé dans le cadre de la théorie des substitutions a par la suite donné lieu à une pratique de réduction canonique que Jordan a utilisé dans des cadres théoriques divers : équations différentielles linéaires à coefficients constants (1871-1872), équations de Fuchs (1874 et 1876-1880), formes quadratiques et bilinéaires (1873-1874), intégration algébriques des équations différentielles linéaires des deuxième, troisième et quatrième ordres (1876-1879), ou encore théorie algébrique et théorie arithmétique des formes algébriques (1879-1882). 

Afin de mieux saisir la spécificité de l'appropriation par Poincaré de la réduction de Jordan, intéressons-nous tout d'abord plus généralement à la réception de cette dernière. Comme nous allons le voir, cette réception a été très limitée durant plusieurs décennies. 

Jusqu'au début du XX\up{e} siècle, la réduction de Jordan a circulé principalement dans le cadre précis des équations différentielles linéaires. À l'exception des travaux de Poincaré, elle n'a au contraire pas plus circulé dans le cadre de la théorie des groupes de substitutions dont elle est pourtant issue que dans les autres domaines où l'a mobilisé Jordan comme celui de la théorie des formes. Durant cette période, les références à la réduction canonique se rapportent d'ailleurs le plus souvent à la version donnée par Jordan dans le troisième volume de son \textit{Cours d'analyse de l'École polytechnique} et non au théorème du \textit{Traité des substitutions}.\footnote{Comme l'illustre le traité sur les équations différentielles de Craig,\cite{Craig:1889}, la référence à Jordan intervenait  souvent dans le cadre du problème de Fuchs.} 

\textit{A contrario}, la réduction de Jordan est non seulement absente des traités d'algèbre mais aussi des manuels d'analyse tournés vers l'enseignement comme ceux de Serret\cite{Serret1886} ou d'Hermann Laurent\cite{Laurent1890}. Dans ces deux types d'ouvrages, le problème de la multiplicité des racines caractéristiques est abordé par un raisonnement très classique consistant à rendre les racines inégales par ajout d'un paramètre et à faire ensuite tendre ce paramètre vers 0.\cite[p.151]{Laurent1890}. Un tel raisonnement a été employé depuis le XVIII\up{e} siècle dans le cadre de la culture algébrique portée par l'équation séculaire. Dans le cas symétrique, ce raisonnement peut être rendu conforme aux critères actuels de rigueur. Le théorème de Bolzano-Weierstrass permet en effet un passage à la limite sur l'ensemble des matrices orthogonales, ensemble fermé et borné de $M_n(\mathbb{R})$. Ce raisonnement n'est cependant pas généralisable au cas des systèmes non diagonalisables dans $\mathbb{C}$. 

À la fin du XIX\up{e} siècle, ce problème n'était cependant généralement pas abordé à l'aide du théorème de Jordan. Au contraire, des présentations comme celles de Gaston Darboux\cite[p.405]{Darboux1889} ou de Louis Sauvage\cite{Sauvage1895} s'appuient sur le théorème des diviseurs élémentaires de Weierstrass pour aborder les cas d'occurrence de racines caractéristiques multiples.

Plusieurs éclairages peuvent être donnés sur les raisons de la réception limitée de la réduction de Jordan. Tout d'abord, le cadre des équations différentielles a offert à Jordan l'occasion d'articuler ses travaux à la culture algébrique commune portée par l'équation séculaire. Dès 1871, Jordan a en effet appliqué ses travaux sur les substitutions au problème mécanique des conditions de stabilité des petites oscillations des systèmes mécaniques. Il a notamment démontré qu'un système linéaire non symétrique de $n$ équations différentielles ne peut se réduire en général à $n$ équations indépendantes. L'intégration d'un tel système peut cependant s'effectuer par réduction canonique. Jordan a ainsi transféré son procédé de réduction canonique des groupes de substitutions de $Gl_n(F_p)$ aux substitutions de $M_n(\mathbb{C}$) des changements de variables de systèmes différentiels. 

La résolution de ce problème a donné un écho à la réduction de Jordan dans le cadre de l'étude des systèmes différentiels. Mais elle a  aussi illustré les limites de cette approche. En effet, la réduction canonique nécessite que soient extraites les racines d'une équation algébrique de degré égal au nombre d'équations du système. Elle ne donne donc pas de méthode effective de résolution de systèmes différentiels de plus de cinq équations. Les valeurs de généralité et de simplicité mises en avant par Jordan lors de l'utilisation de son théorème sont ainsi bien éloignées des besoins d'effectivité des astronomes pour l'étude des petites oscillations des planètes sur leurs orbites. 

La tension entre généralité et effectivité donne un autre éclairage sur la réception limitée de l'approche de Jordan. La réduction canonique s'accompagne en effet d'un idéal, que Jordan a mis en avant à plusieurs reprises depuis sa thèse, et qui attribue aux relations entre classes d'objets un statut essentiel, primant sur l'étude des objets eux mêmes. En tant que tel, cet idéal n'était pas limité à Jordan et a été notamment formulé par de nombreux autres auteurs parmi lesquels Poinsot et Galois. Chez Jordan, il implique une mise au premier plan de propriétés générales sur $n$ variables obtenues par des réductions successives des représentations analytiques jusqu'à des expressions les plus simples.

En 1874, Leopold Kronecker a critiqué publiquement l'usage qu'a voulu faire Jordan de sa réduction canonique en théorie des formes. Il a condamné le caractère faussement général des classes de substitutions de $n$ variables étudiées par Jordan par opposition aux objets de la théorie des formes. Il a surtout souligné le caractère formel car non effectif de la réduite présentée par Jordan comme "la plus simple".\cite{Brechenmacher:2007a}. Plus généralement, Kronecker s'est fermement opposé aux prétentions de Jordan de conférer à des approches algébriques une envergure théorique et a mis en avant un idéal selon lequel le "travail algébrique [...] est effectué au service d'autres disciplines mathématiques dont il reçoit ses fins et dont dépendent ses objectifs". Pour Kronecker, les significations sont portées par les objets de l'arithmétique, comme les classes d'équivalences de formes, et non par les   techniques algébriques d'opérations sur les substitutions.

Surtout, Kronecker a opposé à la réduction canonique de Jordan les méthodes de calculs d'invariants par lesquelles il a donné une formulation effective au théorème des diviseurs élémentaires de Weierstrass.  Le rôle joué par ce théorème à la fin du XIX\up{e} siècle est un autre facteur important de la réception limitée de la réduction canonique de Jordan. 

En 1858, Weierstrass a interprété les procédés algébriques (*) véhiculés par l'équation séculaire comme revenant à transformer simultanément un couple de formes quadratiques ($A(x_1,...,x_n), \ x_1^2+...+x_n^2)$ en des sommes des carrés $(s_1x_1^2+ s_2x_2^2+...+ s_nx_n^2, \ x_1^2+...+x_n^2)$. Nous avons vu que le caractère générique de l'expression (*) avait été critiquée car celle-ci semblait perdre toute signification en cas d'occurrences de racines multiples. Weierstrass a cependant montré que, dans le cas symétrique, une racine d'ordre $p$ annule tous les mineurs d'ordre $p-1$ du déterminant caractéristique. L'expression (*) est donc bien définie. 

En 1868, Weierstrass a généralisé son approche au cas non symétrique des couples non singuliers de formes bilinéaires $(A, B)$ (implicitement à coefficients complexes). Il a introduit des invariants, les diviseurs élémentaires, qui caractérisent les classes d'équivalence de tels couples de formes. 

Les diviseurs élémentaires sont définis par la suite de facteurs linéaires intervenant dans le déterminant caractéristique $[A,B]=|A+sB|$ et par les suites décroissantes d'exposants de chaque facteur linéaire des mineurs successifs de $[A,B]$. Soit $a+sb$ un facteur linéaire d'exposant $l$ de $[A,B]$ et soit $l(x)$ le plus grand exposant pour lequel tous les mineurs d'ordre $n-x$ contiennent le facteur $(a+sb)^{l(x)}$. Alors, soit $l(x-1)>l(x)$ soit $l(x)=0$. Considérons à présent $r$ le plus petit nombre tel que $l(r)=0$, alors  la suite $l, l', l'', ..., l^{(r-1)}$ est une suite strictement positive décroissante. Avec ces notations, et si l'on note $e=l-l'$, $e'=l'-l''$, ... ,$e^{(r-1)}$=$l^{(r-1)}$, les relations de divisibilité entre le déterminant et ses mineurs successifs sont traduites par la décomposition suivante du facteur $(a+sb)^l$ de $[A,B]$ : 
\[
(a+sb)^l=(a+sb)^e(a+sb)^{e'}...(a+sb)^{e^{r-1}}
\]
L'ensemble des termes de cette décomposition, étendue à tous les facteurs linéaires du déterminant 
\[
[A,B]=C(a_1+sb_1)^{l_1}(a_2+sb_2)^{l_2}...(a_p+sb_p)^{l_p}
\]
est appelé ensemble des diviseurs élémentaires de $[A,B]$. 

Deux couples non singuliers de formes bilinéaires $A+sB$ et $A'+sB'$ sont transformables l'un en l'autre par changements de variables si et seulement si leurs ensembles de diviseurs élémentaires sont identiques.\cite[p.21]{Weierstrass1868} Dans ce cas, les deux couples peuvent tous deux être transformés en une forme normale $(\Phi, \Psi)$, équivalente à la forme canonique donnée par Jordan aux substitutions sur des corps finis:
\begin{center}
$
\Phi=\sum (X_\lambda Y_\lambda)_{e_\lambda} ; \
\Psi = \sum c_\lambda(X_\lambda Y_\lambda)_{e_\lambda} + \sum(X_\lambda Y_\lambda)_{e_{\lambda-1}}
$ \\
\end{center}

En 1874, Kronecker a reformulé le théorème de Weierstrass en introduisant des invariants déterminés par des procédés effectifs (en termes actuels, il s'agit des facteurs invariants d'une matrice sur un anneau principal). Dans la suite $\Delta_r(s)$, $\Delta_{r-1}(s)$, ..., $\Delta_1(s)$ des $p.g.c.d.$ des mineurs successifs de $[A,B]$, chaque polynôme est divisible par le précédent. Les quotients correspondants $i_1(s)$, $i_2(s)$, ..., $i_r(s)$ donnent des invariants caractérisant le couple $A+sB$. La décomposition de ces polynômes en facteurs linéaires sur un corps $K$ algébriquement clos donne les diviseurs élémentaires de Weierstrass.

Quelques années plus tard, Frobenius a publié une synthèse unifiant les théories des substitutions linéaires et des formes bilinéaires et quadratiques.\cite{Frobenius:1877} \cite{Frobenius:1879}. Cette théorie est basée sur les notions arithmétiques de classe d'équivalence et d'invariants comme l'a voulu Kronecker. Elle s'organise autour du théorème des diviseurs élémentaires. La réduction canonique de Jordan n'y est qu'un corollaire pour le cas particulier où serait permise l'utilisation de "nombres imaginaires de Galois",\cite[p.544]{Frobenius:1879}, c'est à dire l'extraction de racines d'équations algébriques $mod. p$ (dans des corps finis, en termes actuels). 

Nous avons évoqué au paragraphe précédent le lien entre la pratique de réduction de Jordan et les procédés d'indexation des racines de l'unité ou des imaginaires de Galois. Ces procédés permettent en effet la décomposition de la représentation analytique des substitutions linéaires. Cette décomposition a elle même servi de modèle à la décomposition des groupes linéaires. Cette approche est très spécifique aux travaux de Jordan sur les substitutions. Bien que ce dernier soit parvenu à la transférer à d'autres thématiques, comme celles portées l'équation séculaire, il s'est trouvé en concurrence avec le théorème des diviseurs élémentaires. Or ce théorème permet non seulement de résoudre les mêmes problèmes que celui de Jordan mais les procédés de calculs d'invariants sur lesquels il s'appuie s'inscrivent bien plus naturellement que ceux de Jordan dans le cadre traditionnel de l'équation séculaire.

La synthèse théorique organisée par Frobenius autour du théorème des diviseurs élémentaires s'est imposée à une échelle large jusque dans les années 1920. Elle a notamment joué un rôle important dans la grande majorité des travaux consacrés aux groupes discrets et continus ou aux algèbres associatives. L'utilisation de la réduction canonique de Jordan par Poincaré dans ce cadre présente donc une forte spécificité.

\subsection{Réduction canonique et équations différentielles linéaires (1871-1878)}

Comme nous l'avons vu, la pratique de Jordan de réduction des représentations analytiques des substitutions présente des caractéristiques mettant en jeu le statut de l'algèbre, les rapports entre relations générales et objets particuliers, ainsi que des valeurs épistémiques de simplicité et de généralité. Il s'agit à présent d'en suivre les évolutions dans les travaux de Jordan des années 1870 afin de jeter un éclairage sur sa circulation presque souterraine dans les premiers travaux de Poincaré.

\subsubsection{Réduction de Jordan et équations de Fuchs}

C'est dans le cadre de problèmes associés à l'équation séculaire que Meyer Hamburger a, pour la première fois, rapproché en 1873 le théorème de Jordan du théorème des diviseurs élémentaires de Weierstrass. Il a montré comment l'approche de Jordan sur les équations différentielles linéaires à coefficients constants peut être généralisée aux équations dites de Fuchs à coefficients méromorphes sur une région simplement connexe de $\mathbb{C}$ :
\[
\frac{d^m u}{dz^m} +f_1(z) \frac{d^{m-1}u}{dz^{m-1}} +...+f_m (z)u=0
\]

Ces équations posent notamment des problèmes de prolongements analytiques de solutions locales au voisinage de singularités. Dans la tradition des travaux menés depuis les années 1850 par des auteurs comme Cauchy, Puiseux, Hermite ou Riemann, de tels problèmes ont été abordés en terme de ce que Jordan a dénommé des "groupes de monodromie" dans son \textit{Traité}.\footnote{Un ouvrage de Gray\cite{Gray:2000} développe une histoire sur le temps long du XIX\up{e} siècle de ces problématiques et met notamment en évidence le rôle de modèle joué par l'équation hypergéométrique pour l'objectif d'une compréhension globale des relations entre intégrales définies localement par des séries entières.} Précisons que tous ces travaux concernent des cas où les points singuliers sont dénombrables : ils envisagent des solutions en un ensemble de points réguliers constituant le voisinage d'un unique point singulier $a$. 

Pour des valeurs initiales fixées, soit $u_i(z)$ un système fondamental de $n$ solutions de l'équation différentielle en un point régulier ; si la variable $z$ décrit un lacet autour de $a$, les nouvelles valeurs obtenues pour $u_i(z)$ s'obtiennent à partir des valeurs initiales par une transformation linéaire $A$ de déterminant non nul.\footnote{Cette approche a été développée par Puiseux en 1850 dans le cadre de l'analyse complexe de Cauchy et des travaux d'Hermite sur ce que l'on désignerait aujourd'hui comme la théorie de Galois différentiel en lien avec les fonctions abéliennes. Elle a fait l'objet de la seconde thèse de Jordan. Ce dernier a étudié dans les années 1860 les groupes de monodromie des équations différentielles associées aux fonctions elliptiques et abéliennes, notamment en relation avec les travaux d'Hermite sur l'équation modulaire d'ordre 5. Voir \cite{Tannery1875} pour une présentation du problème de Fuchs et de sa réciproque.} 

En 1866, Fuchs a montré que si l'équation caractéristique de $A$ n'a que des racines simples, les solutions prennent localement la forme $(z-a)\Phi(z-a)$ avec $\Phi$ holomorphe au voisinage de $a$. En cas d'occurrence d'une racine caractéristique multiple d'ordre $k$, Fuchs a exprimé le système d'intégrales associées sous une forme triangulaire. Hamburger s'est appuyé à la fois sur la réduction canonique de Jordan et sur le théorème des diviseurs élémentaires afin de caractériser plus finement les solutions en cas d'occurrence de racines multiples. D'un côté, le théorème de Weierstrass lui a permis de surmonter une difficulté qu'avait passé sous silence Jordan et consistant à distinguer entre l'ordre de multiplicité $\mu$ d'une racine caractéristique et le nombre $\nu$ de solutions indépendantes associées. Comme le formule Hamburger, les $\mu$ fonctions associées à une racine $\omega$ se "répartissent en différents groupes et non pas toutes ensembles".\cite[p.114]{Hamburger1873}.\footnote{Dans sa démonstration, Jordan ne donne en effet aucune manière de déterminer la dimension $\nu$ de chaque espace caractéristique. Au contraire, chez Weierstrass l'invariant $\nu$ est un diviseur élémentaire bien déterminé qui permet de distinguer entre le cas diagonalisable et non diagonalisable.\cite[p.178-192]{Brechenmacher:2006a}}  D'un autre côté, la réduction de Jordan a permis d'exprimer les solutions associées à une racine multiple $\omega_1$ en fonction de $t^r$ (avec $t=z-a$ et $r=\frac{log\omega_1}{2\pi i}$), de $logt$ et de fonctions méromorphes $M_i$, $N_i$  : 
\[
y_o=t^rM_0 ; \;
y_1=t^rM_1(logt+N_1) ; \;
... ; \;
y_k=t^r(M_klog^kt+N_klog^{k-1}t+...) \;
\]
\subsubsection{L'application par Jordan de la théorie des substitutions aux équations différentielles}

Le problème de l'intégration algébrique des équations différentielles linéaires étudiées par Fuchs entre 1865 et 1868 a été repris par de nombreux mathématiciens dans les années 1870. Dans ce contexte, Jordan a élaboré entre 1874 et 1880 une approche spécifique basée sur ses travaux antérieurs sur les groupes de substitutions. Ces travaux ont donné lieu à un mémoire de synthèse publié en 1878 dans le \textit{Journal de Crelle} qui a joué un rôle de transmission important entre les travaux de Jordan et ceux de Poincaré. Ils témoignent d'une généralisation à de nouvelles situations des procédés d'indexation et de représentation analytique des substitutions. Comme nous allons le voir, Jordan a en effet transféré sa pratique de décomposition des groupes linéaires en sous-groupes au cadre des équations différentielles, c'est à dire à ce que nous désignerions aujourd'hui comme la notion de sous-espace vectoriel caractéristique pour l'action d'un opérateur. 

Publié en 1874, le "Mémoire sur une application de la théorie des substitutions à l'étude des équations différentielles linéaires" proclame que le groupe obtenu  en faisant varier $z$ "de toutes les manières possibles de manière à envelopper successivement les divers points singuliers", caractérise "dans ce qu'il a d'essentiel le type de l'équation différentielle qui lui donne naissance, et reflète ses principales propriétés".\cite[p.102]{Jordan1874a} Pour qu'une équation différentielle soit susceptible d'une intégration algébrique, c'est-à-dire, que toutes ses intégrales satisfassent à des équations algébriques ayant pour coefficients des fonctions holomorphes de $z$, le groupe linéaire associé doit être fini. Sur le modèle de l'articulation entre groupes finis et équations algébriques, Jordan en a déduit trois types de problèmes :

\begin{quote}
On voit par ces exemples que dans la théorie des équations différentielles linéaires, comme dans celle des équations algébriques, on aura trois catégories de problèmes à résoudre : \\
1.	L'équation étant donnée, déterminer son groupe. Cette question est du ressort du calcul intégral ; \\
2.	Déterminer les conditions auxquelles le groupe doit satisfaire pour que l'équation donnée jouisse de telle ou telle propriété ; \\
3.	Le groupe étant connu, vérifier s'il satisfait ou non aux conditions requises. Cette dernière question ne dépend plus que de la théorie des substitutions.\cite[p.103]{Jordan1874a}
\end{quote}

Comme de nombreux autres travaux contemporains sur les équations de Fuchs, le mémoire de Jordan se base sur un transfert analogique entre équations algébriques et différentielles. À la même époque, Frobenius s'est lui aussi appuyé sur une analogie avec les équations polynomiales pour développer la notion d'équation différentielle irréductible.\cite[p.56-61]{Gray:2000}\footnote{Pour un panorama des travaux de Thomé et Frobenius sur les intégrales régulières, voir \cite{Floquet:1879}.} De son côté, Jordan a davantage pris modèle sur les problèmes de réduction du degré des équations algébriques provenant de la division des périodes des fonctions elliptiques et abéliennes qu'il a étudiée entre 1868 et 1870. Le mémoire de 1874 vise ainsi à "reconnaître si l'équation différentielle linéaire qui a pour groupe $G$ est satisfaite par les intégrales d'équations analogues d'un ordre inférieur à $n$, et [à] déterminer les groupes de ces équations réduites".  

Ce transfert analogique entre équations algébriques et différentielles est notamment porté par des extensions de procédures opératoires comme la réduction canonique. Comme le formule Jordan, "la notion du faisceau dans la théorie des fonctions linéaires est analogue à celle du groupe dans la théorie des substitutions". Le terme de faisceau, qui avait jusqu'à présent été utilisé pour désigner un sous-groupe de permutations, est ainsi transféré à des systèmes de fonctions linéaires dont "toute combinaison linéaire fait elle-même partie de ce système" et qui contiennent donc "un certain nombre de fonctions linéairement distinctes, en fonction linéaire desquelles on pourra exprimer toutes les autres". Il s'agit, en termes actuels, de sous-espaces vectoriels caractéristiques.\cite[p.106]{Jordan1874a} 

Jordan a alors redéfini de la manière suivante la notion d'équation différentielle irréductible introduite par Frobenius : si aucun faisceau ne peut être formé hormis celui correspondant à l'ensemble des combinaisons linéaires des intégrales indépendantes $y_1, ... ,y_n$, le "groupe $G$ sera dit primaire, et l'équation différentielle correspondante sera irréductible" ; dans le cas contraire, à chacun des  faisceaux correspond une équation différentielle réduite.

\subsubsection{Intégration algébrique et classification des groupes linéaires finis (1876-1879)}

Par la suite, les travaux de Jordan ont visé une classification des types d'équations irréductibles intégrables algébriquement. Cette réorientation semble avoir été impulsée par de nouveaux travaux de Fuchs. En 1875, ce dernier s'est appuyé sur les méthodes de la théorie des invariants et covariants des formes binaires afin d'étudier les équations algébriques liant deux intégrales indépendantes d'une équation différentielle linéaire du second ordre.\cite[p.77-81]{Gray:2000} \footnote{Voir le résumé en français donné par Fuchs\cite{Fuchs:1876c} dans le \textit{journal de Liouville} sous la forme d'une lettre à Hermite.}  A cette époque, Jordan a lui-même consacré plusieurs travaux à la théorie des covariants dans le sillage du résultat de Gordan sur la possibilité d'exprimer les covariants d'un système de formes binaires en fonction entière d'un nombre limité de covariants indépendants. Comme en 1874, il a cependant proposé de traiter les questions abordées par Fuchs "par une méthode toute différente fondée sur la théorie des substitutions". 

Un groupe fini $G$ étant associé à une équation différentielle $E$, la détermination d'une équation algébrique dont $G$ est le groupe de Galois (c'est à dire ce que nous appellerions aujourd'hui le problème inverse de la théorie de Galois) donne une relation algébrique entre les intégrales de $E$, c'est à dire une intégration algébrique de cette dernière équation.\cite[p.607]{Jordan1876a} Autrement dit :
\begin{quote}
Il y a donc identité entre les deux questions suivantes : \\
1.	Énumérer les divers types d'équations différentielles linéaires d'ordre $m$ dont toutes les intégrales soient algébriques. \\
2.	Construire les divers groupes d'ordre fini que contient le groupe linéaire à $m$ variables.\cite[p.605]{Jordan1878}
\end{quote}

Jordan s'est alors appuyé sur les théorèmes de Sylow pour établir une répartition en trois classes des sous-groupes finis de $Gl_2(\mathbb{C})$. Ces classes sont identifiées par la représentation analytique des formes canoniques de leurs substitutions génératrices.\cite[p.607]{Jordan1876a} \footnote{Dans sa première classification, Jordan ne mentionne ni le sous-groupe simple d'ordre 168 ni le sous-groupe simple A6 d'ordre 360 de$Sl_3(\mathbb{C})$. Voir \cite[p.XXV]{Dieudonne:1962}}  

Comme nous allons le voir plus loin, ce résultat a suscité une revendication de priorité de Klein. Face aux objections de ce dernier, Jordan a cependant défendu l'avantage et la nouveauté de sa propre méthode de classification par réduction canonique\cite[p.1035]{Jordan1876b} en raison de sa capacité à généraliser la classification des sous groupes finis de $Gl_2(\mathbb{C})$ aux cas de 3 variables \cite{Jordan1877a} et 4 variables \cite{Jordan1879b}. 

\subsubsection{Le théorème finitude de l'index des sous groupes de $Gl_n(\mathbb{C})$}

Plus encore, la réduction canonique a sous-tendu l'énoncé d'un résultat pour le cas de $n$ variables : le théorème de finitude de l'index d'un groupe linéaire fini par rapport à un sous-groupe abélien normal. \footnote{Il y a une infinité de sous-groupes finis de $Gl_n(\mathbb{C}$) mais il existe une fonction $\phi(n)$ telle que tout groupe fini $G$ de matrices d'ordre $n$ contient un sous-groupe normal $H$ qui est le conjugué d'un groupe de matrices diagonales et tel que l'index $(G: H)$ soit inférieur à $\phi(n)$. Le groupe quotient $G/H$ appartient donc à un système fini de groupes à isomorphisme près. Voir \cite[p.XXIII]{Dieudonne:1962}}  Ce théorème est le résultat central du mémoire de synthèse que Jordan a publié dans le \textit{Journal de Crelle} en 1878.
\begin{quote}
\textbf{Théorème de finitude de Jordan}. Si un groupe $G$ est formé d'un nombre fini de substitutions linéaires à $n$ variables, il contiendra un autre groupe $H$ dont les substitutions seront de la forme simple
\begin{center}
$x_1, ..., x_n \ ; \ a_1x_1, ..., a_nx_n$
\end{center}
et permutable à toutes les substitutions de $G$. L'ordre $g$ de $G$ sera égal à $kh$, $h$ étant l'ordre de $H$, et $k$ un entier inférieur à une limite fixe, assignable a priori pour toute valeur de $n$. \\
Ou, en d'autres termes : Si une équation différentielle linéaire d'ordre $n$
\[
(E)  f(z)u+f_1(z)u'+...+f_{n-1}(z)u^{(n-1)}+u^{(n)}=0
\]
a ses intégrales algébriques, elle admettra $n$ intégrales particulières $x_1$, ..., $x_n$ racines d'équations binômes, dont les seconds membres seront des fonctions rationnelles de $z$ et d'une racine $\gamma$ d'une équation irréductible
\[
F(z, \gamma)=0
\]
dont le degré $k$ sera inférieur à une limite fixe. \\
Ou bien encore, en empruntant le langage de M. Fuchs : Le degré des formes primitives construites avec les intégrales de l'équation $(E)$ sera limité.\cite[p.1036]{Jordan1878}
\end{quote}

Nous allons voir plus loin que ce résultat a, quelques années plus tard, servi de modèle aux travaux de Jordan et Poincaré en théorie des formes. La mobilisation par Jordan de sa réduction canonique dans la démonstration de son théorème\cite[p.90]{Jordan1878} a notamment joué un rôle important pour l'appropriation ultérieure par Poincaré de ce procédé algébrique. Cette preuve s'appuie notamment sur le lemme suivant énonçant que toute substitution "périodique" (d'ordre finie) est diagonalisable :
\begin{quote}
[Si deux racines de l'équation caractéristique en $s$ sont distinctes] $S = (u_1, u_2 ; \ \alpha u_1+ \beta u_2, \gamma u_1+ \delta u_2 )$ [...] sera réduite à la forme canonique $(x, y  ; \ ax, by)$. Nous dirons dans ce cas que $S$ est une substitution de première espèce.
Si les deux racines de l'équation en $s$ se confondent en une seule, $a$, soient $x$ une fonction linéaire que $S$ multiplie par $a$ [...] $S$ prendra la forme $(x, y \ ; \ ax, a(y+ \lambda x)).$ Nous dirons que $S$ est de seconde espèce si $\lambda=0$ ; de troisième espèce, si $\lambda \ne 0$. \\
Soit maintenant $G$ un groupe formé d'un nombre limité de substitutions linéaires. Il ne peut contenir aucune substitution $S$ de troisième espèce. Car il contiendrait ses puissances, qui ont pour formule générale $S^m = (x,  y ; \ a^mx, a^m(y+m\lambda x))$ et sont évidemment en nombre illimité. Quant aux substitutions de première espèce, leurs puissances ont pour formule $S^m = (x, y ; \ a^mx, b^my)$ et seront en nombre limité, lorsque $a$ et $b$ seront des racines de l'unité.\cite[p.93]{Jordan1878} 
\end{quote}

Malgré les références fréquentes dont le théorème de finitude a fait l'objet à la fin du XIX\up{e} siècle, la pratique de réduction canonique sur laquelle ce résultat est basé n'a été que très rarement reprise. Bien que les résultats de  Jordan aient été largement commentés, notamment par Klein, Brioschi, Gordan et Richard Dedekind.\cite{Brechenmacher:2011}, ils ont fait l'objet de reformulations dans d'autres approches des équations différentielles. 

Comme nous l'avons déjà évoqué, Klein a donné en 1875 - soit un an avant Jordan - une classification des sous-groupes finis de $Gl_2(\mathbb{C}$) en cinq classes correspondant aux groupes de symétries des polyèdres réguliers.\cite{Klein1875} \footnote{De telles problématiques étaient notamment inspirées de l'approche géométrique de Clebsch sur les formes binaires, mais aussi des travaux de Schwarz \cite{Schwarz1872} sur les conditions auxquelles une série hypergéométrique de Gauss est une fonction algébrique. Ce dernier a défini une application conforme d'un triangle dont les côtés sont des arcs de cercles dans le disque unité. A la suite des travaux de Klein et Poincaré, la fonction de Schwarz a été considérée comme une fonction automorphe. A propos de la méthode de Schwarz et sur la manière dont ce dernier en a déduit une démonstration du théorème d'application conforme de Riemann et du problème de Dirichlet, voir \cite{Tazzioli1994}} Cette classification vient cependant répondre à des problèmes très différents de ceux abordés par Jordan en 1876. Elle s'inscrit en effet dans le cadre de travaux géométriques sur les transformations linéaires laissant invariante une forme binaire. En outre, bien que Klein se soit appuyé sur des travaux antérieurs de Jordan sur les groupes de mouvements des polyèdres,\cite{Jordan1868b}, il a mis en \oe{}uvre des méthodes très différentes basées sur les covariants de Gordan et l'approche géométrique de Clebsch. 

Les travaux de Jordan sur la classification des équations à intégrales algébriques ont cependant amené Klein à lier ses travaux de 1875 aux équations différentielles. Ce dernier a notamment explicité les équations algébriques vérifiées par les intégrales des cinq types d'équations différentielles du second ordre à intégrales algébriques.\cite{Klein1877} De son côté, Gordan\cite{Gordan:1877} a reformulé la classification donnée par Klein aux groupes linéaires finis de deux variables en évacuant toute considération géométrique au profit d'un approche purement algébrique sur les covariants des formes binaires. Enfin, Brioschi \cite{Brioschi:1877} et Dedekind\cite{Dedekind:1877} ont mis en avant les relations entre les travaux de Klein, de nouveaux travaux de Fuchs\cite{Fuchs:1877} et les travaux d'Hermite et Kronecker des années 1850 sur la résolution de l'équation générale du cinquième degré par les fonctions elliptiques. 

Ces problématiques sont aussi celles pour lesquelles Poincaré et Klein ont introduit quelques années plus tard les fonctions fuchsiennes et modulaires, c'est à dire des fonctions complexes invariantes par certains groupes linéaires comme $PSl_2(\mathbb{R}$) ou $PSl_2(\mathbb{Z}/p\mathbb{Z}$). 

C'est en effet dans ce contexte que Poincaré s'était approprié les travaux de Jordan. Comme nous allons le voir à présent, cette appropriation permet d'éclairer la spécificité de l'approche développée par Poincaré, notamment par rapport à celle de Klein. Elle permet également de mettre en évidence la cohérence des premiers travaux de Poincaré malgré la diversité thématique apparente de ces derniers.

\section{De Jordan à Poincaré via Hermite}
Nous allons étudier à présent la manière par laquelle Poincaré s'est approprié les travaux de Jordan. Nous allons voir que cette appropriation s'est effectuée par la médiation d'un héritage des travaux d'Hermite et permet de jeter un nouvel éclairage sur l'émergence de la célèbre théorie des fonctions fuchsiennes. 

\subsection{Fonctions elliptiques et équations différentielles linéaires}

Rappelons tout d'abord que les fonctions fuchsiennes permettent de généraliser aux équations différentielles linéaires à coefficients algébriques le rôle joué par les fonctions elliptiques dans l'intégration de différentielles algébriques. 

Par exemple, la fonction elliptique de premier ordre de module $k^2$ est définie comme la fonction réciproque $\phi(u,k)$ de l'intégrale :
\[
u(\phi, k) = \int_{0}^{\phi} \frac{dx}{\sqrt{(1-x^2)(1-k^2x^2)}}
\]
Les fonctions elliptiques peuvent être étendues en des fonctions méromorphes sur $\mathbb{C}$ et périodiques sur un réseau tel que, dans l'exemple ci-dessus, le réseau de périodes $K=u(1,k)$ et $K'=u(\frac{1}{k},k)$ (et de pôles $\frac{nK}{4}$). En tant que fonctions du module $k^2$, ces périodes vérifient elles-mêmes une équation différentielle linéaire. Cette équation, dite équation de Legendre, est un cas particulier d'équation hypergéométrique et donc d'équation de Fuchs :
\[
(**)  (1-k^2)\frac{d^2y}{d^2k}+\frac{1-3k^2}{k} \frac{dy}{dk} - y=0
\]
Au XIX\up{e} siècle, les tentatives de généralisation des fonctions elliptiques à des fonctions - dites abéliennes ou hyperelliptiques - obtenues par inversions d'intégrales d'équations différentielles plus générales se sont heurtées au caractère local de telles fonctions. Des problèmes de raccordements des valeurs multiples de ces fonctions ont été au c\oe{}ur des principaux développements de l'analyse complexe.\footnote{Voir \cite[p.449]{BriotetBouquet:1859}} 

Sur les conseils d'Hermite, Fuchs a considéré en 1877 le rapport $\omega=\frac{K}{K'}$, des périodes de la fonction elliptique de premier ordre comme une fonction analytique du module $k$.  La réciproque,  $k= f(\omega)$, est une fonction monodrome de $\omega=x+iy$ pour tout $y$ positif.\cite[p.101-104]{Gray:2000} La monodromie est cependant perdue dans la situation analogue consistant à considérer les rapports des périodes $J$ et $J'$ d'intégrales elliptiques du second ordre :
\[
J = \int_{0}^{1}\frac{k^2x^2dx}{\sqrt{(1-x^2)(1-k^2x^2)}} \ ; \ J'=\int_{0}^{\frac{1}{k}}\frac{k^2x^2dx}{\sqrt{(1-x^2)(1-k^2x^2)}}
\]
et satisfaisant à l'équation différentielle
\[
(***) \ (1-k^2 )  \frac{d^2 y}{d^2 k}+\frac{1-k^2}{k} \frac{dy}{dk}+y=0
\]
Mais comme l'avait déjà montré Hermite, les fonctions des modules sont dans les deux cas invariantes par des substitutions homographiques de déterminant 1 (dites unimodulaires).\footnote{La fonction modulaire est invariante par $Sl_2(\mathbb{Z})$ et l'équation modulaire peut ainsi être introduite par l'équation de transformation des intégrales elliptiques par de telles substitutions.}  Étant donné un système fondamental $K$, $K'$ de solutions de (**) ou (***), le parcours de la variable $k$ sur un lacet autour d'un point singulier donne deux nouvelles solutions, pouvant s'exprimer sous la forme $a_1K +b_1K'$ et $a_2K +b_2K'$, et à partir desquelles Fuchs a formé le quotient 
\[
H=\frac{a_1 K +b_1 K'}{a_2 K +b_2 K'}
\] 

Des prolongements analytiques de $\frac{J}{J'}$ sont ainsi possibles par l'action des substitutions homographiques unimodulaires sur le demi-plan supérieur. Cette approche est à la base des théories développées par Klein et Poincaré au début des années 1880. 

Le problème peut cependant être aussi conçu indépendamment de la notion de groupe. Comme nous allons le voir, avant qu'il ne s'approprie les travaux de Jordan sur les substitutions, Poincaré a initialement abordé le problème de Fuchs selon une approche arithmétique en terme de classes d'équivalence de formes.

\subsection{Fonctions modulaires et classes d'équivalences arithmétiques }

La propriété d'invariance des fonctions modulaires étudiées par Fuchs en 1877 a immédiatement été interprétée par Dedekind dans un cadre arithmétique. Ce dernier a rapproché cette propriété des travaux d'Hermite et Kronecker des années 1850 sur la résolution de l'équation du cinquième degré par les équations modulaires des fonctions elliptiques.\footnote{Les travaux d'Hermite prenaient place dans le cadre d'une approche des nombres algébriques par les formes quadratiques et fractions continues ;\cite{Goldstein:2007} ; ceux de Kronecker visaient une représentation concrète des nombres idéaux introduits par Kummer et s'appuyaient sur la multiplication complexe des fonctions elliptiques.\cite{Petri2004}}  Il en a cependant donné une formulation dans le cadre de sa propre théorie des corps et des idéaux d'entiers complexes. Décrivons rapidement l'approche de Dedekind afin de mettre en évidence la différence de l'approche arithmétique suivie par Poincaré.

Deux points du demi-plan supérieur peuvent être considérés comme "congruents", pour reprendre les termes de Dedekind, si leurs affixes sont liées par une substitution unimodulaire. Le problème de Fuchs se présente alors comme celui de la définition d'une fonction monodrome sur des classes de congruences de nombres complexes.

Dedekind a abordé ce problème pour le cas des équations hypergéométriques. Soit $F$ le domaine fondamental des points $\omega$ \{$-1/2<Re(\omega)<1/2$ ; $N(\omega)>1$ \} du demi-plan supérieur, situés strictement entre les droites d'équation $Re(\omega)=-1/2$ et $Re(\omega)=1/2$ et strictement à l'extérieur du cercle unité, auquel on ajoute les points d'abscisse $Re(\omega)=-1/2$ ainsi que les points du cercle d'abscisse $-1/2 \leqslant x \leqslant 0$. Ce domaine fournit un système complet de représentants des classes de congruence : tout point du demi plan supérieur est congruent à un unique point du domaine fondamental. La "fonction de valence" définie par Dedekind sur le domaine fondamental prend ainsi ses valeurs sur les classes de nombres complexes du demi plan supérieur. 

Dès 1878, l'approche de Dedekind a été envisagée par Klein sous l'angle très différent de l'action du groupe $PSl_2(\mathbb{Z}/p\mathbb{Z}$) sur le pavage obtenu à partir du domaine $F$.\footnote{Dedekind a également introduit à cette occasion sa fonction $\eta$ \cite[p.107-115]{Gray:2000}}  Quant à Poincaré, sa méconnaissance de l'approche de Dedekind est attestée par certains échanges avec Klein à l'époque de la controverse suscitée par le choix du nom des fonctions fuchsiennes.\cite[p.184-188]{Gray:2000} 
 
Comme nous allons le voir, Poincaré s'est appuyé sur la théorie de la réduction des formes d'Hermite pour aborder le problème de l'extension de la fonction modulaire considérée par Fuchs. Cette théorie donne en effet une alternative à l'approche de Dedekind sur les nombres algébriques comme à celle de Klein sur les groupes de transformations. Elle a d'ailleurs été développée par Hermite dans le cadre du programme de ce dernier visant à caractériser les nombres algébriques via la réduction des formes et les fractions continues.\cite{Goldstein:2007}.

 C'est dans ce contexte que Poincaré s'est progressivement approprié certaines approches de Jordan sur les substitutions. Un premier point de rencontre a eu lieu lors des travaux complémentaires menés par les deux mathématiciens sur les formes algébriques entre 1880 et 1881. Il a été l'occasion pour Poincaré de se familiariser avec les méthodes utilisées par Jordan pour étudier les formes quadratiques invariantes par des groupes linéaires. Surtout, les travaux de Jordan ont progressivement amené Poincaré à modifier son approche arithmétique initiale pour envisager des fonctions invariantes par l'action du groupe $PSl_2(\mathbb{R})$.
 
La considération par Poincaré de groupes linéaires de substitutions à coefficients réels (ou complexes) et non seulement entiers comme chez Klein ou Dedekind témoigne de l'héritage mêlé des travaux de Jordan et d'Hermite. L'utilisation de variables continues en théorie des nombres, avec notamment la réduction continue des formes quadratiques, avait en effet amené ce dernier à introduire une distinction entre les classes d'équivalences arithmétique et algébrique, selon que les substitutions opérant sur les formes quadratiques ont leurs coefficients entiers ou réels. Dans le cadre de l'approche de Jordan, ces classes d'équivalences peuvent être associées à des groupes linéaires laissant stable une fonction, c'est d'ailleurs l'approche qui avait déjà été suivie dans la thèse ce dernier comme nous l'avons vu au paragraphe 2.
 
 \subsection{La chronologie des publications de Poincaré entre 1879 et 1881}
 
 Il nous faut être attentif à l'évolution des approches mises en \oe{}uvre par Poincaré au tournant des années 1870-1880 afin de saisir le rôle joué par Jordan dans le basculement entre une approche arithmétique centrée sur les classes d'équivalences de formes et une approche algébrique basée sur la notion de groupe.
 
Rappelons que la théorie des fonctions fuchsiennes a émergé de travaux menés par Poincaré en vue du Grand prix des sciences mathématiques de 1880 : "Perfectionner en quelque point important la théorie des équations différentielles linéaires à une seule variable indépendante". Afin de comprendre comment les publications de Poincaré s'articulent les unes aux autres, il nous faut à présent suivre la chronologie fine de ces publications dans l'intervalle séparant le 10 mars 1879, date de mise au concours du Grand prix, et le 1er juin 1880, date limite de réception des mémoires. 
 
Poincaré a tout d'abord consacré deux notes aux formes quadratiques.\cite{Poincar1879b} \cite{Poincar1879a} qui ont donné lieu au mémoire "Sur un nouveau mode de représentation géométrique des formes quadratiques définies ou indéfinies", publié dans le \textit{Journal de l'École polytechnique}.\cite{Poincar1880a} 
 
 La note "Sur les courbes définies par les équations différentielles" $\frac{dx}{X}=\frac{dy}{Y}$ ($X, Y$ polynômes réels) \cite{Poincar1880b},\footnote{Christian Gilain\cite[p.240]{Gilain:1991} donne une analyse détaillée du rôle joué par la géométrie dans la conception pluraliste de Poincaré sur la théorie qualitative des équations différentielles.} publiée le 22 mars 1880, a quant à elle constitué la base d'un premier mémoire déposé à l'Académie pour le Grand prix mais retiré ultérieurement. L'approche est basée sur la thèse de Poincaré sur les équations aux dérivées partielles, ainsi que sur une publication de 1878 sur l'étude des propriétés des fonctions définies par des équations différentielles à coefficients méromorphes (dont les fonctions elliptiques sont un cas particulier) dans le sillage des travaux de Charles Briot et Claude Bouquet.
 
 Un autre mémoire a été déposé le 28 mai pour le Grand prix. Dès le 29 mai, Poincaré a débuté une correspondance avec Fuchs à propos de travaux récemment publiés par ce dernier\cite{Fuchs:1880} sur les fonctions réciproques des quotients $H$ de deux intégrales indépendantes.\cite[p.174-177]{Gray:2000}. Trois suppléments ont par la suite été adressés à l'Académie les 28 juin, 6 septembre et 20 décembre.\cite{GrayWalter:1997}
 
Entre temps, et peu après le dépôt de son mémoire, Poincaré a adressé à l'Académie le 7 juin une note sur les formes cubiques ternaires.\cite{Poincar1880c}. Plus tard, le 22 novembre, il a abordé le cas particulier des cubiques décomposables en une forme quadratique et une forme linéaire.\cite{Poincar1880d} (voir aussi \cite{Poincar1886b}). La première partie d'un mémoire de synthèse sur la "théorie algébrique des formes" a été publiée dans le \textit{Journal de l'École polytechnique} en février 1881, soit à l'époque où est parue la première note sur la théorie des fonctions fuchsiennes.\cite{Poincar1881b}. L'auteur a aussi précisé la place de sa nouvelle théorie par rapport aux travaux de Jordan sur l'intégration algébrique des équations différentielles.\cite{Poincar1881c}
 
Les notes consacrées aux fonctions fuchsiennes se sont alors égrenées dans les \textit{Comptes rendus}. De premiers mémoires sont parus dans \textit{Acta mathematica} en 1882 en même temps que la seconde partie du mémoire sur la "théorie arithmétique des formes" a été publiée dans le \textit{Journal de l'École polytechnique}, à la suite d'un mémoire de Jordan sur le même sujet.\cite{Jordan1882a}. 

\subsection{Formes algébriques, réseaux et fonctions fuchsiennes}
 Envisagé selon des découpages disciplinaires, cet ensemble de préoccupations peut sembler hétéroclite. Comme nous allons le voir, le tout s'avère pourtant très cohérent sous l'angle du problème de l'extension des fonctions réciproques considérées par Fuchs. 
 
Il est bien connu que l'introduction des fonctions fuchsiennes a été présentée par Poincaré comme une généralisation du pavage du plan complexe par le réseau des fonctions elliptiques. Pour de telles fonctions, "la connaissance de la fonction à l'intérieur de l'un des parallélogrammes entraine sa connaissance dans tout le plan".\cite[p.43]{Poincar1921}. Afin de prolonger les fonctions considérées par Fuchs à partir d'une certaine région, ou "polygone curviligne", Poincaré a construit "les polygones voisins, puis les polygones voisins de ceux-ci, et ainsi de suite". De cette construction résulte un pavage du disque unité - envisagé comme un modèle du plan hyperbolique - par l'action d'un sous-groupe de $PSl_2(\mathbb{R}$) sur des polygones non euclidiens.\footnote{La généralisation de ces résultats aux fonctions invariantes par des "groupes kleinéens", i.e. des sous-groupes discrets de $PSl_2(\mathbb{C})$, allait nécessiter l'étude des polyèdres de l'espace hyperbolique et conduire à des problèmes de classification des variétés fermées à trois dimensions.}  Ainsi, dans son premier mémoire de synthèse sur les "fonctions uniformes qui se reproduisent par des substitutions linéaires", le géomètre insiste sur les 
 \begin{quote}
 [...] plus grandes analogies des fonctions fuchsiennes avec les fonctions elliptiques et modulaires qui n'en sont que des cas particuliers. [...] Nous allons rechercher s'il existe une fonction uniforme $F(\zeta)$ qui ne change pas quand on applique à $\zeta$ l'une des substitutions linéaires en nombre infini
 \[
S_i=\frac{\alpha_i \zeta + \beta_i}{\gamma_i  \zeta+ \delta_i }
 \] 
[...] Il est clair que les substitutions $S_i$ devront former un groupe et un groupe discontinu, c'est à dire que la portion du plan où la fonction $F$ existe peut être divisée en une infinité de régions $R_0$, $R_1$, ..., $R_i$, ... telles que quand $\zeta$ parcourt $R_0$, $\frac{\alpha_i \zeta + \beta_i}{\gamma_i  \zeta+ \delta_i }$ parcourt $R_i$. Ces diverses régions formeront, comme dans le cas des fonctions elliptiques, une sorte de damier dont l'ensemble ne varie pas, mais dont les cases se permuteront quand on appliquera à $\zeta$ l'une des substitutions $S_i$.\cite[p.553]{Poincar1882d}
 \end{quote}
 Dès sa note du 14 février 1881, Poincaré a insisté sur la proximité entre les groupes de substitutions unimodulaires à coefficients réels laissant invariantes les fonctions fuchsiennes et les "groupes de substitutions à coefficients entiers reproduisant une forme quadratique ternaire indéfinie à coefficients entiers". Selon lui, cette proximité fait "ressortir les liens intimes qui unissent la théorie des nombres à la question analytique qui nous occupe".\cite[p.335]{Poincar1881a} 
 
Les "dimensions respectivement algébriques et arithmétiques"\cite[p.341]{Poincar1886b} des liens entre théorie des formes et fonctions fuchsiennes ont été explicitées par Poincaré quelques années après les travaux menés pour le Grand prix.\cite{Poincar1884f} \cite{Poincar1886a} \cite{Poincar1886b} \cite{Poincar1886c} Au groupe de substitutions semblables (à coefficients entiers ou rationnels) d'une forme quadratiques ternaire est associé un "groupe fuchsien arithmétique" définissant des fonctions fuchsiennes vérifiant une propriété analogue au théorème d'addition des fonctions elliptiques.\cite{Poincar1886d} \cite{Poincar1887}
 
 Cette importance donnée aux substitutions laissant invariantes des formes algébriques résulte de ce que, jusqu'à l'été 1880, Poincaré a mené ses travaux hors du cadre de la théorie des groupes. De fait, les substitutions employées dans le premier mémoire déposé pour le Grand prix ne désignent encore de manière traditionnelle que des procédés de changements de variables et ne sont pas encore envisagées comme formant des groupes à la manière de Jordan.
 
 Comme en témoignent les trois suppléments à ce mémoire, le statut des substitutions dans les travaux de Poincaré a basculé durant l'été 1880 en même temps que ce dernier a commencé à s'appuyer sur la géométrie non euclidienne. Dans les deux cas, cette évolution est liée à des travaux sur la théorie des formes. Le deuxième supplément du 6 septembre 1880 envisage ainsi un pavage du disque par des polygones non euclidiens transformés les uns en les autres par des classes d'équivalence de substitutions laissant invariante une forme quadratique donnée. 
 
Plus tard, lorsque Poincaré a donné un compte rendu rétrospectif de son appropriation de la géométrie non euclidienne,\cite{Poincar1908a}, il a insisté sur le rôle qu'a joué le fait qu'une même représentation analytique des transformations homographiques est impliquée dans la géométrie de Lobatchevski et les transformations des formes quadratiques ternaires indéfinies (et donc des fonctions modulaires de réseaux associés à de telles formes). Plus encore, Poincaré a mis en avant le rôle joué par les interprétations arithmétiques pour se "libérer" du cas particulier de l'équation hypergéométrique sur lequel il était resté bloqué jusqu'à la fin du mois de juillet 1880.\cite[p.52]{Poincar1908a}
 
 Poincaré a en effet d'abord suivi une approche arithmétique en envisageant le plan complexe sous l'angle des réseaux représentant les classes d'équivalences de formes quadratiques. Ainsi, les deux notes successives publiées en août et novembre 1879 considèrent des systèmes de nombres complexes "existants" associés à des "nombres idéaux" par des classes d'équivalence de formes quadratiques binaires, le tout représenté par des réseaux du plan.\footnote{L'utilisation de réseaux pour représenter des formes quadratiques a été développée par Gauss (1840) et Dirichlet (1850). Indépendamment de la théorie des formes, les réseaux du plan ont aussi été étudiés par Jordan en relation avec la cristallographie et les symétries des polygones et polyèdres.\cite{Brechenmacher:2011b}. Sur les liens entre fonctions elliptiques, réseaux, formes quadratiques et nombres idéaux, voir \cite[p.897]{Poincar1879b} et \cite[p.29]{Poincar1908b}}  "Cette théorie", insistait Poincaré, "se rattache directement à celle des fonctions elliptiques, et la même méthode qui a permis de calculer les nombres corrélatifs par des intégrales définies permet d'exprimer, à l'aide d'une intégrale définie, les fonctions doublement périodiques".\cite[p.190]{Poincar1879b} 
 
Rappelons qu'il existe deux représentations géométriques des classes d'équivalence de formes quadratiques binaires. La première associe au domaine fondamental la forme réduite représentant une classe de formes : soit $\omega$ une racine (de partie imaginaire positive) de l'équation algébrique $a\omega^2+b\omega+c=0$ associée à la forme binaire (définie) $ax^2+2bxy+cy^2$, une telle forme est réduite si le point d'affixe $\omega$ est située dans le domaine fondamental.\footnote{Pour une présentation de la représentation géométrique des formes quadratiques, voir \cite[p.116]{Cahen:1908}. Sur certains épisodes de l'histoire de la réduction des formes quadratiques, voir \cite[p.488]{Schwermer2007}}  La seconde donne une représentation géométrique des classes de formes par des réseaux de parallélogrammes du plan.\cite[p.490]{Schwermer2007}. Chaque système de parallélogramme d'un réseau peut ainsi s'interpréter comme un représentant d'une classe de formes. La théorie de la réduction des formes  peut être présentée dans ce cadre géométrique ainsi que l'a notamment montré Hermite.
 
 C'est dans ce cadre hermitien de la réduction des formes que Poincaré a d'abord envisagé le problème de l'extension de la fonction considérée par Fuchs. Soit un réseau de base $(K, K')$, les deux périodes d'une fonction elliptique. Le rapport des périodes ne dépend pas du choix de la base et la fonction de Fuchs est donc définie modulo le réseau. Elle est en particulier invariante par les changements de base du réseau par l'action de substitutions unimodulaires.

Loin de toute considération en terme de groupes de substitutions, Poincaré a ainsi initialement interprété la similitude des mailles d'un réseau en termes de multiplication de nombres complexes ou de composition de formes quadratiques. Les classes d'équivalence des formes étant déterminées par des invariants, la généralisation des fonctions elliptiques doublement périodiques a été envisagée comme revenant à obtenir des fonctions admettant des "invariants arithmétiques" définis comme des séries sur un réseau.\cite[p.194]{Poincar1879b} \footnote{Ces séries et invariants jouent un rôle analogue à ceux employés pour la représentation des fonctions $P$ de Weierstrass.}  A la suite d'Hermite, Poincaré a notamment envisagé la réduction d'une forme quadratique comme généralisant l'algorithme des fractions continues\footnote{Voir \cite[p.393]{Goldstein:2007}} qu'il interprète comme un déplacement dans les mailles du réseau donné par les points d'intersections de ce dernier avec une suite de "triangles ambigus".\cite{Poincar1880c} 
 
 \subsection{L'appropriation des groupes de substitutions dans le cadre de la théorie des formes}
 C'est seulement dans un second temps, et dans le cadre d'une théorie des formes ancrée dans l'héritage hermitien, que Poincaré s'est approprié les travaux de Jordan sur les groupes de substitutions. 
 
 Comme nous l'avons déjà évoqué, Jordan a publié le 5 mai 1879 une note sur les classes d'équivalences de formes algébriques. Ce travail concerne donc directement l'objectif que s'est fixé Poincaré dans sa note du 11 août 1879 de "reconnaître si deux formes données sont équivalentes, et par quels moyens on peut passer de l'une à l'autre".\cite[p.345]{Poincar1879a} Après que Poincaré ait adressé en janvier 1880 à Jordan un courrier à ce sujet, ce dernier a conseillé à son jeune "camarade" de concentrer son attention sur les groupes de substitutions laissant une forme donnée invariante : 
 \begin{quote}
 Cette question est assez à l'ordre du jour en ce moment. M. Klein y a consacré de nombreux mémoires dans les \textit{Mathematische Annalen}. Mais il se borne aux groupes d'un nombre fini de substitutions entre deux variables. Il trouve que ces groupes se réduisent à ceux qui superposent à lui-même une pyramide régulière, une double pyramide régulière, ou l'un des polyèdres réguliers. Quant aux formes binaires correspondantes, elles se réduisent (au moins pour les groupes des polyèdres réguliers) aux fonctions entières d'un petit nombre de variables indépendantes. De mon côté, j'ai réussi à former les groupes finis pour trois variables, dans le journal de Borchardt de l'année dernière [i.e. \cite{Jordan1878} sur l'intégration algébrique des équations différentielles] ; mais je ne sais rien sur les formes correspondantes, non plus que sur les groupes où le nombre de substitutions serait infinie.\cite{Jordan-Poincar}
 \end{quote}
 
 Comme nous l'avons vu, Poincaré n'a cependant pas adopté avant l'été 1880 une approche des équations différentielles en terme de groupes de substitutions. Dans l'intervalle, Jordan a consacré le 15 mars une note sur la réduction des substitutions linéaires associées aux classes d'équivalence des formes algébriques qu'il avait étudiées l'année précédente. Poincaré s'est approprié cette problématique dans sa note sur les formes cubiques du 7 juin tandis que Jordan a publié une nouvelle note à ce sujet le 14 juin. Jusqu'à la parution de leurs deux mémoires de synthèse en 1882, les deux géomètres se sont alors focalisés sur les cas des formes quadratiques et de la représentation des nombres par des formes (\cite{Poincar1880d}, \cite{Poincar1881d}, \cite{Jordan1881a}, \cite{Jordan1881b}, \cite{Jordan1881c}).\footnote{Au sujet du problème de la représentation des nombres par les formes, envisagé dès 1879 en terme de recherche d'idéaux de norme donnée et de représentation des formes par des réseaux, voir les commentaires de Châtelet à \cite{Poincar1886a}} 
 
 \subsection{Les travaux de Jordan et Poincaré sur les formes algébriques}
 
 Considérons à présent les problèmes abordés par les travaux de Jordan et Poincaré sur les formes algébriques. En mai 1879, la première note de Jordan a généralisé des travaux d'Hermite\cite{Hermite1854} aux formes à $n$ variables et de degré $m$. Elle est basée sur la théorie de la réduction des formes de ce dernier et notamment sur la distinction entre classes d'équivalence algébrique (i.e. modulo $Sl_n(\mathbb{R}$)) et arithmétique (modulo $Sl_n(\mathbb{Z}$) ou $Sl_n(\mathbb{Z}(i))$) des formes à coefficients entiers. 
 
 A une classe de formes (quadratiques, binaires ou décomposables en facteurs linéaires) est associée une forme réduite dont les coefficients sont définis par récurrence comme des minimums successifs de valeurs prises par les formes d'une même classe d'équivalence.\footnote{Sur la théorie de la réduction chez Hermite, voir \cite[p.391]{Goldstein:2007}.}  La borne donnée par Hermite pour de tels minimums a été affinée en 1873 par Aleksandr Korkin et Egor Zolotarev pour les cas des formes ternaires et quaternaires. Elle a permis à Jordan d'énoncer un résultat de finitude du nombre de classes d'équivalences de formes $n$-aires de discriminant non nul \cite[p.1422]{Jordan1880a} sur le modèle de son théorème de finitude des groupes linéaires finis.\footnote{Comme Hermite, Korkin et Zolotarev introduisent des réduites en maximisant ou minimisant les coefficients de formes sur des classes d'équivalence par des raisonnements analytiques. A la classe des matrices $U^*FU$ Hermite a associé deux invariants à partir desquels il mesure la grandeur de la classe : le déterminant commun $D$ et la plus petite valeur $\mu$ prise par la forme hermitienne correspondante aux points distinct de l'origine ayant des coordonnées entières. Hermite montre alors qu'il existe une matrice de la classe, la réduite, dont les éléments sont bornés en valeur absolue par $D$, $\mu$ et $n$.
 
 Plus tard, Hurwitz (1884) a inscrit les travaux de Jordan et Poincaré dans le cadre de l'approche de Kronecker sur la théorie des nombres de classes tandis que les travaux de Minkowski (1885) poursuivent les questions de réduction des formes à la suite d'Hermite, Korkin, Zolotarev, Jordan et Poincaré. Au sujet du rôle joué par les travaux de Minkowski sur la représentation des formes par des réseaux dans la genèse de la géométrie des nombres, voir \cite[p.490]{Schwermer2007}. Sur le développement de la géométrie des nombres à partir du problème de l'estimation du minimum des formes quadratiques, voir plus généralement \cite{Gauthier:2009}.}  
 
 Plus précisément, Jordan a associé à chaque classe d'équivalence algébrique de formes $\Phi$ à coefficients entiers (réels ou complexes) de degré $m$ et de $n$ variables (correspondant aux transformations $\Phi \to \Phi U$, avec $U$ une substitution unimodulaire) une classe d'équivalence arithmétique de formes quadratiques (ou hermitiennes) $F= \Phi^* \Phi$ (modulo les transformations $F \to U^*FU$).\footnote{En termes actuels, Jordan étudiait les classes de matrices inversibles complexes qui se déduisent les unes des autres par multiplication à droite par une matrice unimodulaire, c'est à dire $Gl_n(\mathbb{C})/Sl_n(\mathbb{Z}(i)$).\cite[p.XV]{Dieudonne:1962}}  Il a dénommé "substitution réduite" la substitution $U$ transformant $F$ en sa réduite $G$ au sens de Korkin et Zolotarev. La notion de réduite a ainsi été transférée des classes d'équivalences arithmétiques des formes quadratiques aux classes d'équivalences algébriques des formes de degré $n$.\cite[p.135]{Jordan1880b}. Ce transfert a permis à Jordan d'employer pour la classification des formes ses travaux sur les groupes linéaires, en particulier sa pratique de réduction canonique.
 
 Jordan s'est aussi appuyé sur les travaux d'Hermite sur les substitutions dites "semblables" qui laissent une forme quadratique invariante. Il a notamment mené une étude des substitutions transformant une réduite de degré $n$ en elle-même ou en une autre réduite. Dans le cas des classes d'équivalences algébriques, il a eu recours à sa réduction canonique afin d'obtenir des décompositions de substitutions en produits de substitutions infinitésimales. L'objectif affiché était en effet d'assurer la possibilité d'une réduction effective des formes : "on voit par là qu'on n'aura qu'un nombre limité d'essais à faire pour constater l'équivalence de deux formes réduites et trouver les transformations de l'une dans l'autre".\cite[p.133]{Jordan1880b}
 
 Si cette problématique manifeste un héritage de principes d'effectivité d'Hermite, Jordan s'est néanmoins davantage intéressé aux propriétés générales des groupes de substitutions de $n$ variables qu'aux objets hermitiens eux-mêmes. Les travaux menés en parallèle par Poincaré manifestent une fidélité bien plus forte aux principes hermitiens. Ce dernier s'attache notamment à mener des calculs permettant des résultats explicites sur des objets spécifiques. Ainsi, tandis que Jordan a énoncé un théorème pour des formes de degré $m$, Poincaré s'est attelé au cas des formes quadratiques et cubiques de discriminant nul. Il a aussi cherché à exprimer explicitement les équations différentielles linéaires à intégrales algébriques correspondant à la classification donnée par Jordan aux sous-groupes finis de $Gl_n(\mathbb{C})$ ($n \leqslant 4$).
 
 \subsection{L'appropriation par Poincaré de procédés algébriques de Jordan}
 
 Dans ses premières interventions, Poincaré a généralisé le théorème de finitude de Jordan aux formes quadratiques de discriminant nul. Il a aussi démontré que, au contraire du cas générique considéré par Jordan, le cas particulier des cubiques de discriminant nul présente une infinité de classes d'équivalences. 
 C'est à cette occasion que Poincaré s'est approprié l'approche de Jordan. Il a en effet associé à chaque cubique son groupe des substitutions semblables afin de ramener le problème de l'équivalence algébrique à celui de la classification des substitutions linéaires. Sur le modèle des travaux de Jordan, cette classification est basée sur la réduction canonique de substitutions. En termes actuels, il s'agit de déterminer les sous-groupes de $Gl_n(\mathbb{Z}$) laissant invariante à gauche une forme homogène de degré 3 à 3 ou 4 variables et, réciproquement, de déterminer toutes les formes laissées invariantes par un sous-groupe $G$ de $Gl_n(\mathbb{Z})$.
 
 En même temps qu'il s'est approprié la pratique de réduction de Jordan, Poincaré a adopté le point de vue selon lequel "le fondement de [ses] recherches sur les fonctions fuchsiennes est l'étude des groupes discontinus contenus dans le groupe linéaire à une variable".\cite[p.840]{Poincar1882c}
 
 Le vocabulaire employé témoigne du rôle joué par les travaux sur l'intégration algébrique que Jordan avait signalé à Poincaré dans sa lettre de janvier 1880. En effet, Poincaré désigne comme forme canonique la seule forme diagonale comme Jordan l'a fait en 1878 pour son étude des sous-groupes finis de $Gl_3(\mathbb{C})$. Plus encore, le premier s'appuie à plusieurs reprises, en 1883 et en 1884, sur l'argument développé par le second en 1878 pour montrer que les substitutions d'ordre fini sont diagonalisables.\footnote{En 1883, par exemple, Poincaré a prouvé qu'un groupe continu $G$ contient toujours une infinité de sous-groupes par le fait que si $G$ contient une substitution admettant une certaine équation caractéristique, alors une telle substitution se trouve dans tous les sous-groupes de $G$. Comme dans la note "Sur les nombres complexes" de 1884, le cas non commutatif est associé à l'occurrence de racines caractéristiques multiples et Poincaré s'est appuyé sur les réductions canoniques données par Jordan aux substitutions infinitésimales laissant invariante une forme algébrique.} 
 
 \section{Les pratiques algébriques de Poincaré}
 Tout en s'appropriant certains travaux de Jordan, Poincaré a élaboré des pratiques algébriques qui lui sont propres et dont nous allons à présent donner les principales caractéristiques.
 
 \subsection{Algèbre, arithmétique et géométrie}
 
À la différence de Jordan, Poincaré articule systématiquement ses considérations algébriques sur les groupes de substitutions à des problématiques arithmétiques sur les classes d'équivalences de formes ainsi qu'à des interprétations géométriques. 
 
Si $x_1$, $x_2$, $x_3$ sont les coordonnées d'un point du plan, la cubique $F=0$ représente une courbe $C$. Une substitution $T$ pour laquelle $F$ est reproductible peut s'interpréter comme une homographie entre deux points de $C$ :\footnote{Voir en particulier l'importance de la représentation géométrique dans la généralisation par Poincaré du théorème de Jordan de finitude des classes.\cite{Poincar1882a}}
 \begin{quote}
 
 Supposons que $T$ soit de la première ou de la deuxième catégorie, c'est-à-dire que l'équation [caractéristique] ait quatre racines distinctes, $\lambda_1, \lambda_2, \lambda_3, \lambda_4$. Il y a alors quatre plans reproductibles par $T$. Imaginons que l'on fasse un changement de coordonnées $\Sigma$ en prenant pour nouveau tétraèdre de référence le tétraèdre formé par ces quatre plans. Il est clair que la transformée de $T$ par $\Sigma$ est
 \[
\Sigma^{-1}T\Sigma=
 \begin{vmatrix}
\lambda_1 & 0 & 0 & 0 \\
0 & \lambda_2 & 0 & 0   \\
0 & 0 & \lambda_3 & 0 \\
0 & 0 & 0 & \lambda_4
\end{vmatrix}
\]
[...] Il suit de là que si $T$ est de première ou deuxième catégorie, on peut choisir $\Sigma$ de telle sorte que $\Sigma^{-1}T\Sigma$ soit par conséquent canonique  [...]. Les transformations ternaires (canoniques) de la troisième catégorie se classent en deux types :
 \[
 Type \ A...
  \begin{vmatrix}
\alpha & 0 & 0 \\
0 & \alpha & 0   \\
0 & 0 & \beta 
\end{vmatrix}
\ , \ Type B ...
  \begin{vmatrix}
\alpha & 0 & 0 \\
0 & \alpha & 0   \\
0 & 0 & \alpha 
\end{vmatrix}
\]
 [...] La question qui se pose est de trouver les points, les droites et les plans reproductibles par la transformation $T$ et de discuter complètement le problème [...]. Passons maintenant aux transformations de la quatrième catégorie ; on ne peut pas les réduire à la forme canonique, mais on peut choisir $\Sigma$ de façon à ramener $\Sigma^{-1}T\Sigma$ à sa forme la plus simple. Ainsi, les transformations ternaires se divisent en deux types dont je donne ici les formes les plus simples : 
 \[
 Type \ A_1...
  \begin{vmatrix}
\beta& 0 & 0 \\
0 & \alpha & 0   \\
0 & \gamma & \alpha 
\end{vmatrix}
\ , \ Type B_1 ...
  \begin{vmatrix}
\alpha & 0 & 0 \\
\beta & \alpha & 0   \\
\gamma & \delta & \alpha 
\end{vmatrix}  \cite[p.34]{Poincar1882e}
\]
\end{quote}
Les quatre catégories de formes canoniques de substitutions quaternaires sont mises en correspondance avec sept familles de cubiques, chacune représentée par une forme réduite ou "forme la plus simple du type considéré".\cite[p.31]{Poincar1881d} A chacune de ces formes canoniques sont aussi associés des invariants et covariants au sens des travaux de Sigfried Aronhold et d'Alfred Clebsch. 

Par exemple, les quatre premières familles correspondent aux formes non décomposables. Poincaré y distingue deux couples de deux familles selon que le discriminant de la forme soit nul ou non, ce dernier cas correspondant à l'absence de points doubles. Chaque famille est caractérisée par une forme canonique comme la forme $F$ ci-dessous : 
\begin{quote}
La forme $F=6xyz+\beta(x^3+y^3+z^3)$ est reproductible par les transformations suivantes, qui appartiennent à la deuxième catégorie (nous n'écrivons que les transformations réelles) :
\[
\begin{vmatrix}
0 & 1 & 0 \\
0 & 0 & 1 \\
1 & 0 & 0
\end{vmatrix}
,
\begin{vmatrix}
0 & 0 & 1 \\
1 & 0 & 0 \\
0 & 1 & 0
\end{vmatrix}
 ,
 \begin{vmatrix}
0 & 1 & 0 \\
1 & 0 & 0 \\
0 & 0 & 0
\end{vmatrix}
,
\begin{vmatrix}
0 & 0 & 1 \\
0 & 1 & 0 \\
1 & 0 & 0
\end{vmatrix}
 ,
 \begin{vmatrix}
1 & 0 & 0 \\
0 & 0 & 1 \\
0 & 1 & 0
\end{vmatrix}
 \]
 toutes ces transformations se réduisent à des permutations entre les lettres $x, y, z$. C'est là un résultat qu'il était aisé de prévoir en effet, le système des trois droites $x=0$, $y=0$, $z=0$, est le seul système de trois droites réelles sur lesquelles se distribuent les neuf points d'inflexion. Toute transformation réelle qui reproduit la forme proposée doit donc reproduire le système des trois droites : elle doit donc se ramener à une permutation entre ces trois droites.\cite[p.57]{Poincar1882a}
 \end{quote}
 
 \subsection{Entre Jordan et Hermite : formes reproductibles et groupes de substitutions semblables}

 Une autre caractéristique de la réduction canonique de Poincaré est que celle-ci articule constamment formes reproductibles et substitutions linéaires. Ceci témoigne du rôle joué par l'héritage hermitien dans l'utilisation par Poincaré des travaux de Jordan sur les groupes de substitutions linéaires. 
 
Chez Poincaré, une notion de réduite générale englobe en réalité la réduction canonique de Jordan. Cette notion est avant tout envisagée comme un leg des travaux d'Hermite sur les formes quadratiques. Elle joue ainsi un rôle clé dans le texte ci-dessous dans lequel Poincaré présente son approche sur les formes cubiques :
 \begin{quote}
 Les divers problèmes qui se rattachent à la théorie des formes quadratiques binaires ont été résolus depuis longtemps, grâce à la notion de réduite [...]. La notion de réduite s'étend sans peine aux formes quadratiques définies par un nombre quelconque de variables. [...] Généraliser une idée aussi utile, trouver des formes jouant, dans le cas général, le même rôle que les réduites remplissent dans le cas des formes quadratiques définies, tel est le problème qui se pose naturellement et que M. Hermite a résolu de la façon la plus élégante dans divers Mémoires insérés dans les Tomes 44 et 47 du Journal de Crelle.
 
M. Hermite s'est bornée à l'étude des formes quadratiques définies ou indéfinies et des formes décomposables en facteurs linéaires ; mais sa méthode peut s'étendre sans difficulté au cas le plus général. [...]. Les résultats auxquels je suis arrivé s'appliquent à une forme quelconque; mais, ne voulant pas sacrifier la clarté à la généralité, je me suis restreint aux formes qui sont les plus simples parmi celles que M. Hermite avait laissé de côté. [...] Les plus simples de toutes les formes, après les formes quadratiques et les formes décomposables en facteurs linéaires, sont les formes cubiques ternaires.\cite[p.28]{Poincar1881d}
 \end{quote}
 Mais que recouvre au juste cette notion de réduite ? Les travaux de Poincaré n'en donnent aucune définition précise. Les réduites prennent au contraire des formes et des significations variables selon les problèmes particuliers considérés. L'idée générale de réduite n'est quant à elle normée que par un idéal de simplicité maximale dont la signification est relative à chaque problème considéré: 
 \begin{quote}
 On choisira dans chaque type ou dans chaque sous-type, pour le représenter, une des formes de ce type ou de ce sous-type que l'on appellera la forme canonique $H$. Le choix de la forme $H$ est peu près arbitraire ; toutefois on sera conduit, dans la plupart des cas, à choisir de préférence la forme la plus simple du type considéré.\cite[p.203]{Poincar1881a}
 \end{quote}
 Cette norme mise en avant par Poincaré peut sembler proche de l'idéal de simplicité attaché par Jordan à sa réduction canonique. Elle en diffère cependant. Chez Jordan, la simplicité est effet intrinsèquement associée à la décomposition maximale des substitutions ou des groupes de substitutions. Elle présente ainsi un caractère algébrique et a d'ailleurs été critiqué comme telle par Kronecker qui a reproché le caractère formel de l'expression forme canonique car celle-ci ne revêt aucun sens précis mais prend au contraire des significations variées selon le problème considéré.
 
Contrairement à Jordan, Poincaré ne s'appuie pas en 1881 sur un idéal attribuant un caractère essentiel à la notion de groupe. Il associe au contraire la notion de réduite à la problématique des substitutions semblables des formes reproductibles. Il articule ainsi constamment arithmétique des classes d'équivalences des formes, discussions algébriques sur la multiplicité des racines caractéristiques, interprétations géométriques et raisonnements analytiques de majorations ou minorations. 
 
Ces articulations sont explicitement valorisées par Poincaré. Elles manifestent un héritage des idéaux d'Hermite sur l'unité des mathématiques.\cite[p.399]{Goldstein:2007} 

Considérons l'exemple des substitutions dites hyperfuchsiennes. Ces substitutions de trois variables ont été introduites par Picard dans des travaux sur des formes hermitiennes ternaires. Tout comme ceux de Poincaré, ces travaux sont marqués par l'héritage d'Hermite. La réponse de Poincaré à Picard est basée sur l'étude de formes semblables et vise à montrer que des "considérations, arithmétiques, algébriques ou géométriques, permettent d'obtenir une infinité de groupes discontinus".\cite[p.840]{Poincar1882c} 
Nous avons vu que la note de 1884 "Sur les nombres complexes" s'inscrit elle aussi dans les problématiques de généralisation qui ont accompagné les travaux  sur les groupes hyperfuchsiens. Tout comme de nombreux autres travaux contemporains de Poincaré, cette note inscrit les procédés algébriques de Jordan au c\oe{}ur d'articulations entre différents domaines.

Ces articulations s'appuient néanmoins sur une distinction explicite entre les "questions purement algébriques" pour lesquelles Poincaré reprend les travaux de Jordan et l'"étude arithmétique" basée sur des méthodes d'Hermite comme la réduction continuelle. Plus précisément, l'approche algébrique - et en particulier la réduction canonique des substitutions - est fréquemment présentée comme relevant de "systèmes de définitions qui seront nécessaires par la suite".\cite[p.203]{Poincar1881a} Elle regroupe notamment un ensemble de préliminaires communs à de nombreux travaux de Poincaré : notation des substitutions linéaires par des tableaux, équations aux multiplicateurs, réduction canonique des substitutions et des formes, points doubles, formes reproductibles par des substitutions semblables etc.\footnote{ On trouve ainsi de tels préliminaires dans des études des formes,\cite[p.34]{Poincar1882a} des fonctions fuchsiennes,\cite[p.108]{Poincar1882e} des groupes continus et équations aux dérivées partielles,\cite{Poincar1883a} des nombres complexes,\cite{Poincar1884d} de l'intégration algébrique,\cite[p.300]{Poincar1884e},  de l'intégration par des séries,\cite[p.316]{Poincar1886g} de l'arithmétique des fonctions fuchsiennes,\cite[p.463-505]{Poincar1887}, des "exposants caractéristiques" du problème des trois corps,\cite[pp.293;337-343]{Poincar1890} des homologies de l'Analysis Situs,\cite[p.342-345]{Poincar1900} ou encore les groupes continus.\cite[p.216-252]{Poincar1901}}

\subsection{Le calcul des Tableaux}
Les paragraphes précédents se sont attachés à décrire la spécificité des pratiques algébriques élaborées par Poincaré au début des années 1870-1880 à partir d'un double héritage des travaux de Jordan et d'Hermite. Il nous reste encore à aborder une caractéristique importante de ces pratiques que nous avons notamment vu se manifester dans la note de 1884 : le  calcul des Tableaux.

Le terme de Tableau a été utilisé en France tout au long du XIX\up{e} siècle. Il ne s'agit cependant pas d'une simple notation. Ce terme a en effet généralement été utilisé pour désigner la "forme" obtenue par le regroupement d'un ensemble d'objets du même type. Il s'accompagne ainsi le plus souvent de l'expression "formons à présent le Tableau de tels objets", les objets en question pouvant être des nombres, des substitutions, etc. Dès le début du XIX\up{e} siècle, cette forme a présenté un caractère opératoire. Elle permet, d'une part, de manipuler des collections d'objets comme un tout et, d'autre part, d'en extraire des sous-collections. Les Tableaux ont notamment été fortement  mobilisés dans le cadre du développement de la théorie des formes quadratiques. En tant que tout, ils permettent de représenter une classe d'équivalence par un représentant canonique. En tant que collection, un tel représentant fournit en son sein une liste d'éléments invariants par équivalence. 

Dans les années 1850, des auteurs comme Hermite, Sylvester et Cayley ont développé l'usage qu'avait fait Cauchy des Tableaux lors de ses travaux sur l'équation séculaire. Plus tard, les Tableaux ont notamment circulé chez des épigones d'Hermite. Chez de nombreux auteurs, l'expression "calcul des Tableaux" identifie alors un ensemble de procédés de décompositions, comme dans le cas des matrices et des mineurs. Avant le début du XX\up{e} siècle, ces procédés n'étaient cependant pas envisagés comme des méthodes pouvant s'intégrer à un cadre théorique. 

Par ailleurs, au cours du XIX\up{e} siècle, de nombreuses situations qui avaient d'abord été représentées par des Tableaux ont donné lieu à des notions spécifiques, comme dans le cas des formes quadratiques, des groupes de substitutions, des algèbres associatives etc. Le terme de Tableau a cependant subsisté en tant que tel pour désigner la forme algébrique prise en général par un regroupement d'objets. Cet usage est proche de la définition qu'avait donné Sylvester en 1851 à la notion de matrice. Rappelons que ce dernier avait distingué de la notion de déterminant la "forme carrée" dont peuvent être extraits les mineurs. En France, de nombreux auteurs utilisaient à la fois les termes Tableau et matrice. Ils réservaient cependant ce dernier aux problèmes dans lesquels intervenaient des calculs de déterminants et de mineurs. Le calcul des Tableaux présentait quant à lui un caractère plus général et plus transversal.

Compte tenu de cette longue tradition, la présence des Tableaux dans les travaux de Poincaré n'est pas surprenante. Nous allons voir cependant que l'usage qu'en fait ce dernier présente un caractère spécifique. En effet, lors de leurs interactions sur la théorie des formes au début des années 1880, Jordan et Poincaré ont doté les Tableaux d'un nouveau caractère opératoire.

Les Tableaux ont fait leur apparition au sein de l'\oe{}uvre de Jordan à la fin des années 1870, lors des premiers travaux publiés par ce dernier dans le cadre hermitien de la réduction des formes. A partir de cette époque, ils ont cohabité avec les deux autres notations que Jordan a employées depuis sa thèse pour représenter les substitutions :
\begin{center}
$S=
\begin{vmatrix}
x_1 & a_{11}x_1+a_{12}x_2+... \\
x_2 & a_{21}x_1+a_{22}x_2+... \\
x_3 & a_{31}x_1+a_{32}x_2+... 
\end{vmatrix}
=
\begin{vmatrix}
a_{11} & a_{12} & ... \\
a_{21} & a_{22} & ... \\
... & ... & ... 
 \end{vmatrix} $
 \cite[p.117]{Jordan1880b}
 \end{center}

La notation par une lettre, $S$, est dotée d'une opérativité symbolique et permet notamment de représenter la relation d'équivalence $S'=USV$. La représentation analytique en lignes permet, quant à elle, de représenter simultanément des décompositions des variables $x_i$ en blocs et de la substitution en substitutions partielles laissant invariant chaque bloc. Cette seconde notation est par conséquent un élément constitutif des procédés de réduction de Jordan. L'énoncé du théorème de réduction canonique en donne un bon exemple. 

En 1880, Jordan a transféré aux Tableaux les caractères opératoires des deux autres formes de représentations. La multiplication symbolique par des substitutions unimodulaires s'articule ainsi à des opérations arithmétiques sur les lignes et les colonnes.\cite[p.152]{Jordan1880c} \footnote{L'étude des procédés tabulaires arithmétiques développés par Smith (1861) et par Kronecker (1884-1891) permet, par contraste, de mettre en valeur l'originalité des procédés de Jordan et Poincaré. Voir \cite{Brechenmacher:2006a}.} Surtout, des multiplications de Tableaux entre eux permettent de décomposer la notation tabulaire en sous-Tableaux représentant les étapes successives du procédé de réduction canonique. Les Tableaux donnent ainsi une représentation de la pratique qui consiste à réduire un problème général (représenté par le Tableau principal) en une chaîne de maillons les plus simples (les sous-Tableaux), eux même dotés d'une opérativité algébrique.\cite[p.130]{Jordan1880b}

Poincaré a développé une pratique des Tableaux semblable à l'époque de ses interactions avec les travaux de Jordan sur les formes algébriques. Contrairement à ce dernier, il ne recourt cependant pas à la représentation de substitutions de $n$ variables mais en écrit la généralité sur des paradigmes de Tableaux canoniques. Il emploie ainsi le plus petit nombre de variables permettant une présentation exhaustive des cas possibles.\footnote{Au sujet du mode d'écriture de la généralité par paradigmes chez Poincaré, voir \cite{Robadey2004}.}

Dès le début des années 1880, la réduction des Tableaux a été intégrée aux préliminaires algébriques de Poincaré. Elle a par la suite joué un rôle dans des thèmes variés comme les groupes des équations linéaires,\cite[p.301]{Poincar1884e} les formes, nombres complexes idéaux et réseaux,\cite[pp.375;389-393]{Poincar1886b} ou les systèmes de périodes des fonctions abéliennes.\cite[p.334]{Poincar1884g} \cite{Poincar1886e} \cite[pp.319-330;360-361]{Poincar1886f}

Nous avons déjà évoqué que les Tableaux présentaient traditionnellement en France un caractère transversal au sens où ils n'étaient pas identifiés à un objet ou à un cadre théorique spécifique. Dans les travaux de Poincaré, le calcul des Tableaux permet des exploitations parallèles de domaines d'interprétations et de problèmes. Il supporte ainsi à plusieurs reprises des transferts réciproques d'un champ à un autre par des analogies portées par des procédés opératoires.\footnote{Au sujet des transferts analogiques portés par des procédés opératoires, voir \cite{Durand-Richard:2008}.} Un premier exemple en est donné par les transferts entre périodes des fonctions abéliennes et idéaux, formes binaires ou réseaux, tous représentés par des formes réduites de Tableaux à coefficients entiers modulo des opérations par des Tableaux unimodulaires.\cite[p.406-431]{Poincar1886a} \footnote{Voir à ce sujet les notes de Châtelet à l'édition des \oe{}uvres complètes de \cite{Poincar1886a}. Étant donné une forme binaire $f(x,y)$ telle que $\Phi(x) = f(x,1)$ est un polynôme unitaire et soit $\alpha_1$ un zéro de ce polynôme. L'anneau des entiers algébriques $x_0+\alpha_1x_1+ a_1^2x_2+...+a_1^{m-1}x^{m-1}$ ($x_i$ entiers) est isomorphe à un anneau de matrices $A$ à coefficients entiers qui représentent les tables de multiplication des entiers par les $m$ premières puissances de $\alpha_1$, i.e. par exemple $(1, \alpha_1, ..., \alpha_1^{m-1}) \times \alpha_1 =(1, \alpha_1, ..., \alpha_1^{m-1}) \times A$. Un sous module de l'anneau peut être défini par une base de $m$ entiers, que l'on peut écrire $(1, \alpha_1, ... , \alpha_1^{m-1})\times S$ , où $S$ est une matrice à termes entiers définie modulo la multiplication à droite par une matrice unimodulaire. Pour que ce module soit un idéal, il est nécessaire et suffisant que $S^{-1} \times A \times S$ soit à termes entiers. L'étude des idéaux premiers se ramène alors à celle des formes de matrices $S$.} D'autres exemples de transferts incluent l'introduction des déterminants d'ordre infini,\cite[p.100]{Poincar1886c} les longs développements consacrés à la réduction des Tableaux dans le mémoire "Les fonctions fuchsiennes et l'arithmétique",\cite[p.485-499]{Poincar1887} et leur utilisation ultérieure dans la définition des nombres de Betti de l'Analysis situs.\cite[p.310-328]{Poincar1899} \cite[p.342-345]{Poincar1900}

La décomposition d'un Tableau donne une forme à la pratique algébrique de classification par réduction canonique qu'a élaboré Poincaré à la suite des travaux de Jordan et Hermite.

Au-delà des travaux de Poincaré, cette pratique de décomposition des Tableaux a par la suite circulé au sein d'un réseau de textes spécifique entre 1880 et 1914. L'une des caractéristiques principales de ce réseau est la référence commune aux travaux d'Hermite et Jordan. Les textes y sont en outre majoritairement publiés par des auteurs français - parmi lesquels Poincaré, Jordan, Léon Autonne, Jean-Armand de Séguier ou Albert Châtelet. Ce groupe ne se réduit cependant pas à une tradition nationale comme l'illustre la présence de textes d'auteurs cultivant l'héritage d'Hermite, comme Hermann Minkowski, ou celui de Jordan, comme Leonard Dickson. La notion de réduite joue également un rôle important dans ce réseau. Comme chez Poincaré, cette notion n'y fait cependant pas l'objet d'une définition explicite. Elle présente au contraire le caractère polysémique des "formes les plus simples" en fonction du problème considéré. 

Les travaux de Poincaré ont joué un rôle de modèle pour certains auteurs importants de ce réseau de textes. Autonne a notamment suivi ce modèle en introduisant ses travaux sur les équations différentielles par des préliminaires sur la réduction des Tableaux et ses interprétations algébriques, géométriques et arithmétiques.\cite{Autonne:1885a} \footnote{Poincaré a d'ailleurs personnellement assuré l'encadrement des travaux de thèse d'Autonne à l'École polytechnique sur les intégrales algébriques des équations différentielles linéaires entre janvier et septembre 1881. Il a notamment incité ce dernier à se plonger dans le \textit{Traité} de Jordan.\cite{Autonne-Poincar} Prolongeant les travaux de Poincaré sur l'intégration algébrique, Autonne a cherché à expliciter la forme des équations algébriques admettant pour racines les intégrales des types d'équations différentielles associées aux différents types de sous-groupes finis de $Gl_3(\mathbb{C})$ donnés par Jordan (à la même époque, Goursat a consacré ses premiers travaux à un problème semblable). L'échange épistolaire entre Autonne et Poincaré est essentiellement consacré à des questions liées aux racines caractéristiques de substitutions en relation avec les travaux de Jordan (en particulier \cite[p.200]{Jordan1878}.} Les nombreux traités sur la théorie des formes, les groupes, les équations différentielles et la géométrie algébrique que ce dernier a publié au début du XX\up{e} siècle ont institutionnalisé l'usage de tels préliminaires.\cite[p.599-628]{Brechenmacher:2006a}. Ces travaux, tout comme ceux de Séguier dans le cadre algébrique des groupes finis ou ceux de Châtelet en théorie des nombres, ont par la suite joué un rôle important dans le développement de la théorie des matrices canoniques des années 1920.

Le rôle de modèle joué par les préliminaires algébriques de Poincaré se manifeste aussi par l'adoption au sein du réseau du calcul des Tableaux de l'expression forme canonique pour désigner une forme diagonale. Au contraire, des travaux qui, comme ceux de Cartan,\cite[p.639]{Cartan:1894} ne s'appuient pas sur les Tableaux mais sur une référence directe aux premiers travaux de Jordan sur les groupes de substitutions font un usage important de la réduction canonique au sens de ce dernier.

\section*{Conclusion}

Nous avons vu que la note "Sur les nombres complexes" de 1884 est loin d'être isolée dans les travaux de Poincaré. Elle s'insère non seulement dans un ensemble de travaux autour des fonctions fuchsiennes mais témoigne aussi plus généralement de pratiques algébriques élaborées par Poincaré à partir d'héritages de travaux d'Hermite et de Jordan. Rappelons les principales caractéristiques de ces pratiques algébriques. 

Tout d'abord, nous avons vu que ces pratiques s'ancrent dans une culture algébrique commune portée par l'équation séculaire. Elles s'appuient ainsi sur des procédés traditionnels de manipulation des systèmes linéaires par la décomposition polynomiale de leur équation caractéristique. Depuis la fin du XVIII\up{e} siècle, ces procédés étaient indissociables d'un enjeu de généralité que l'on trouve encore chez Poincaré et consistant à généraliser des situations étudiées pour 2 variables à 3 ou $n$ variables. Cet enjeu de généralité participe de la nature algébrique de ces procédés. En, effet, les significations géométriques, mécaniques ou arithmétiques associés à ces derniers pour les cas de 2 variables se trouvent bien souvent perdues, ou du moins changées, dans le cas de $n$ variables. D'un point de vue algébrique, le problème de la généralisation à plus de deux variables de situations de proportionnalité s'est notamment posé en terme d'ordre de multiplicité des racines caractéristique, problème qui s'est avéré lié à la non commutativité de la multiplication matricielle.

Durant le dernier tiers du XIX\up{e} siècle, cette culture commune s'est cependant disloquée progressivement et des lignes divergentes se sont développées à des niveaux locaux. C'est cependant toujours par l'intermédiaire du problème traditionnel de l'équation séculaire que la réduction canonique de Jordan est entrée en contact avec les diviseurs élémentaires de Weierstrass. De même, des auteurs comme Scheffers ou Klein ont réagi aux travaux de Poincaré sans pour autant s'approprier la pratique de réduction canonique qui y joue pourtant un rôle essentiel ; réciproquement, ce dernier ne reprend jamais les diviseurs élémentaires auxquels recourent systématiquement Scheffers et Klein.

Ces différentes pratiques algébriques locales ne se limitent pas à des procédés techniques mais présentent au contraire des aspects culturels qui engagent le statut de l'algèbre, des philosophies internes de généralité ou des valeurs épistémiques de simplicité. Les espaces pertinents pour identifier de telles cultures algébriques ne correspondent cependant souvent pas plus à des aires géographiques ou institutionnelles qu'à des théories ou disciplines mathématiques. Nous avons vu au contraire qu'une lecture de la note de 1884 en terme de la rencontre entre une théorie des algèbres anglo-américaine et une théorie des groupes continus européennes a eu pour effet de brouiller les dimensions collectives de ce texte. L'approche que nous avons suivi dans cet article a consisté à proposer différents éclairages sur la note de 1884 en plaçant cette dernière dans plusieurs espaces intertextuels, ou réseaux de textes, à différentes échelles : temps long de la culture commune de l'équation séculaire, temps moyen de la réception des travaux de Jordan, temps court des travaux de Poincaré qui amènent à l'émergence des fonctions fuchsiennes. Si l'articulation de tels réseaux de textes avec des espaces sociaux ou institutionnels est problématique, nous avons vu que ces réseaux éclairent tout à la fois sur des normes collectives et sur la position singulière d'auteurs comme Poincaré. 

L'un des apports cruciaux de l'ensemble de travaux dans lequel s'insère la note de 1884 est l'articulation entre les héritages de Jordan et Hermite. 

D'un côté, nous avons vu le rôle important joué par l'héritage des travaux de Jordan sur les substitutions pour l'adoption progressive par Poincaré d'une approche des équations différentielles centrée sur les groupes et non plus sur le cadre hermitien des classes d'équivalences de formes. Le rôle joué par cet héritage a pourtant été ignoré de l'historiographie de la théorie des fonctions fuchsiennes. En effet, les travaux de Jordan ont souvent été rattachés au contexte de l'émergence  de la théorie des groupes comme domaine mathématique autonome. Ils ont pour cette raison été fondus dans un paysage englobant les travaux de Klein, Lie, Frobenius, etc. Or, les pratiques  des groupes sont très différentes chez ces mathématiciens. Par ailleurs, le statut de la notion de groupe a précisément évolué au  tournant des années 1870-1880 dans le cadre des travaux sur les équations différentielles. Mais contrairement au rôle que leur a souvent attribué l'historiographie, les travaux de Jordan ne représentent qu'une approche des groupes parmi d'autres. Ils n'ont en particulier qu'été très partiellement appropriés par Klein et n'ont joué qu'un rôle limité dans une évolution qui doit s'envisager comme résultant de travaux variés menés sur le temps long du XIX\up{e} siècle. Là où une lecture contemporaine reconnaîtrait une même notion de groupe, il est alors nécessaire de distinguer des pratiques différentes dont les circulations peuvent créer des dynamiques collectives. Ces différences sont importantes pour saisir la spécificité de l'approche de Poincaré sur les équations différentielles dans le double héritage des travaux d'Hermite et de Jordan.

Car, d'un autre côté, Poincaré n'a que partiellement adopté l'approche de Jordan, et par l'intermédiaire du prisme hermitien de la théorie des formes. Contrairement à ce dernier, Poincaré ne se satisfait en effet jamais de l'étude de groupes de substitutions à $n$ variables. A la suite d'Hermite, il envisage toujours des substitutions opérant sur des formes algébriques et articule ainsi non seulement arithmétique, analyse et algèbre comme ce dernier mais aussi des interprétations géométriques. Comme Hermite, et contrairement à Jordan, Poincaré cherche notamment à obtenir des solutions effectives pour les cas de deux ou trois variables où celles-ci pouvaient prendre des significations non uniquement algébriques. 

La pratique algébrique qui a résulté de cet héritage mêlé se manifeste notamment dans le rôle qu'attribuent les préliminaires algébriques de Poincaré aux groupes jordaniens de substitutions linéaires à $n$ variables tout en en écrivant la généralité sur des exemples paradigmatiques permettant les calculs effectifs dictés par les principes hermitiens. Elle a indéniablement contribué à la capacité de Poincaré d'intervenir dans un spectre large allant de l'arithmétique aux équations différentielles en passant par les nombres complexes. 

Nous avons vu notamment comment la réduction de la représentation analytique des substitutions a supporté des transferts d'analogies opératoires entre théorie arithmétique des formes, équations différentielles linéaires et groupes de substitutions. Nous avons vu aussi que Poincaré a envisagé cette pratique dans le cadre plus général de l' "idée de réduite" qu'il a attribué à Hermite. Cette notion joue un rôle clé dans ses travaux sans pour autant revêtir de définition précise dans un cadre disciplinaire donné. La nature algébrique des liens entre différents domaines n'est ainsi pas explicitée sous une forme structurelle ou abstraite comme nous le ferions aujourd'hui. Au contraire, les procédés algébriques ne sont jamais thématisés de manière autonomes mais toujours en articulation avec les significations qui peuvent leur être attachés dans différents cadres théoriques. Un exemple additionnel de tels transferts est donné par la manière dont Poincaré aborde le problème de la stabilité du problème des trois corps.\cite{Brechenmacher:2012c}\footnote{A propos de la conception de la stabilité chez Poincaré et de sa réception, voir \cite{Roque2011}.}  Ce dernier ajoute aux solutions périodiques des petites variations permettant de linéariser les équations et, par là, de réduire celles-ci à une forme canonique permettant d'engager une discussion de la stabilité du système sur le modèle de la conception de la stabilité des systèmes différentiels linéaires à coefficients constants qui avait été traditionnellement liée à la multiplicité des racines caractéristiques de l'équation séculaire.

Tout comme les situations liées aux fonctions fuchsiennes que nous avons rencontré dans cet article, l'exemple du problème des trois corps témoigne bien du rôle important joué par les pratiques algébriques de réduction canoniques élaborées par Poincaré bien que ces dernières ne jouent jamais le premier rôle dans une publication donnée. En raison des relations analogiques qu'elles supportent, ces dernières participent cependant à la fois de la définition des approches développées par Poincaré et des résultats visés par ce dernier.

Nous avons vu aussi que l'héritage des travaux de Jordan et Hermite a donné lieu à l'élaboration de procédés opératoires spécifiques. Ainsi, le calcul des Tableaux s'inscrit dans l'héritage des travaux d'Hermite sur les formes algébriques. Au début des années 1880, les travaux conjoints de Jordan et Poincaré sur ce sujet en ont enrichi la dimension opératoire afin de donner une représentation des étapes successives de leur pratique de réduction. La reprise de ces procédés par d'autres mathématiciens témoigne de ce que l'inscription des travaux de Jordan sur les groupes linéaires dans un cadre hermitien en a favorisé la circulation ultérieure.

\

La réévaluation des rôles et finalités de la note de 1884 que nous avons proposée dans cet article ouvre plus largement des questionnements relatifs à l'histoire de l'algèbre linéaire. Des difficultés similaires à celles que nous avons illustrées par le problème du rôle joué par les matrices dans la note de Poincaré se présentent en effet pour de nombreux autres travaux considérés comme fondateurs de notions algébriques tout en apparaissant isolés dans les \oe{}uvres de leurs auteurs. En outre, comme pour le cas du rôle attribué à Poincaré dans l'histoire des matrices, de nombreux travaux ont proposé une structuration en trois étapes de l'émergence de notions ou théories algébriques: à une première période caractérisée par la présence implicite de ces notions succède une étape d'explicitation d'objets étudiés pour eux-mêmes et de redécouvertes multiples puis une troisième période de théorisations rigoureuses. La note de 1884 nous invite pourtant à substituer différentes échelles de temporalités et des espaces non homogènes aux chronologies basées sur des découpages disciplinaires rétrospectifs. 

\nocite{Cartan:1898}

\bibliography{bibliographie}

\begin{thebibliography}{}
\expandafter\ifx\csname fonteauteurs\endcsname\relax
\def\fonteauteurs{\scshape}\fi
\expandafter\ifx\csname url\endcsname\relax
  \def\url#1{{\tt #1}}%
    \message{You should include the url package}\fi

\bibitem[Archibald, 2011]{Archibald:2011}
\bgroup\fonteauteurs\bgroup Archibald\egroup\egroup{}, T. (2011).
\newblock Differential equations and algebraic transcendents: French efforts at
  the creation of a galois theory of differential equations (1880-1910).
\newblock {\em Revue d'histoire des math{\'e}matiques},
  17(2)\string:\penalty500\relax 371--399.

\bibitem[Autonne, 1881]{Autonne-Poincar}
\bgroup\fonteauteurs\bgroup Autonne\egroup\egroup{}, L. (1881).
\newblock Correspondance d'{A}utonne (l{\'e}on) {\`a} {P}oincar{\'e} (henri).
\newblock Collection priv{\'e}e 75017. Reproduite {\`a} l'adresse
  http://www.univ-nancy2.fr/poincare/chp/.

\bibitem[Autonne, 1885]{Autonne:1885a}
\bgroup\fonteauteurs\bgroup Autonne\egroup\egroup{}, L. (1885).
\newblock Recherches sur les int{\'e}grales alg{\'e}briques des {\'e}quations
  diff{\'e}rentielles lin{\'e}aires {\`a} coefficients rationnels (second
  m{\'e}moire).
\newblock {\em Journal de l'{\'E}cole polytechnique},
  54\string:\penalty500\relax 1--30.

\bibitem[Boucard, 2011]{Boucard:2011}
\bgroup\fonteauteurs\bgroup Boucard\egroup\egroup{}, J. (2011).
\newblock Louis {P}oinsot et la th{\'e}orie de l'ordre : un cha{\^\i}non
  manquant entre {G}auss et {G}alois ?
\newblock {\em Revue d'histoire des math{\'e}matiques},
  17(1)\string:\penalty500\relax 41--138.

\bibitem[Brechenmacher, 2006a]{Brechenmacher:2006b}
\bgroup\fonteauteurs\bgroup Brechenmacher\egroup\egroup{}, F. (2006a).
\newblock A controversy and the writing of a history: the discussion of 'small
  oscillations' (1760-1860) from the standpoint of the controversy between
  {J}ordan and {K}ronecker (1874).
\newblock {\em Bulletin of the Belgian Mathematical Society},
  13\string:\penalty500\relax 941--944.

\bibitem[Brechenmacher, 2006b]{Brechenmacher:2006a}
\bgroup\fonteauteurs\bgroup Brechenmacher\egroup\egroup{}, F. (2006b).
\newblock {\em Histoire du th{\'e}or{\`e}me de Jordan de la d{\'e}composition
  matricielle (1870-1930). Formes de repr{\'e}sentations et m{\'e}thodes de
  d{\'e}compositions.}
\newblock Th\`ese de doctorat, Th{\`e}se de doctorat. Ecole des hautes
  {\'e}tudes en sciences sociales.

\bibitem[Brechenmacher, 2007a]{Brechenmacher:2007a}
\bgroup\fonteauteurs\bgroup Brechenmacher\egroup\egroup{}, F. (2007a).
\newblock La controverse de 1874 entre {C}amille {J}ordan et {L}eopold
  {K}ronecker.
\newblock {\em Revue d'Histoire des Math{\'e}matiques},
  13\string:\penalty500\relax 187--257.

\bibitem[Brechenmacher, 2007b]{Brechenmacher:2007b}
\bgroup\fonteauteurs\bgroup Brechenmacher\egroup\egroup{}, F. (2007b).
\newblock L'identit{\'e} alg{\'e}brique d'une pratique port{\'e}e par la
  discussion sur l'{\'e}quation {\`a} l'aide de laquelle on d{\'e}termine les
  in{\'e}galit{\'e}s s{\'e}culaires des plan{\`e}tes (1766-1874).
\newblock {\em Sciences et techniques en perspective},
  1\string:\penalty500\relax 5--85.

\bibitem[Brechenmacher, 2010]{Brechenmacher:2010a}
\bgroup\fonteauteurs\bgroup Brechenmacher\egroup\egroup{}, F. (2010).
\newblock Une histoire de l'universalit{\'e} des matrices math{\'e}matiques.
\newblock {\em Revue de Synth{\`e}se}, 4\string:\penalty500\relax 569--603.

\bibitem[Brechenmacher, 2011]{Brechenmacher:2011}
\bgroup\fonteauteurs\bgroup Brechenmacher\egroup\egroup{}, F. (2011).
\newblock Self-portraits with {{\'E}}variste {G}alois (and the shadow of
  {C}amille {J}ordan).
\newblock http://hal.archives-ouvertes.fr/aut/Frederic+Brechenmacher/.

\bibitem[Brechenmacher, 2012a]{Brechenmacher:2012c}
\bgroup\fonteauteurs\bgroup Brechenmacher\egroup\egroup{}, F. (2012a).
\newblock Linear groups in galois fields. a case study of tacit circulation of
  explicit knowledge.
\newblock {\em Oberwolfach Reports}, 4-2012\string:\penalty500\relax 48--54.

\bibitem[Brechenmacher, 2012b]{Brechenmacher:2011b}
\bgroup\fonteauteurs\bgroup Brechenmacher\egroup\egroup{}, F. (2012b).
\newblock On {J}ordan's measurements.
\newblock http://hal.archives-ouvertes.fr/aut/Frederic+Brechenmacher/.

\bibitem[Brechenmacher et Ehrhardt, 2010]{Brechenmacher:2010b}
\bgroup\fonteauteurs\bgroup Brechenmacher\egroup\egroup{}, F. et
  \bgroup\fonteauteurs\bgroup Ehrhardt\egroup\egroup{}, C. (2010).
\newblock On the identities of algebra in the 19th century.
\newblock {\em Oberwolfach Reports}, 7(1)\string:\penalty500\relax 24--31.

\bibitem[Brioschi, 1877]{Brioschi:1877}
\bgroup\fonteauteurs\bgroup Brioschi\egroup\egroup{}, F. (1877).
\newblock La th{\'e}orie des formes dans l'int{\'e}gration des {\'e}quations
  diff{\'e}rentielles lin{\'e}aires du second ordre.
\newblock {\em Mathematische Annalen}, 11\string:\penalty500\relax 401--412.

\bibitem[Briot et Bouquet, 1859]{BriotetBouquet:1859}
\bgroup\fonteauteurs\bgroup Briot\egroup\egroup{}, C. et
  \bgroup\fonteauteurs\bgroup Bouquet\egroup\egroup{}, C. (1859).
\newblock {\em Th{\'e}orie des fonctions doublement p{\'e}riodiques et, en
  particulier, des fonctions elliptiques}.
\newblock Paris.

\bibitem[Burkhardt \emph{et~al.}, 1916]{Burkhardt:1916}
\bgroup\fonteauteurs\bgroup Burkhardt\egroup\egroup{}, H.,
  \bgroup\fonteauteurs\bgroup Maurer\egroup\egroup{}, L. et
  \bgroup\fonteauteurs\bgroup Vessiot\egroup\egroup{}, E. (1916).
\newblock Groupes de transformations continus.
\newblock \emph{In} {\em Encyclop{\'e}die des sciences math{\'e}matiques},
  volume~4, pages 161--240. Gauthier-Villars, Paris.

\bibitem[Cahen et Vahlen, 1908]{Cahen:1908}
\bgroup\fonteauteurs\bgroup Cahen\egroup\egroup{}, E. et
  \bgroup\fonteauteurs\bgroup Vahlen\egroup\egroup{}, K.~T. (1908).
\newblock Th{\'e}orie arithm{\'e}tique des formes.
\newblock \emph{In} {\em Encyclop{\'e}die des sciences math{\'e}matiques},
  volume~3. Gauthier-Villars, Paris.

\bibitem[Cartan, 1894]{Cartan:1894}
\bgroup\fonteauteurs\bgroup Cartan\egroup\egroup{}, {\'E}. (1894).
\newblock Sur la r{\'e}duction de la structure d'un groupe {\`a} sa forme
  canonique.
\newblock {\em Comptes rendus hebdomadaires des s{\'e}ances de l'Acad{\'e}mie
  des sciences}, 119\string:\penalty500\relax 639--642.

\bibitem[Cartan, 1898]{Cartan:1898}
\bgroup\fonteauteurs\bgroup Cartan\egroup\egroup{}, {\'E}. (1898).
\newblock Les groupes bilin{\'e}aires et les syst{\`e}mes de nombres complexes.
\newblock {\em Annales de la Facult{\'e} des Sciences de Toulouse},
  12(2)\string:\penalty500\relax 65--99.

\bibitem[Cartan et Study, 1908]{Cartan:1908}
\bgroup\fonteauteurs\bgroup Cartan\egroup\egroup{}, {\'E}. et
  \bgroup\fonteauteurs\bgroup Study\egroup\egroup{}, E. (1908).
\newblock {\em Encyclop{\'e}die des sciences math{\'e}matiques}, volume~1,
  chapitre Nombres complexes, pages 329--468.
\newblock Gauthier-Villars, PAris.

\bibitem[Cifoletti, 1995]{Cifoletti:1995}
\bgroup\fonteauteurs\bgroup Cifoletti\egroup\egroup{}, G. (1995).
\newblock The creation of the history of algebra in the sixteenth century.
\newblock {\em in [Goldstein, Gray, Ritter 1995, p. 123-144]}.

\bibitem[Craig, 1889]{Craig:1889}
\bgroup\fonteauteurs\bgroup Craig\egroup\egroup{}, T. (1889).
\newblock {\em A treatise on linear differential equations. Vol. I. Equations
  with uniform coefficients}.
\newblock J. Wiley and Sons.

\bibitem[Darboux, 1889]{Darboux1889}
\bgroup\fonteauteurs\bgroup Darboux\egroup\egroup{}, G. (1889).
\newblock {\em Le{\c c}ons sur la th{\'e}orie g{\'e}n{\'e}rale des surfaces et
  les applications g{\'e}om{\'e}triques du calcul infinit{\'e}simal}.
\newblock Gauthier-Villars, Paris.

\bibitem[Dedekind, 1877]{Dedekind:1877}
\bgroup\fonteauteurs\bgroup Dedekind\egroup\egroup{}, R. (1877).
\newblock Schreiben an herrn borchardt {\"u}ber die theorie der elliptischen
  modul-functionen.
\newblock {\em Journal f{\"u}r die reine und angewandte Mathematik},
  83\string:\penalty500\relax 265--292.

\bibitem[Dieudonn{\'e}, 1962]{Dieudonne:1962}
\bgroup\fonteauteurs\bgroup Dieudonn{\'e}\egroup\egroup{}, J. (1962).
\newblock {\em Notes sur les travaux de Camille Jordan relatifs {\`a}
  l'alg{\`e}bre lin{\'e}aire et multilin{\'e}aire et la th{\'e}orie des
  nombres}.
\newblock in [Jordan {\OE}uvres, 3, p. V-XX].

\bibitem[Durand-Richard, ]{Durand-Richard:1996}
\bgroup\fonteauteurs\bgroup Durand-Richard\egroup\egroup{}, M.-J.
\newblock L'{\'e}cole alg{\'e}brique anglaise : les conditions conceptuelles et
  institutionnelles d'un calcul symbolique comme fondement de la connaissance.
\newblock \emph{In} {\em [GOLDSTEIN, GRAY, RITTER 1996]}, pages 445--477.

\bibitem[Durand-Richard, 2008]{Durand-Richard:2008}
\bgroup\fonteauteurs\bgroup Durand-Richard\egroup\egroup{}, M.-J., \'editeur
  (2008).
\newblock {\em L'analogie dans la d{\'e}marche scientifique}.
\newblock l'Harmattan.

\bibitem[Floquet, 1879]{Floquet:1879}
\bgroup\fonteauteurs\bgroup Floquet\egroup\egroup{}, G. (1879).
\newblock Sur la th{\'e}orie des {\'e}quations diff{\'e}rentielles
  lin{\'e}aires.
\newblock {\em Annales de l'{\'E}cole normale sup{\'e}rieure},
  8(2)\string:\penalty500\relax 3--132.

\bibitem[Frobenius, 1877]{Frobenius:1877}
\bgroup\fonteauteurs\bgroup Frobenius\egroup\egroup{}, G. (1877).
\newblock Ueber linear substitutionen und bilineare formen.
\newblock {\em Journal f{\"u}r die reine und angewandte Mathematik},
  84\string:\penalty500\relax 1--63.

\bibitem[Frobenius, 1879]{Frobenius:1879}
\bgroup\fonteauteurs\bgroup Frobenius\egroup\egroup{}, G. (1879).
\newblock Theorie der bilinearen formen mit ganzen coefficienten.
\newblock {\em Journal f{\"u}r die reine und angewandte Mathematik},
  86\string:\penalty500\relax 147--208.

\bibitem[Fuchs, 1876]{Fuchs:1876c}
\bgroup\fonteauteurs\bgroup Fuchs\egroup\egroup{}, L. (1876).
\newblock Extrait d'une lettre adress{\'e}e {\`a} m. {H}ermite.
\newblock {\em Journal de math{\'e}matiques pures et appliqu{\'e}es},
  2(3)\string:\penalty500\relax 158--160.

\bibitem[Fuchs, 1877]{Fuchs:1877}
\bgroup\fonteauteurs\bgroup Fuchs\egroup\egroup{}, L. (1877).
\newblock Sur quelques propri{\'e}t{\'e}s des int{\'e}grales des {\'e}quations
  diff{\'e}rentielles auxquelles satisfont les modules de p{\'e}riodicit{\'e}
  des int{\'e}grales elliptiques des deux premi{\`e}res esp{\`e}ces.
\newblock {\em Journal f{\"u}r die reine und angewandte Mathematik},
  83\string:\penalty500\relax 13--38.

\bibitem[Fuchs, 1880]{Fuchs:1880}
\bgroup\fonteauteurs\bgroup Fuchs\egroup\egroup{}, L. (1880).
\newblock Ueber die functionen, welche durch umkehrung der integrale von
  l{\"o}sungen der linearen differentialgleichungen entstehen, g{\"o}ttingen
  nachrichten, 1880, p. 445-453.
\newblock {\em Bulletin des sciences math{\'e}matiques},
  4(2)\string:\penalty500\relax 328--336.

\bibitem[Gauthier, 2009]{Gauthier:2009}
\bgroup\fonteauteurs\bgroup Gauthier\egroup\egroup{}, S. (2009).
\newblock La g{\'e}om{\'e}trie dans la g{\'e}om{\'e}trie des nombres : histoire
  de discipline ou histoire de pratiques {\`a} partir des exemples de
  {M}inkowski, {M}ordell et {D}avenport.
\newblock {\em Revue d'histoire des math{\'e}matiques},
  15(2)\string:\penalty500\relax 215--242.

\bibitem[Gilain, 1991]{Gilain:1991}
\bgroup\fonteauteurs\bgroup Gilain\egroup\egroup{}, C. (1991).
\newblock La th{\'e}orie qualitative de poincar{\'e} et le probl{\`e}me de
  l'int{\'e}gration des {\'e}quations diff{\'e}rentielles.
\newblock \emph{In} {\em [GISPERT 1991]}, pages 215--242.

\bibitem[Goldstein, 1995]{Goldstein:1995}
\bgroup\fonteauteurs\bgroup Goldstein\egroup\egroup{}, C. (1995).
\newblock {\em Un th{\'e}or{\`e}me de Fermat et ses lecteurs}.
\newblock PUV, Saint-Denis.

\bibitem[Goldstein, 1999]{Goldstein:1999}
\bgroup\fonteauteurs\bgroup Goldstein\egroup\egroup{}, C. (1999).
\newblock Sur la question des m{\'e}thodes quantitatives en histoire des
  math{\'e}matiques : le cas de la th{\'e}orie des nombres en france (1870-
  1914).
\newblock {\em Acta historiae rerum necnon technicarum},
  3\string:\penalty500\relax 187--214.

\bibitem[Goldstein, 2007]{Goldstein:2007}
\bgroup\fonteauteurs\bgroup Goldstein\egroup\egroup{}, C. (2007).
\newblock The hermitian form of reading the disquisitiones.
\newblock \emph{In} {\em [Goldstein, Schappacher, Schwermer, 2007}, pages
  377--410.

\bibitem[Goldstein, 2009]{Goldstein:2009}
\bgroup\fonteauteurs\bgroup Goldstein\egroup\egroup{}, C. (2009).
\newblock L'arithm{\'e}tique de {F}ermat dans le contexte de la correspondance
  de {M}ersenne : une approche micro-sociale.
\newblock {\em Annales de la Facult{\'e} des sciences de Toulouse},
  XVIII\string:\penalty500\relax 25--57.

\bibitem[Goldstein, 2011]{Goldstein:2011}
\bgroup\fonteauteurs\bgroup Goldstein\egroup\egroup{}, C. (2011).
\newblock Charles {H}ermite's {S}troll through the {G}alois fields.
\newblock {\em Revue d'histoire des math{\'e}matiques},
  17\string:\penalty500\relax 135--152.

\bibitem[Goldstein et Schappacher, 2007]{GoldsteinSchappa:2007c}
\bgroup\fonteauteurs\bgroup Goldstein\egroup\egroup{}, C. et
  \bgroup\fonteauteurs\bgroup Schappacher\egroup\egroup{}, N. (2007).
\newblock Several disciplines and a book (1860--1901).
\newblock \emph{In} {\em [Goldstein, Schappacher, Schwermer, 2007]}, pages
  67--104.

\bibitem[Gordan, 1877]{Gordan:1877}
\bgroup\fonteauteurs\bgroup Gordan\egroup\egroup{}, P. (1877).
\newblock {\"U}ber endliche gruppen linearer transformationen einen
  veranderlichen.
\newblock {\em Mathematische Annalen}, 12\string:\penalty500\relax 23--46.

\bibitem[Gray, 2000]{Gray:2000}
\bgroup\fonteauteurs\bgroup Gray\egroup\egroup{}, J. (2000).
\newblock {\em Linear differential equations and group theory from Riemann to
  Poincar{\'e}}, volume 2nd ed.
\newblock Gray, Jeremy, Boston.

\bibitem[Hamburger, 1873]{Hamburger1873}
\bgroup\fonteauteurs\bgroup Hamburger\egroup\egroup{}, M. (1873).
\newblock Bemerkung {\"u}ber die form der integrale der linearen
  differentialgleichungen mit ver{\"a}nderlicher coefficienten.
\newblock {\em Journal f{\"u}r die reine und angewandte Mathematik},
  76\string:\penalty500\relax 113--125.

\bibitem[Hawkins, 1972]{Hawkins1972}
\bgroup\fonteauteurs\bgroup Hawkins\egroup\egroup{}, T. (1972).
\newblock Hypercomplex numbers, lie groups, and the creation of group
  representation theory.
\newblock {\em Archive for History of Exact Sciences,},
  8\string:\penalty500\relax 243--87.

\bibitem[Hawkins, 2000]{Hawkins2000}
\bgroup\fonteauteurs\bgroup Hawkins\egroup\egroup{}, T. (2000).
\newblock {\em Emergence of the theory of Lie groups an essay in the history of
  mathematics, 1869-1926}.
\newblock Springer, New York, Berlin, Barcelone.

\bibitem[Hermite, 1854]{Hermite1854}
\bgroup\fonteauteurs\bgroup Hermite\egroup\egroup{}, C. (1854).
\newblock Sur la th{\'e}orie des formes quadratiques.
\newblock {\em Journal f{\"u}r die reine und angewandte Mathematik}, 47.

\bibitem[Jordan, 1860]{Jordan1860}
\bgroup\fonteauteurs\bgroup Jordan\egroup\egroup{}, C. (1860).
\newblock {\em Sur le nombre des valeurs des fonctions. Th{\`e}ses
  pr{\'e}sent{\'e}es {\`a} la Facult{\'e} des sciences de {P}aris par {C}amille
  {J}ordan, 1re th{\`e}se}.
\newblock Mallet-Bachelier, Paris.

\bibitem[Jordan, 1868a]{Jordan1868b}
\bgroup\fonteauteurs\bgroup Jordan\egroup\egroup{}, C. (1868a).
\newblock M{\'e}moire sur les groupes des mouvements.
\newblock {\em Annali di Matematica}, 2\string:\penalty500\relax 167--215,
  322--345.

\bibitem[Jordan, 1868b]{Jordan1868a}
\bgroup\fonteauteurs\bgroup Jordan\egroup\egroup{}, C. (1868b).
\newblock Sur la r{\'e}solution alg{\'e}brique des {\'e}quations primitives de
  degr{\'e} $p^2$.
\newblock {\em Journal de math{\'e}matiques pures et appliqu{\'e}es},
  32(2)\string:\penalty500\relax 111--135.

\bibitem[Jordan, 1870]{Jordan1870}
\bgroup\fonteauteurs\bgroup Jordan\egroup\egroup{}, C. (1870).
\newblock {\em Trait{\'e} des substitutions et des {\'e}quations
  alg{\'e}briques}.
\newblock Gauthier-Villars, Paris.

\bibitem[Jordan, 1874]{Jordan1874a}
\bgroup\fonteauteurs\bgroup Jordan\egroup\egroup{}, C. (1874).
\newblock Sur une application de la th{\'e}orie des substitutions aux
  {\'e}quations diff{\'e}rentielles lin{\'e}aires.
\newblock {\em Comptes rendus hebdomadaires des s{\'e}ances de l'Acad{\'e}mie
  des sciences}, 78\string:\penalty500\relax 614--617.

\bibitem[Jordan, 1876a]{Jordan1876b}
\bgroup\fonteauteurs\bgroup Jordan\egroup\egroup{}, C. (1876a).
\newblock Sur la d{\'e}termination des groupes form{\'e}s d'un nombre fini de
  substitutions lin{\'e}aires.
\newblock {\em Comptes rendus hebdomadaires des s{\'e}ances de l'Acad{\'e}mie
  des sciences}, 83\string:\penalty500\relax 1035--1037.

\bibitem[Jordan, 1876b]{Jordan1876a}
\bgroup\fonteauteurs\bgroup Jordan\egroup\egroup{}, C. (1876b).
\newblock Sur les {\'e}quations lin{\'e}aires du second ordre dont les
  int{\'e}grales sont alg{\'e}briques.
\newblock {\em Comptes rendus hebdomadaires des s{\'e}ances de l'Acad{\'e}mie
  des sciences}, 82\string:\penalty500\relax 605--607.

\bibitem[Jordan, 1877]{Jordan1877a}
\bgroup\fonteauteurs\bgroup Jordan\egroup\egroup{}, C. (1877).
\newblock D{\'e}termination des groupes form{\'e}s d'un nombre fini de
  substitutions lin{\'e}aires.
\newblock {\em Comptes rendus hebdomadaires des s{\'e}ances de l'Acad{\'e}mie
  des sciences}, 84\string:\penalty500\relax 1446--1448.

\bibitem[Jordan, 1878]{Jordan1878}
\bgroup\fonteauteurs\bgroup Jordan\egroup\egroup{}, C. (1878).
\newblock M{\'e}moire sur les {\'e}quations diff{\'e}rentielles lin{\'e}aires
  {\`a} int{\'e}grale alg{\'e}brique.
\newblock {\em Journal f{\"u}r die reine und angewandte Mathematik},
  84\string:\penalty500\relax 89--215.

\bibitem[Jordan, 1879]{Jordan1879b}
\bgroup\fonteauteurs\bgroup Jordan\egroup\egroup{}, C. (1879).
\newblock Sur la d{\'e}termination des groupes d'ordre fini contenus dans le
  groupe lin{\'e}aire.
\newblock {\em Atti della Accademia di Napoli}, 8(11).

\bibitem[Jordan, 1880a]{Jordan-Poincar}
\bgroup\fonteauteurs\bgroup Jordan\egroup\egroup{}, C. (1880a).
\newblock Correspondance de {J}ordan (camille) {\`a} {P}oincar{\'e} (henri).
\newblock Collection priv{\'e}e 75017. Reproduite {\`a} l'adresse
  http://www.univ-nancy2.fr/poincare/chp/.

\bibitem[Jordan, 1880b]{Jordan1880c}
\bgroup\fonteauteurs\bgroup Jordan\egroup\egroup{}, C. (1880b).
\newblock M{\'e}moire sur l'{\'e}quivalence des formes.
\newblock {\em Journal de l'{\'E}cole polytechnique},
  28\string:\penalty500\relax 111--150.

\bibitem[Jordan, 1880c]{Jordan1880a}
\bgroup\fonteauteurs\bgroup Jordan\egroup\egroup{}, C. (1880c).
\newblock Sur la r{\'e}duction des substitutions lin{\'e}aires.
\newblock {\em Comptes rendus hebdomadaires des s{\'e}ances de l'Acad{\'e}mie
  des sciences}, 90\string:\penalty500\relax 598--601.

\bibitem[Jordan, 1880d]{Jordan1880b}
\bgroup\fonteauteurs\bgroup Jordan\egroup\egroup{}, C. (1880d).
\newblock Sur l'{\'e}quivalence des formes.
\newblock {\em Comptes rendus hebdomadaires des s{\'e}ances de l'Acad{\'e}mie
  des sciences}, 90\string:\penalty500\relax 1422--1423.

\bibitem[Jordan, 1881a]{Jordan1881a}
\bgroup\fonteauteurs\bgroup Jordan\egroup\egroup{}, C. (1881a).
\newblock Sur la r{\'e}duction des formes quadratiques.
\newblock {\em Comptes rendus hebdomadaires des s{\'e}ances de l'Acad{\'e}mie
  des sciences}, 93\string:\penalty500\relax 113--117.

\bibitem[Jordan, 1881b]{Jordan1881c}
\bgroup\fonteauteurs\bgroup Jordan\egroup\egroup{}, C. (1881b).
\newblock Sur la repr{\'e}sentation d'un nombre ou d'une forme quadratique par
  une autre forme quadratique.
\newblock {\em Comptes rendus hebdomadaires des s{\'e}ances de l'Acad{\'e}mie
  des sciences}, 93\string:\penalty500\relax 234--237.

\bibitem[Jordan, 1881c]{Jordan1881b}
\bgroup\fonteauteurs\bgroup Jordan\egroup\egroup{}, C. (1881c).
\newblock Sur l'{\'e}quivalence des formes quadratiques.
\newblock {\em Comptes rendus hebdomadaires des s{\'e}ances de l'Acad{\'e}mie
  des sciences}, 93\string:\penalty500\relax 181--185.

\bibitem[Jordan, 1882a]{Jordan1882b}
\bgroup\fonteauteurs\bgroup Jordan\egroup\egroup{}, C. (1882a).
\newblock {\em Cours d'Analyse de l'Ecole Polytechnique}.
\newblock Gauthier-Villars, Paris.

\bibitem[Jordan, 1882b]{Jordan1882a}
\bgroup\fonteauteurs\bgroup Jordan\egroup\egroup{}, C. (1882b).
\newblock Sur la th{\'e}orie arithm{\'e}tique des formes quadratiques.
\newblock {\em Journal de l'{\'E}cole polytechnique},
  31\string:\penalty500\relax 1--43.

\bibitem[Klein, 1875]{Klein1875}
\bgroup\fonteauteurs\bgroup Klein\egroup\egroup{}, F. (1875).
\newblock Ueber bin{\"a}re formen mit linearen transformationen in sich selbst.
\newblock {\em Mathematische Annalen}, 9\string:\penalty500\relax 183--208.

\bibitem[Klein, 1877]{Klein1877}
\bgroup\fonteauteurs\bgroup Klein\egroup\egroup{}, F. (1877).
\newblock Ueber lineare differentialgleichungen.
\newblock {\em Mathematische Annalen}, 11, 12\string:\penalty500\relax 115--119
  , 167--180.

\bibitem[Laurent, 1890]{Laurent1890}
\bgroup\fonteauteurs\bgroup Laurent\egroup\egroup{}, H. (1890).
\newblock {\em Trait{\'e} d'analyse. Tome V. Calcul int{\'e}gral. {\'E}quations
  diff{\'e}rentielles ordinaires}.
\newblock Gauthier-Villars, Paris.

\bibitem[Nabonnand, 2000]{Nabonnand2000}
\bgroup\fonteauteurs\bgroup Nabonnand\egroup\egroup{}, P. (2000).
\newblock Les recherches sur l'oeuvre de poincar{\'e}.
\newblock {\em La gazette des math{\'e}maticiens}, 85\string:\penalty500\relax
  34--54.

\bibitem[Neumann, 2007]{Neumann2007}
\bgroup\fonteauteurs\bgroup Neumann\egroup\egroup{}, O. (2007).
\newblock {\em The Disquisitiones Arithmeticae and the Theory of Equations}.
\newblock in [Goldstein, Schappacher, Schwermer, 2007].

\bibitem[Neumann, 2006]{Neumann2006}
\bgroup\fonteauteurs\bgroup Neumann\egroup\egroup{}, P.~M. (2006).
\newblock The concept of primitivity in group theory and the second memoir of
  galois.
\newblock {\em Archive for History of Exact Sciences},
  60\string:\penalty500\relax 379--429.

\bibitem[Parshall, 1985]{Parshall1985}
\bgroup\fonteauteurs\bgroup Parshall\egroup\egroup{}, K.~H. (1985).
\newblock J. {H}.{M}. {W}edderburn and the structure theory of algebras.
\newblock {\em Archive for History of Exact Sciences},
  32\string:\penalty500\relax 223--349.

\bibitem[Petri et Norbert, 2004]{Petri2004}
\bgroup\fonteauteurs\bgroup Petri\egroup\egroup{}, B. et
  \bgroup\fonteauteurs\bgroup Norbert\egroup\egroup{}, S. (2004).
\newblock {\em From Abel to Kronecker. Episodes from 19th Century Algebra},
  pages 227--266.
\newblock in LAUDAL (Olav Arnfinn), PIENE (Ragni), The Legacy of Niels Henrik
  Abel ; the Abel Bicentennial, Oslo 2002, Springer Verlag, Berlin, 2004.

\bibitem[Picard, 1883]{Picard1884}
\bgroup\fonteauteurs\bgroup Picard\egroup\egroup{}, {\'E}. (1883).
\newblock Sur certaines substitutions lin{\'e}aires.
\newblock {\em Comptes rendus hebdomadaires des s{\'e}ances de l'Acad{\'e}mie
  des sciences}, 48\string:\penalty500\relax 416--417.

\bibitem[Poincar{\'e}, 1879a]{Poincar1879b}
\bgroup\fonteauteurs\bgroup Poincar{\'e}\egroup\egroup{}, H. (1879a).
\newblock Sur les formes quadratiques.
\newblock {\em Comptes rendus hebdomadaires des s{\'e}ances de l'Acad{\'e}mie
  des sciences}, 89\string:\penalty500\relax 897--899.

\bibitem[Poincar{\'e}, 1879b]{Poincar1879a}
\bgroup\fonteauteurs\bgroup Poincar{\'e}\egroup\egroup{}, H. (1879b).
\newblock Sur quelques propri{\'e}t{\'e}s des formes quadratiques.
\newblock {\em Comptes rendus hebdomadaires des s{\'e}ances de l'Acad{\'e}mie
  des sciences}, 89\string:\penalty500\relax 897--899.

\bibitem[Poincar{\'e}, 1880a]{Poincar1880d}
\bgroup\fonteauteurs\bgroup Poincar{\'e}\egroup\egroup{}, H. (1880a).
\newblock Sur la r{\'e}duction simultan{\'e}e d'une forme quadratique et d'une
  forme lin{\'e}aire.
\newblock {\em Comptes rendus hebdomadaires des s{\'e}ances de l'Acad{\'e}mie
  des sciences}, 91\string:\penalty500\relax 844--846.

\bibitem[Poincar{\'e}, 1880b]{Poincar1880b}
\bgroup\fonteauteurs\bgroup Poincar{\'e}\egroup\egroup{}, H. (1880b).
\newblock Sur les courbes d{\'e}finies par les {\'e}quations
  diff{\'e}rentielles.
\newblock {\em Comptes rendus hebdomadaires des s{\'e}ances de l'Acad{\'e}mie
  des sciences}, 90\string:\penalty500\relax 673--675.

\bibitem[Poincar{\'e}, 1880c]{Poincar1880c}
\bgroup\fonteauteurs\bgroup Poincar{\'e}\egroup\egroup{}, H. (1880c).
\newblock Sur les formes cubiques ternaires.
\newblock {\em Comptes rendus hebdomadaires des s{\'e}ances de l'Acad{\'e}mie
  des sciences}, 90\string:\penalty500\relax 1336--1339.

\bibitem[Poincar{\'e}, 1880d]{Poincar1880a}
\bgroup\fonteauteurs\bgroup Poincar{\'e}\egroup\egroup{}, H. (1880d).
\newblock Sur un mode nouveau de repr{\'e}sentation g{\'e}om{\'e}trique des
  formes quadratiques d{\'e}finies ou ind{\'e}finies.
\newblock {\em Journal de l'{\'E}cole polytechnique},
  47\string:\penalty500\relax 177--245.

\bibitem[Poincar{\'e}, 1881a]{Poincar1881d}
\bgroup\fonteauteurs\bgroup Poincar{\'e}\egroup\egroup{}, H. (1881a).
\newblock Sur la repr{\'e}sentation des nombres par les formes.
\newblock {\em Comptes rendus hebdomadaires des s{\'e}ances de l'Acad{\'e}mie
  des sciences}, 92\string:\penalty500\relax 333--335.

\bibitem[Poincar{\'e}, 1881b]{Poincar1881c}
\bgroup\fonteauteurs\bgroup Poincar{\'e}\egroup\egroup{}, H. (1881b).
\newblock Sur les {\'e}quations diff{\'e}rentielles lin{\'e}aires {\`a}
  int{\'e}grales alg{\'e}briques.
\newblock {\em Comptes rendus hebdomadaires des s{\'e}ances de l'Acad{\'e}mie
  des sciences}, 92\string:\penalty500\relax 698--701.

\bibitem[Poincar{\'e}, 1881c]{Poincar1881b}
\bgroup\fonteauteurs\bgroup Poincar{\'e}\egroup\egroup{}, H. (1881c).
\newblock Sur les fonctions fuchsiennes.
\newblock {\em Comptes rendus hebdomadaires des s{\'e}ances de l'Acad{\'e}mie
  des sciences}, 92\string:\penalty500\relax 333--335.

\bibitem[Poincar{\'e}, 1881d]{Poincar1881a}
\bgroup\fonteauteurs\bgroup Poincar{\'e}\egroup\egroup{}, H. (1881d).
\newblock Sur les formes cubiques ternaires et quaternaires.
\newblock {\em Journal de l'{\'E}cole polytechnique},
  50\string:\penalty500\relax 199--253.

\bibitem[Poincar{\'e}, 1882a]{Poincar1882d}
\bgroup\fonteauteurs\bgroup Poincar{\'e}\egroup\egroup{}, H. (1882a).
\newblock Sur les fonctions uniformes qui se reproduisent par des substitutions
  lin{\'e}aires.
\newblock {\em Mathematische Annalen}, 19\string:\penalty500\relax 553--564.

\bibitem[Poincar{\'e}, 1882b]{Poincar1882a}
\bgroup\fonteauteurs\bgroup Poincar{\'e}\egroup\egroup{}, H. (1882b).
\newblock Sur les formes cubiques ternaires et quaternaires. deuxi{\`e}me
  partie.
\newblock {\em Journal de l'{\'E}cole polytechnique},
  51\string:\penalty500\relax 45--91.

\bibitem[Poincar{\'e}, 1882c]{Poincar1882c}
\bgroup\fonteauteurs\bgroup Poincar{\'e}\egroup\egroup{}, H. (1882c).
\newblock Sur les groupes discontinus.
\newblock {\em Comptes rendus hebdomadaires des s{\'e}ances de l'Acad{\'e}mie
  des sciences}, 94\string:\penalty500\relax 840--843.

\bibitem[Poincar{\'e}, 1882d]{Poincar1882e}
\bgroup\fonteauteurs\bgroup Poincar{\'e}\egroup\egroup{}, H. (1882d).
\newblock Th{\'e}orie des groupes fuchsiens,.
\newblock {\em Acta Mathematica}, 1\string:\penalty500\relax 1--62.

\bibitem[Poincar{\'e}, 1883a]{Poincar1883b}
\bgroup\fonteauteurs\bgroup Poincar{\'e}\egroup\egroup{}, H. (1883a).
\newblock Sur la reproduction des formes.
\newblock {\em Comptes rendus hebdomadaires des s{\'e}ances de l'Acad{\'e}mie
  des sciences}, 97\string:\penalty500\relax 949--951.

\bibitem[Poincar{\'e}, 1883b]{Poincar1883a}
\bgroup\fonteauteurs\bgroup Poincar{\'e}\egroup\egroup{}, H. (1883b).
\newblock Sur les groupes des {\'e}quations lin{\'e}aires.
\newblock {\em Comptes rendus hebdomadaires des s{\'e}ances de l'Acad{\'e}mie
  des sciences}, 96\string:\penalty500\relax 691--694.

\bibitem[Poincar{\'e}, 1884a]{Poincar1884g}
\bgroup\fonteauteurs\bgroup Poincar{\'e}\egroup\egroup{}, H. (1884a).
\newblock Sur la r{\'e}duction des p{\'e}riodes des int{\'e}grales
  ab{\'e}liennes.
\newblock {\em Bulletin de la soci{\'e}t{\'e} math{\'e}matique de France},
  12\string:\penalty500\relax 124--143.

\bibitem[Poincar{\'e}, 1884b]{Poincar1884e}
\bgroup\fonteauteurs\bgroup Poincar{\'e}\egroup\egroup{}, H. (1884b).
\newblock Sur les groupes des {\'e}quations lin{\'e}aires.
\newblock {\em Acta Mathematica}, 4\string:\penalty500\relax 202--311.

\bibitem[Poincar{\'e}, 1884c]{Poincar1884b}
\bgroup\fonteauteurs\bgroup Poincar{\'e}\egroup\egroup{}, H. (1884c).
\newblock Sur les groupes hyperfuchsiens.
\newblock {\em Comptes rendus hebdomadaires des s{\'e}ances de l'Acad{\'e}mie
  des sciences}, 98\string:\penalty500\relax 503--504.

\bibitem[Poincar{\'e}, 1884d]{Poincar1884d}
\bgroup\fonteauteurs\bgroup Poincar{\'e}\egroup\egroup{}, H. (1884d).
\newblock Sur les nombres complexes.
\newblock {\em Comptes rendus hebdomadaires des s{\'e}ances de l'Acad{\'e}mie
  des sciences}, 99\string:\penalty500\relax 740--742.

\bibitem[Poincar{\'e}, 1884e]{Poincar1884a}
\bgroup\fonteauteurs\bgroup Poincar{\'e}\egroup\egroup{}, H. (1884e).
\newblock Sur les substitutions lin{\'e}aires.
\newblock {\em Comptes rendus hebdomadaires des s{\'e}ances de l'Acad{\'e}mie
  des sciences}, 98\string:\penalty500\relax 349--352.

\bibitem[Poincar{\'e}, 1884f]{Poincar1884f}
\bgroup\fonteauteurs\bgroup Poincar{\'e}\egroup\egroup{}, H. (1884f).
\newblock Sur une g{\'e}n{\'e}ralisation des fractions continues.
\newblock {\em Comptes rendus hebdomadaires des s{\'e}ances de l'Acad{\'e}mie
  des sciences}, 99\string:\penalty500\relax 1014--1016.

\bibitem[Poincar{\'e}, 1886a]{Poincar1886b}
\bgroup\fonteauteurs\bgroup Poincar{\'e}\egroup\egroup{}, H. (1886a).
\newblock R{\'e}duction simultan{\'e}e d'une forme quadratique et d'une forme
  lin{\'e}aire.
\newblock {\em Journal de l'{\'E}cole polytechnique},
  56\string:\penalty500\relax 79--142.

\bibitem[Poincar{\'e}, 1886b]{Poincar1886e}
\bgroup\fonteauteurs\bgroup Poincar{\'e}\egroup\egroup{}, H. (1886b).
\newblock Sur la r{\'e}duction des p{\'e}riodes des int{\'e}grales
  ab{\'e}liennes.
\newblock {\em Comptes rendus hebdomadaires des s{\'e}ances de l'Acad{\'e}mie
  des sciences}, 102\string:\penalty500\relax 915--916.

\bibitem[Poincar{\'e}, 1886c]{Poincar1886a}
\bgroup\fonteauteurs\bgroup Poincar{\'e}\egroup\egroup{}, H. (1886c).
\newblock Sur la repr{\'e}sentation des nombres par les formes.
\newblock {\em Bulletin de la soci{\'e}t{\'e} math{\'e}matique de France},
  18\string:\penalty500\relax 162--194.

\bibitem[Poincar{\'e}, 1886d]{Poincar1886c}
\bgroup\fonteauteurs\bgroup Poincar{\'e}\egroup\egroup{}, H. (1886d).
\newblock Sur les d{\'e}terminants d'ordre infini.
\newblock {\em Bulletin de la soci{\'e}t{\'e} math{\'e}matique de France},
  14\string:\penalty500\relax 77--90.

\bibitem[Poincar{\'e}, 1886e]{Poincar1886f}
\bgroup\fonteauteurs\bgroup Poincar{\'e}\egroup\egroup{}, H. (1886e).
\newblock Sur les fonctions ab{\'e}liennes.
\newblock {\em American Journal of mathematics}, 8\string:\penalty500\relax
  283--342.

\bibitem[Poincar{\'e}, 1886f]{Poincar1886d}
\bgroup\fonteauteurs\bgroup Poincar{\'e}\egroup\egroup{}, H. (1886f).
\newblock Sur les fonctions fuchsiennes et les formes quadratiques
  ind{\'e}finies.
\newblock {\em Comptes rendus hebdomadaires des s{\'e}ances de l'Acad{\'e}mie
  des sciences}, 102\string:\penalty500\relax 735--737.

\bibitem[Poincar{\'e}, 1886g]{Poincar1886g}
\bgroup\fonteauteurs\bgroup Poincar{\'e}\egroup\egroup{}, H. (1886g).
\newblock Sur les int{\'e}grales irr{\'e}guli{\`e}res des {\'e}quations
  lin{\'e}aires.
\newblock {\em Acta Math{\'e}matica}, 8\string:\penalty500\relax 295--344.

\bibitem[Poincar{\'e}, 1887]{Poincar1887}
\bgroup\fonteauteurs\bgroup Poincar{\'e}\egroup\egroup{}, H. (1887).
\newblock Les fonctions fuchsiennes et l'arithm{\'e}tique.
\newblock {\em Journal de math{\'e}matiques pures et appliqu{\'e}es}, (4)
  3\string:\penalty500\relax 405--464.

\bibitem[Poincar{\'e}, 1890]{Poincar1890}
\bgroup\fonteauteurs\bgroup Poincar{\'e}\egroup\egroup{}, H. (1890).
\newblock Sur le probl{\`e}me des trois corps et les {\'e}quations de la
  dynamique.
\newblock {\em Acta Mathematica}, 13\string:\penalty500\relax 1--270.

\bibitem[Poincar{\'e}, 1899]{Poincar1899}
\bgroup\fonteauteurs\bgroup Poincar{\'e}\egroup\egroup{}, H. (1899).
\newblock Compl{\'e}ment {\`a} l'analysis situs.
\newblock {\em Rendiconti del circolo matematico del Palermo},
  13\string:\penalty500\relax 285--343.

\bibitem[Poincar{\'e}, 1900]{Poincar1900}
\bgroup\fonteauteurs\bgroup Poincar{\'e}\egroup\egroup{}, H. (1900).
\newblock Second compl{\'e}ment {\`a} l'analysis situs.
\newblock {\em Proceedings of the London Mathematical Society},
  32\string:\penalty500\relax 277--308.

\bibitem[Poincar{\'e}, 1901]{Poincar1901}
\bgroup\fonteauteurs\bgroup Poincar{\'e}\egroup\egroup{}, H. (1901).
\newblock Quelques remarques sur les groupes continus.
\newblock {\em Rendiconti del circolo matematico del Palermo},
  15\string:\penalty500\relax 321--368.

\bibitem[Poincar{\'e}, 1908a]{Poincar1908b}
\bgroup\fonteauteurs\bgroup Poincar{\'e}\egroup\egroup{}, H. (1908a).
\newblock {\em L'avenir math{\'e}matique}, pages 19--42.
\newblock in Science et M{\'e}thode, Flammarion : Paris, 1947.

\bibitem[Poincar{\'e}, 1908b]{Poincar1908a}
\bgroup\fonteauteurs\bgroup Poincar{\'e}\egroup\egroup{}, H. (1908b).
\newblock {\em L'invention math{\'e}matique, conf{\'e}rence fait {\`a}
  l'Institut g{\'e}n{\'e}ral psychologique}, pages 43--63.
\newblock in Science et M{\'e}thode, Flammarion : Paris, 1947.

\bibitem[Poincar{\'e}, 1921]{Poincar1921}
\bgroup\fonteauteurs\bgroup Poincar{\'e}\egroup\egroup{}, H. (1921).
\newblock Analyse des travaux scientifiques de henri poincar{\'e} faite par
  lui-m{\^e}me.
\newblock {\em Acta Mathematica}, 38\string:\penalty500\relax 1--135.

\bibitem[Poincar{\'e}, 1997]{GrayWalter:1997}
\bgroup\fonteauteurs\bgroup Poincar{\'e}\egroup\egroup{}, H. (1997).
\newblock {\em Trois suppl{\'e}ments sur la d{\'e}couverte des fonctions
  fuchsiennes}.
\newblock Gray, Jeremy and Walter, Scott, Paris / Berlin.

\bibitem[Robadey, 2004]{Robadey2004}
\bgroup\fonteauteurs\bgroup Robadey\egroup\egroup{}, A. (2004).
\newblock Exploration d'un mode d'{\'e}criture de la g{\'e}n{\'e}ralit{\'e}:
  l'article de poincar{\'e} sur les lignes g{\'e}od{\'e}siques des surfaces
  convexes (1905).
\newblock {\em Revue d'histoire des math{\'e}matiques},
  10(2)\string:\penalty500\relax 257--318.

\bibitem[Roque, 2011]{Roque2011}
\bgroup\fonteauteurs\bgroup Roque\egroup\egroup{}, T. (2011).
\newblock Stability of trajectories from poincar{\'e} to birkhoff : approaching
  a qualitative definition.
\newblock {\em Archive for History of Exact Sciences},
  65\string:\penalty500\relax 295--342.

\bibitem[Sauvage, 1895]{Sauvage1895}
\bgroup\fonteauteurs\bgroup Sauvage\egroup\egroup{}, L. (1895).
\newblock {\em Th{\'e}orie g{\'e}n{\'e}rale des syst{\`e}mes d'{\'e}quations
  diff{\'e}rentielles lin{\'e}aires et homog{\`e}nes}.
\newblock Gauthier-Villars, Paris.

\bibitem[Scheffers, 1891]{Scheffers1891}
\bgroup\fonteauteurs\bgroup Scheffers\egroup\egroup{}, G. (1891).
\newblock Zuruckf{\"u}hrung complexer zahlensysteme auf typische formen.
\newblock {\em Mathematische Annalen}, 39\string:\penalty500\relax 293--390.

\bibitem[Schmid, 1982]{Schmid1982}
\bgroup\fonteauteurs\bgroup Schmid\egroup\egroup{}, W. (1982).
\newblock Poincar{\'e} and {L}ie groups.
\newblock {\em Bulletin of the American Mathematical society},
  6-2\string:\penalty500\relax 175--186.

\bibitem[Schwarz, 1873]{Schwarz1872}
\bgroup\fonteauteurs\bgroup Schwarz\egroup\egroup{}, H. (1873).
\newblock Ueber diejenigen f{\"a}lle, in welchen die gaussische
  hypergeometrische reihe eine algebraische function ihres vierten elementes
  darstellt.
\newblock {\em Journal f{\"u}r die reine und angewandte Mathematik},
  75\string:\penalty500\relax 292--335.

\bibitem[Schwermer, 2007]{Schwermer2007}
\bgroup\fonteauteurs\bgroup Schwermer\egroup\egroup{} (2007).
\newblock {\em Reduction Theory of Quadratic Forms : Toward R{\"a}umliche
  Anschauung in Minkowski's Early Work}.
\newblock in [Golstein, Schappacher, Schwermer 2007], p. 483-505.

\bibitem[Serret, 1886]{Serret1886}
\bgroup\fonteauteurs\bgroup Serret\egroup\egroup{}, J.-A. (1886).
\newblock {\em Cours de calcul diff{\'e}rentiel et int{\'e}gral, 3e
  {\'e}dition}.
\newblock Gauthier-Villars.

\bibitem[Study, 1889]{Study1889}
\bgroup\fonteauteurs\bgroup Study\egroup\egroup{}, E. (1889).
\newblock Complexe zahlen und transformationsgruppen.
\newblock {\em Berichte {\"u}ber die Verhandlungen der Kgl. S{\"a}chsischen
  Gesellschaft der Wissenschaften zu Leipzi}, pages 177--227.

\bibitem[Sylvester, 1851]{Sylvester1851}
\bgroup\fonteauteurs\bgroup Sylvester\egroup\egroup{}, J.~J. (1851).
\newblock Enumeration of the contacts of lines and surfaces of the second
  order.
\newblock {\em Philosophical Magazine}, pages 119--140.

\bibitem[Sylvester, 1852]{Sylvester1852}
\bgroup\fonteauteurs\bgroup Sylvester\egroup\egroup{}, J.~J. (1852).
\newblock Sur une propri{\'e}t{\'e} nouvelle de l'{\'e}quation qui sert {\`a}
  d{\'e}terminer les in{\'e}galit{\'e}s s{\'e}culaires des plan{\`e}tes.
\newblock {\em Nouvelles Annales de Math{\'e}matiques}, pages 438--440.

\bibitem[Sylvester, 1882a]{Sylvester1882a}
\bgroup\fonteauteurs\bgroup Sylvester\egroup\egroup{}, J.~J. (1882a).
\newblock Sur les puissances et les racines de substitutions lin{\'e}aires.
\newblock {\em Comptes rendus hebdomadaires des s{\'e}ances de l'Acad{\'e}mie
  des sciences}, 94(55-59).

\bibitem[Sylvester, 1882b]{Sylvester1882b}
\bgroup\fonteauteurs\bgroup Sylvester\egroup\egroup{}, J.~J. (1882b).
\newblock Sur les racines des matrices unitaires.
\newblock {\em Comptes rendus hebdomadaires des s{\'e}ances de l'Acad{\'e}mie
  des sciences}, 94\string:\penalty500\relax 396--99.

\bibitem[Sylvester, 1883]{Sylvester1883}
\bgroup\fonteauteurs\bgroup Sylvester\egroup\egroup{}, J.~J. (1883).
\newblock On the equation to the secular inequalities in the planetary theory.
\newblock {\em Philosophical Magazine}, (16) 100\string:\penalty500\relax
  267--269.

\bibitem[Tannery, 1875]{Tannery1875}
\bgroup\fonteauteurs\bgroup Tannery\egroup\egroup{}, J. (1875).
\newblock Sur une {\'e}quation diff{\'e}rentielle lin{\'e}aire du second ordre.
\newblock {\em Annales scientifiques de l'{\'E}cole normale sup{\'e}rieure},
  (2) 8\string:\penalty500\relax 169--194.

\bibitem[Tazzioli, 1994]{Tazzioli1994}
\bgroup\fonteauteurs\bgroup Tazzioli\egroup\egroup{}, R. (1994).
\newblock Schwarz's critique and interpretation of the riemann representation
  theorem.
\newblock {\em Rendiconti del circolo matematico del Palermo}, (2)
  34\string:\penalty500\relax 95--132.

\bibitem[Ton-That et Tran, 1999]{Ton-That1999}
\bgroup\fonteauteurs\bgroup Ton-That\egroup\egroup{}, T. et
  \bgroup\fonteauteurs\bgroup Tran\egroup\egroup{}, T.-D. (1999).
\newblock Poincar{\'e}'s proof of the so-called {B}irkoff-{W}itt theorem.
\newblock {\em Revue d'histoire des math{\'e}matiques},
  5(2)\string:\penalty500\relax 249--284.

\bibitem[Wedderburn, 1934]{Wedderburn1934}
\bgroup\fonteauteurs\bgroup Wedderburn\egroup\egroup{}, J.~H. (1934).
\newblock {\em Lectures on Matrices}, volume~17.
\newblock American Mathematical Society Colloquium Publications, New York.

\bibitem[Weierstrass, 1868]{Weierstrass1868}
\bgroup\fonteauteurs\bgroup Weierstrass\egroup\egroup{}, K. (1868).
\newblock Zur theorie der quadratischen und bilinearen formen.
\newblock {\em Monatsberichte der K{\"o}niglich Preussischen Akademie der
  Wissenschaften zu Berlin}, pages 310--338.

\bibitem[WEYR, 1890]{Weyr1890}
\bgroup\fonteauteurs\bgroup WEYR\egroup\egroup{}, E. (1890).
\newblock Zur theorie der bilinearen formen.
\newblock {\em Monatshefte f{\"u}r Mathematik und Physik},
  1\string:\penalty500\relax 161--235.

\end{thebibliography}
\bibliographystyle{apalike-fr}

\end{document}